\documentclass[10pt,reqno, colorlinks]{amsart}

\usepackage{amsmath, amsfonts, amssymb, amscd, enumerate, stmaryrd}
\usepackage{hyperref}
\newtheorem{theo}{Theorem}[section]

\newtheorem{example}[theo]{Example}
\newtheorem{definition}[theo]{Definition}
\newenvironment{pf}{\noindent{\it Proof. }}{$\square$\par\medskip}

\newtheorem{lemma}[theo]{Lemma}
\newtheorem{theorem}[theo]{Theorem}
\newtheorem{corollary}[theo]{Corollary}
\newtheorem{proposition}[theo]{Proposition}

\newtheorem{remark}[theo]{Remark}

\newcommand{\beq}{\begin{equation}}
\newcommand{\eeq}{\end{equation}}

\renewcommand{\i}{\iota}

\renewcommand{\L}{\Lambda}

\newcommand{\bR}{\mathbb{R}}
\newcommand{\bZ}{\mathbb{Z}}

\newcommand{\bN}{\mathbb{N}}

\renewcommand{\gg}{\mathfrak{g}}
\newcommand{\gh}{\mathfrak{h}}
\newcommand{\gk}{\mathfrak{k}}

\newcommand{\gm}{\mathfrak{m}}

\newcommand{\gr}{\mathfrak{r}}

\newcommand{\so}{\mathfrak{so}}

\newcommand{\spin}{\mathfrak{spin}}

\newcommand{\ggl}{\mathfrak{gl}}

\newcommand\GL{\mathrm{GL}}
\newcommand\SL{\mathrm{SL}}
\newcommand\SO{\mathrm{SO}}

\newcommand\Spin{\mathrm{Spin}}

\renewcommand\sl{\mathfrak{sl}}

\newcommand{\cA}{\mathcal{A}}

\newcommand{\cC}{\mathcal{C}}

\newcommand{\cL}{\mathcal{L}}

\newcommand{\cP}{\mathcal{P}}

\newcommand{\cS}{\mathcal{S}}

\newcommand{\cU}{\mathcal{U}}

\def\Map{\mathop\text{\rm Map}\nolimits}

\renewcommand{\square}{\kern1pt\vbox
{\hrule height 0.6pt\hbox{\vrule width 0.6pt\hskip 3pt
\vbox{\vskip 6pt}\hskip 3pt\vrule width 0.6pt}\hrule height0.6pt}\kern1pt}

\DeclareMathOperator\ad{ad}

\DeclareMathOperator\CoDer{CoDer\;}
\DeclareMathOperator\End{End\;}
\DeclareMathOperator\Ker{Ker\;}
\DeclareMathOperator\Hom{Hom\;}
\DeclareMathOperator\Aut{Aut\;}
\DeclareMathOperator\Mor{Mor\;}

\DeclareMathOperator\Ad{Ad}

\DeclareMathOperator\Id{Id}
\DeclareMathOperator{\sgn}{sgn}

\DeclareMathOperator{\Der}{Der\;}

\DeclareMathOperator{\Span}{Span}
\DeclareMathOperator{\Conn}{Conn}
\DeclareMathOperator{\Bil}{Bil}
\DeclareMathOperator{\Nom}{Nom}
\DeclareMathOperator{\Ric}{Ric}

\renewcommand\Im{\operatorname{Im}}

\newcommand{\ev}{{\operatorname{ev}}}

\newcommand{\ol}{\overline}
\newcommand{\0}{\overline{0}}

\newcommand{\ou}{\overline{1}}

\newcommand{\be}{\begin{equation}}
\newcommand{\ee}{\end{equation}}

\def\<#1,#2>{\langle\,#1,\,#2\,\rangle}
\newcommand{\arr}{\begin{array}{rlll}}
\newcommand{\ea}{\end{array}}
\newcommand{\bea}{\begin{eqnarray}}
\newcommand{\eea}{\end{eqnarray}}
\newcommand{\bean}{\begin{eqnarray*}}
\newcommand{\eean}{\end{eqnarray*}}

\def\sideremark#1{\ifvmode\leavevmode\fi\vadjust{
\vbox to0pt{\hbox to 0pt{\hskip\hsize\hskip1em
\vbox{\hsize3cm\tiny\raggedright\pretolerance10000
\noindent #1\hfill}\hss}\vbox to8pt{\vfil}\vss}}}

\newcounter{ssig}
\newcounter{ttig}
\begin{document}
\title[Superization of homogeneous spin manifolds]{
Superization of homogeneous spin manifolds\\ and geometry of homogeneous supermanifolds}
\author[Andrea Santi]{Andrea \ Santi}
%
\maketitle
\setcounter{equation}{0}
\bigskip
\bigskip
{\bf Abstract:} 
Let $M_{0}=G_{0}/H$ be a (pseudo)-Riemannian homogeneous spin manifold, with reductive decomposition $\gg_{0}=\gh+\gm$ and let $S(M_{0})$ be the spin bundle defined by the spin representation $\tilde{\Ad}:H\rightarrow\GL_{\bR}(S)$ of the stabilizer $H$. 
This article studies the superizations of $M_{0}$, \textit{i.e.} its extensions to a homogeneous supermanifold $M=G/H$ whose sheaf of superfunctions is isomorphic to $\Lambda(S^{*}(M_{0}))$. Here $G$ is the Lie supergroup associated with a certain extension of the Lie algebra of symmetry $\gg_{0}$ to an algebra of supersymmetry $\gg=\gg_{\0}+\gg_{\ou}=\gg_{0}+S$ via the Kostant-Koszul construction. Each algebra of supersymmetry naturally determines a flat connection $\nabla^{\cS}$ in the spin bundle $S(M_{0})$. Killing vectors together with generalized Killing spinors (\textit{i.e.} $\nabla^{\cS}$-parallel spinors) are interpreted as the values of appropriate geometric symmetries of $M$, namely even and odd Killing fields. An explicit formula for the Killing representation of the algebra of supersymmetry is obtained, generalizing some results of Koszul. The generalized spin connection $\nabla^{\cS}$ defines a superconnection on $M$, via the super-version of a theorem of Wang. 
\tableofcontents
\section*{Introduction}
Lie superalgebras have played an important role in modern physics since the idea of supersymmetry arose. Complex and real simple Lie superalgebras were classified by Kac (\cite{Kac, Sc}). This classification was used to describe algebras of supersymmetry, \textit{i.e.} extensions of a Lie algebra of symmetry $\gg_{0}$ to a Lie superalgebra $\gg=\gg_{\0}+\gg_{\ou}=\gg_{0}+\gg_{\ou}$.
In this spirit, the classification of extensions of Poincare' Lie algebras, in all signatures and dimensions, was achieved in \cite{AC, ACDV}. 
\\
Lorentzian symmetric spin manifolds $M_{0}=G_{0}/H$ appear in constructions of maximally supersymmetric solutions of $11$-dimensional supergravity (\cite{CJS, F1, F2, FP, HKS}). In this context, some special, physically relevant algebras of supersymmetry $\gg=\gg_{\0}+\gg_{\ou}$ have been considered: the action of the even part $\gg_{\0}=\gg_{0}=\gh+\gm$ on the odd part $\gg_{\ou}=S$ is a spin representation, \textit{i.e.} an extension $\ad|_{\gg_{0}}:\gg_{0}\rightarrow\ggl_{\bR}(S)$ of the spin representation $\ad|_{\gh}:\gh\rightarrow\ggl_{\bR}(S)$ of the stability subalgebra $\gh$. 
In this paper, every such algebra of supersymmetry is called {\it adapted} to the spin manifold $(M_{0},g,S(M_{0}))$, where $S(M_{0})$ is the spin bundle
defined by the spin representation $\tilde{\Ad}:H\rightarrow\GL_{\bR}(S)$ of the stabilizer $H$. We remark that notions of generalized Killing spinors appear naturally in the context of supergravity theories (\cite{F}).
\\
"Even" adapted algebras of supersymmetry are used in \cite{Co} and \cite{F3}. \cite{Co} obtains a unified description of homogeneous quaternionic K\"ahler manifolds of solvable group by means of extended Poincare' Lie algebras (\cite{AC, ACDV}). This construction has a natural mirror in the setting of supergeometry and leads to the construction of homogeneous quaternionic K\"ahler supermanifolds, which is also described in \cite{Co}. 
\cite{F3} constructs compact real forms of the exceptional Lie algebras $F\sb 4$ and $E\sb 8$ in terms of even adapted algebras of supersymmetry of the spheres $S^{8}$ and $S^{15}$.
\\

This paper deals with a geometric representation of the algebra of supersymmetry as an algebra of Killing fields on a supermanifold. Recall that the spin bundle $S(M_{0})$ of a (pseudo)-Riemannian spin manifold $(M_{0},g,S(M_{0}))$ canonically defines a supermanifold $M=(M_{0},\mathcal{A}_{M})$ whose sheaf of superfunctions $\mathcal{A}_{M}$ is isomorphic to the sheaf of sections of the exterior algebra $\Lambda(S^{*}(M_{0}))$ of $S^{*}(M_{0})$. Supermanifolds of this type have been studied in \cite{ACDS, K1, K2}. If $M_{0}=G_{0}/H$ is a homogeneous reductive (pseudo)-Riemannian spin manifold 
then every adapted algebra of supersymmetry $\gg=\gg_{\0}+\gg_{\ou}=\gg_{0}+S$ 
defines a structure of homogeneous supermanifold $G/H$ on $M$. Here $G$ is the Lie supergroup associated with the super Harish-Chandra pair $(G_{0},\gg)$ via the Kostant-Koszul construction (\cite{Kt, Kz}). The associated {\it Killing representation} of the algebra of supersymmetry 
\be
\label{purelygeor}
\hat{\varphi}:\gg\rightarrow\Der_{\bR}(\Lambda(\cS^{*}(M_{0})))
\ee
recognizes Killing {\it vectors} and generalized Killing {\it spinors} as "values" of even and odd Killing {\it fields}. An explicit description of the representation (\ref{purelygeor}) is obtained. Moreover, we give the super-version of the classical theory of invariant connections on a homogeneous manifold and we show that the generalized spin connection $\nabla^{\cS}$ defines a $G$-invariant superconnection on $M=G/H$. This could be of interest for future research along the lines of \cite{ACDS, K1, K2}. 
\\

The paper is structured as follows.
\\

The first section recalls the basic notions of supergeometry; in particular even and odd {\bf evaluations}
$$
\ev_{\ol{0}}:\mathcal{T}_{M}(U)\rightarrow\Gamma(U,TM_{0})\qquad,\qquad
\ev_{\ol{1}}:\mathcal{T}_{M}(U)\rightarrow\Gamma(U,(TM)_{\ol{1}})
$$
are defined, where $\mathcal{T}_{M}$ and $TM=(TM)_{\0}+(TM)_{\ou}$ are respectively the {\it tangent sheaf} and the {\it tangent bundle} of a supermanifold $M=(M_{0},\cA_{M})$. The even part $(TM)_{\0}$ of the tangent bundle is naturally isomorphic to the tangent bundle $TM_{0}$ of the underlying smooth manifold $M_{0}$ while the odd part $(TM)_{\ou}$, in the case $\cA_{M}\cong\Lambda(S^{*}(M_{0}))$, is the spin bundle $S(M_{0})$.
\\ 

The second section describes the Lie-Kostant-Koszul theory of Lie supergroups. A supermanifold $G=(G_{0},\mathcal{A}_{G})$ is a Lie supergroup if it is a group object in the category of supermanifolds. In the seminal paper \cite{Kt}, the notion of super Harish-Chandra pair (shortly sHC pair) is introduced and proved to be equivalent to the notion of Lie supergroup. A sHC pair is a pair $(G_{0},\gg)$, where $G_{0}$ is a Lie group and $\gg=\gg_{\0}+\gg_{\ou}$ a Lie superalgebra satisfying some consistency conditions. While this notion is probably the most efficient to prove theorems in the category of Lie supergroups and homogeneous supermanifolds, it has the disadvantage of obscuring the geometric meaning of the sheaf of superfunctions $\cA_{G}$ of the associated Lie supergroup $G=(G_{0},\cA_{G})$. \cite{BS, Kz} describe how to reconstruct the structure sheaf $\cA_{G}$ together with its Hopf superalgebra structure; in particular \cite{Kz} shows the existence of canonical isomorphisms
$$
\cA(G)\cong\Hom_{\cU(\gg_{\0})}(\cU(\gg),\mathcal{C}^{\infty}(G_{0}))\cong\mathcal{C}^{\infty}(G_{0})\otimes\Lambda(\gg_{\ou}^{*})
$$
and explicitly describes the representation of the Lie superalgebra $\gg=\gg_{\0}+\gg_{\ou}$
$$
\varphi:\gg\rightarrow\Der_{\bR}(\mathcal{C}^{\infty}(G_{0})\otimes\Lambda(\gg_{\ou}^{*}))
$$
by {\it left-invariant} vector fields 
(see Proposition \ref{koszulpari} and Theorem \ref{Koszulsx}). Theorem \ref{Koszulsx} has been proved in \cite{DP} with the aid of coalgebra theory. Using this approach, we obtain an analogous description for the representation 
$$
\hat{\varphi}:\gg\rightarrow\Der_{\bR}(\mathcal{C}^{\infty}(G_{0})\otimes\Lambda(\gg_{\ou}^{*}))
$$
by {\it right-invariant} vector fields (see Proposition \ref{parikill} and Theorem \ref{Andrea}).
\\

Section \ref{section3} recalls the basic definitions of action of a Lie supergroup and homogeneous supermanifold and it is an introduction to section \ref{section4} where the natural generalization to the category of supermanifolds of a classical theorem of Wang (\cite{KN1, KN2}) is obtained. Wang's theorem estabilishes a natural bijective correspondence between the set $\Conn(M_{0})^{G_{0}}$ of $G_{0}$-invariant linear connections on a homogeneous manifold $M_{0}=G_{0}/H_{0}$ with reductive decomposition $\gg_{0}=\gh_{0}+\gm_{0}$ and the set $\Hom_{\bR}(\gm_{0},\ggl_{\bR}(\gm_{0}))^{H_{0}}$ of Nomizu maps. 
The analogous result in the super-setting is the following.
\begin{theorem}
Let $M=G/H=(G_{0}/H_{0},\cA_{G/H})$ be a homogeneous supermanifold with reductive decomposition $\gg=\gh+\gm$. 
There is a bijective correspondence
$$
\Conn(M)^{G}\cong\Hom_{\bR}(\gm,\ggl_{\bR}(\gm))^{H}
$$
between the set $\Conn(M)^{G}$ of $G$-invariant linear connections on $M$ and the set $\Hom_{\bR}(\gm,\ggl_{\bR}(\gm))^{H}$ of Nomizu maps.
\end{theorem}
This result has already been stated in the literature in the case of even stability subgroup (\cite{Co}), but, to our knowledge, a proof is missing.
The dimension of the space of invariant linear connections on a Poincare' superspacetime in signature $(r,s)$ is determined.
\begin{theorem}
The dimension $D$ of the vector space of invariant connections on a Poincare' superspacetime of signature $(r,s)$ depends on $r-s \mod 8$ as follows
$$\left[
\begin{array}{c|c|c|c|c|c|c|c|c}
\ r-s \mod 8 & 1 & 2 & 3 & 4 & 5 & 6 & 7 & 8 \\
\hline
\ D & 12 & 24 & 12 & 24 & 12 & 6 & 3 & 6 \\
\end{array}
\right]$$
\end{theorem}
\phantom{c}
Section \ref{sezione5} deals with a general method to climb from the level of classical geometry to the level of supergeometry. This procedure is called {\bf superization}, which is the inverse of evaluation. 
To every reductive homogeneous pseudo-Riemannian spin manifold $M_{0}=G_{0}/H$, together with an adapted algebra of supersymmetry 
\be
\label{adapted}
\gg=\gg_{\0}+\gg_{\ou}=\gg_{0}+S
\ee
superization associates a homogeneous supermanifold $M=G/H$ whose sheaf of superfunctions $\cA_{M}$ is isomorphic to $\Lambda(S^{*}(M_{0}))$.
Each algebra (\ref{adapted}) determines a flat connection $\nabla^{\cS}$ in the spin bundle and generalized Killing spinors are defined as $\nabla^{\cS}$-parallel spinors.  
The Lie superalgebra (\ref{adapted}) admits a geometric representation (\ref{purelygeor}) whose underlying classical geometry is described by Killing vector/spinor maps
\be
\label{macheneso3}
\,\ev_{\0}\circ\hat{\varphi}|_{\gg_{\0}}:\gg_{\0}\rightarrow\mathcal{T}(M_{0})\qquad,\qquad x\mapsto\hat{\varphi}_{0}(x)
\ee
\be
\label{macheneso2}
\ev_{\ou}\circ\hat{\varphi}|_{\gg_{\ou}}:\gg_{\ou}\rightarrow\cS(M_{0})\qquad,\qquad s\mapsto\psi^{s}\quad.
\ee
We remark that (\ref{macheneso2}) does not depend on the odd-odd bracket $[S,S]$ but its lift to the anti-homomorphism (\ref{purelygeor}) does.
For example, for every $s\in S$, the supervector field $\hat{\varphi}(s)$ coincides with the algebraic interior product $\i_{\psi^{s}}\in\Der(\Lambda(\cS^{*}(M_{0})))$ only in the odd-commutative case, \textit{i.e.} when $[S,S]=0$. On the other hand, the action of $\hat{\varphi}(x)$ is given by the canonical connection $\nabla^{can}_{\hat{\varphi}_{0}(x)}\in\End(\cS^{*}(M_{0}))$ of the reductive homogeneous manifold $M_{0}=G_{0}/H$ considered as a derivation of $\Lambda(\cS^{*}(M_{0}))$. Each adapted algebra of supersymmetry (\ref{adapted}) also determines an embedding
\be
\label{siamoallafine}
\mathcal{T}(M_{0})\oplus\cS(M_{0})\hookrightarrow\Der_{\bR}(\Lambda(\cS^{*}(M_{0})))
\ee
$$
X+\psi\mapsto\mathbf{X}+\Psi
$$
such that $\ev_{\0}+\ev_{\ou}$ composed with (\ref{siamoallafine}) is the identity. The action of $\mathcal{T}(M_{0})$ is given by the flat spin connection $\nabla^{\cS}$ considered as derivation of $\Lambda(\cS^{*}(M_{0}))$ and $\Psi=\i_{\psi}$ if and only if odd-commutativity holds. The flat spin connection $\nabla^{\cS}$ can be interpreted, via the super-version of Wang's theorem, as a $G$-invariant superconnection on the homogeneous supermanifold $G/H$ whose associated Nomizu map is given by (\ref{supersymmetryconnection}).
We check that our setting is compatible with the construction of the Killing superalgebra in 11-dimensional bosonic supergravity, whose recipe is explained in \cite{F, F1, FP, T}. Our construction gives an interpretation of the generalized Killing spinors as values of geometric "infinitesimal odd symmetries" as hinted in \cite{T}. The study of adapted algebras of supersymmetry (with possibly as odd part a $\gh$-submodule of $S$) could become the primary object of interest in the search of supergravity solutions with many Killing spinors. Indeed each algebra of supersymmetry encodes the relevant object of the theory: the generalized spin connection $\nabla^{\cS}$.
\\

{\bf Notation.}
We mainly deal with $\bZ_{2}$-graded real vector spaces. We refer to them by "graded vector spaces" or "supervector spaces". On the other hand, a graded algebra is called a "superalgebra" only if supercommutativity holds. The symbol of exterior algebra is used {\bf only in the non graded} sense, while the same is not true for the symbol of symmetric algebra.
The set of global sections $\mathcal{F}_{M}(M_{0})$ of a sheaf $\mathcal{F}_{M}$ over $M_{0}$ is denoted by $\mathcal{F}(M)$. For example, the superalgebra of global superfunctions $\mathcal{A}_{M}(M_{0})$ of a supermanifold $M=(M_{0},\mathcal{A}_{M})$ is denoted by $\mathcal{A}(M)$. 
\section{Preliminaries of supergeometry} \label{LSG and HSM}
\setcounter{equation}{0}
This section recalls the basic notions of the theory of supermanifolds. 
For more details see \cite{DM, Kt, L, V}. 
\subsection{Linear superalgebra, Hopf algebra and convolution algebra}
\begin{definition}
\rm{A non-zero element $x\in V$ of a supervector space $V=V_{\overline{0}}+V_{\overline{1}}$ is {\it homogeneous} of parity $|x|\in\bZ_{2}=\left\{\0,\ou\right\}$ if $x\in V_{|x|}$. The dimension of $V$ is defined as the pair $\dim V:=\dim V_{\overline{0}}|\dim V_{\overline{1}}$. 
}\end{definition}
Let $V=V_{\0}+V_{\ou}$, $W=W_{\0}+W_{\ou}$ be two supervector spaces.
\begin{definition}
\rm{
A linear map $\varphi\in\Hom_{\bR}(V,W)$ is an {\it even} (resp. {\it odd}) {\it morphism} of supervector spaces if $\varphi(V_{\alpha})\subseteq W_{\alpha}$ (resp. $\varphi(V_{\alpha})\subseteq W_{\alpha+\ou}$), $\alpha\in\bZ_{2}$. 
The following decomposition holds
$$
\Hom_{\bR}(V,W)=\Hom_{\bR}(V,W)_{\0}\oplus \Hom_{\bR}(V,W)_{\ou}\quad.
$$
}\end{definition}
\begin{definition}
\rm{
The {\it tensor product} $V\otimes W$ is a supervector space with gradation $$(V\otimes W)_{\0}=(V_{\0}\otimes W_{\0})\oplus (V_{\ou}\otimes W_{\ou})\qquad,\qquad (V\otimes W)_{\ou}=(V_{\0}\otimes W_{\ou})\oplus (V_{\ou}\otimes W_{\0})\quad.$$
}\end{definition}
\begin{definition}
\rm{A supervector space $A=A_{\0}+A_{\ou}$ is a {\it (graded) algebra} if there exist {\it multiplication} $m_{A}\in\Hom_{\bR}(A\otimes A,A)_{\0}$ and {\it unity} $1_{A}\in\Hom_{\bR}(\bR,A)_{\0}$ such that $(A,m_{A},1_{A})$ is an algebra.
{\it Left} and {\it right multiplication} by $a\in A$ are given by
$$L_{a}:A\rightarrow A\qquad,\qquad R_{a}:A\rightarrow A$$
$$
\phantom{ccccccccccccccccccc}b\mapsto ab\phantom{ccccccccccccccc}\, b\mapsto (-1)^{|b||a|}ba\phantom{cccc}.
$$
The {\it supercommutator} is the even morphism $[\cdot,\cdot]\in\Hom_{\bR}(A\otimes A,A)_{\0}$ defined by $$[a,b]:=L_{a}b-R_{a}b=ab-(-1)^{|a||b|}ba\quad.$$
$(A,m_{A},1_{A})$ is a supercommutative superalgebra (shortly a {\it superalgebra}) if $[\cdot,\cdot]\equiv0$. 
}\end{definition}
\begin{example}
\rm{The exterior algebra $\Lambda_{\bR} V=\Lambda_{\bR}^{even}V+\Lambda_{\bR}^{odd}V$ over an even vector space $V=V_{\ol{0}}$ is a superalgebra.
}\end{example}
\begin{definition}
\label{defLie}
\rm{A {\it Lie superbracket} on a supervector space $\gg=\gg_{\ol{0}}+\gg_{\ol{1}}$ is an even morphism $[\cdot,\cdot]\in\Hom_{\bR}(\gg\otimes\gg,\gg)_{\0}$ satisfying
\be
[x,y]=-(-1)^{|x||y|}[y,x]
\ee
\be
\label{Jacobi}
[x,[y,z]]=[[x,y],z]+(-1)^{|x||y|}[y,[x,z]]
\ee
for all $x,y,z\in\gg$. The pair $(\gg,[\cdot,\cdot])$ is called a {\bf Lie superalgebra}.
}\end{definition}
The Lie superalgebra of endomorphisms $\End_{\bR}(V)$ of $V$ is denoted by $\ggl_{\bR}(V)$. 
The definitions of morphism of zero parity ({\it morphism} in short) and derivation of Lie superalgebras are the natural ones. 
\begin{definition}
\rm{Let $\gg=\gg_{\ol{0}}+\gg_{\ol{1}}$ be a Lie superalgebra. A {\it representation} $\varphi$ of $\gg$ on a supervector space $V=V_{\ol{0}}+V_{\ol{1}}$ is a morphism of Lie superalgebras $\varphi:\gg\rightarrow\ggl_{\bR}(V)$.
}\end{definition}
An anti-homomorphism $\varphi:\gg\rightarrow\ggl_{\bR}(V)$ is also called a "representation".
\begin{definition}[\cite{Sw}]
\label{coalgebra}
\rm{
A supervector space $C=C_{\0}+C_{\ou}$ is a (graded) {\bf coalgebra} if there exist {\it comultiplication} $\Delta_{C}\in\Hom_{\bR}(C,C\otimes C)_{\0}$ and {\it counity} $\epsilon_{C}\in\Hom_{\bR}(C,\bR)_{\0}$ satisfying
$$
\label{Coassoc}
\textbf{Coassociativity}\qquad\qquad (\Delta_{C}\otimes\Id_{C})\circ\Delta_{C}=(\Id_{C}\otimes\Delta_{C})\circ\Delta_{C}\quad,
$$
$$
\label{Counit}
\textbf{Counity\,\,}\qquad\qquad\quad (\epsilon_{C}\otimes\Id_{C})\circ\Delta_{C}=(\Id_{C}\otimes\epsilon_{C})\circ\Delta_{C}=\Id_{C}\quad.
$$
A coalgebra is denoted by the triple $(C,\Delta_{C},\epsilon_{C})$.
}\end{definition}
Sometimes the so-called {\bf sigma-notation} of \cite{Sw} is used:
$$
\Delta_{C}(c):=\sum_{i}c_{(1)}^{i}\otimes c_{(2)}^{i} 
$$
for every $c\in C$. Symbolically we write
$$
\Delta_{C}(c)=c_{(1)}\otimes c_{(2)}\quad,\quad
(\Delta_{C}\otimes\Id_{C})\circ\Delta_{C}(c):=c_{(1)}\otimes c_{(2)}\otimes c_{(3)}\quad,\quad\cdot\cdot\cdot
$$
A coalgebra $(C,\Delta_{C},\epsilon_{C})$ is {\it cocommutative} if
$
c_{(1)}\otimes c_{(2)}=(-1)^{|c_{(1)}||c_{(2)}|}c_{(2)}\otimes c_{(1)}
$.
\begin{definition}[\cite{R}]
\rm{
A {\it coderivation} of a coalgebra $(C,\Delta_{C},\epsilon_{C})$ is an endomorphism $\Phi\in\End_{\bR}(C,C)$ such that $\Delta_{C}\circ\Phi=(\Phi\otimes \Id+\Id\otimes \Phi)\circ\Delta_{C}$. The space of coderivations is denoted by $\CoDer(C)\ni\Phi$.
}\end{definition}
Let $(A,m_{A},1_{A})$ be a graded algebra and $(C,\Delta_{C},\epsilon_{C})$ a graded coalgebra.
\begin{definition}[\cite{Sw}]
\rm{
The space $\Hom_{\bR}(C,A)$ is a graded algebra, over the graded algebra $(C^{*},\Delta_{C}^{*},\epsilon_{C})$ dual to $C$, 
with respect to the 
{\it convolution product}
$$
F*G:=m_{A}\circ (F\otimes G)\circ \Delta_{C}\qquad\qquad F,G\in \Hom_{\bR}(C,A)\quad.
$$
The triple $(\Hom_{\bR}(C,A),*,\epsilon_{C})$ is the {\it convolution algebra} produced from $A$ and $C$.
}\end{definition}
The importance of coderivations lies in the following fact. The transpose action of a coderivation $\Phi\in\CoDer(C)$ is a derivation of the algebra $(\Hom_{\bR}(C,A),*,\epsilon_{C})$. Alternative descriptions of $\CoDer(C)$ are important; with some extra assumptions this space can be identified with another one, called space of {\it formal vector fields} (this is the content of Theorem \ref{R}).
\begin{definition}[\cite{Sw}]
\rm{
A (graded) {\bf bialgebra} $(B,m_{B},1_{B},\Delta_{B},\epsilon_{B})$ is a supervector space $B=B_{\0}+B_{\ou}$ which has the structure of a graded algebra $(B,m_{B},1_{B})$ and of a graded coalgebra $(B,\Delta_{B},\epsilon_{B})$ such that $\Delta_{B}$ and $\epsilon_{B}$ are morphisms of algebras. 
}\end{definition}
\begin{definition}[\cite{P}]
\rm{
The space $\cP(B)$ of {\it primitive elements} of a graded bialgebra $(B,m_{B},1_{B},\Delta_{B},\epsilon_{B})$ is defined by 
$\cP(B):=\left\{b\in B | \Delta_{B}(b)=1_{B}\otimes b + b\otimes 1_{B}\right\}$. The space $\Hom_{\bR}(B,\cP(B))$ is called the space of {\it formal vector fields} of $\cP(B)$.
}\end{definition}
\begin{definition}[\cite{Sw}]
\rm{An {\bf antipode} of a graded bialgebra $(H,m_{H},1_{H},\Delta_{H},\epsilon_{H})$ is a convolution inverse $\delta_{H}\in\Hom_{\bR}(H,H)$ of the identity $\Id_{H}$. A graded bialgebra with antipode $(H,m_{H},1_{H},\Delta_{H},\epsilon_{H},\delta_{H})$ is called a {\bf Hopf (graded) algebra}.
}\end{definition}
\begin{example}[\cite{Sc, Sw}]
\label{uni}
\rm{The universal enveloping algebra of a Lie superalgebra
$$
\cU(\gg):=T(\gg)/\left\langle x\otimes y-(-1)^{|x||y|}y\otimes x-[x,y]\quad| x,y\in\gg\right\rangle
$$
has a natural structure of a cocommutative Hopf graded algebra. In particular, the symmetric algebra $S(\gg)$ is a cocommutative Hopf superalgebra. The unique extension to $\cU(\gg)$ of a derivation $D\in\Der_{\bR}(\gg,\gg)$ of $\gg=\gg_{\0}+\gg_{\ou}$ is a coderivation. {\bf Notation}: the action of the antipode $\delta_{\cU(\gg)}$ on $u\in\mathcal{U}(\gg)$ is denoted by $u\mapsto\overline{u}$.
}\end{example}
The space $\CoDer(H)$ of coderivations of a cocommutative Hopf graded algebra $(H,m_{H},1_{H},\Delta_{H},\epsilon_{H},\delta_{H})$ is identified with the space of formal vector fields. Indeed
\begin{theorem}[\cite{R}]
\label{R}
The maps
$$
\Hom_{\bR}(H,\cP(H))\rightarrow \CoDer(H)
$$
\be
\label{allora}
\phantom{ccccccccccccccccc}\,\phi\mapsto \phantom{c}c_{\phi}:=\Id_{H}*\phi
\ee
\be
\label{allora2}
\phantom{cccccccccccccccccc}\phi\mapsto\phantom{c}_{\phi}c:=\phi*\Id_{H}
\ee
are bijections satisfying $\delta_{H}*\phantom{c}c_{\phi}=\phantom{c}_{\phi}c*\delta_{H}=\phi$. Moreover $(\ref{allora})=(\ref{allora2})$ whenever $(H,m_{H},1_{H})$ is a superalgebra. 
\end{theorem}
Let $(\gg,[\cdot,\cdot])$ be a Lie superalgebra and $(C,\Delta_{C},\epsilon_{C})$ a graded cocommutative coalgebra. 
\begin{definition}[\cite{P}]
\rm{
The space $\Hom_{\bR}(C,\gg)$ is a Lie superalgebra, over the graded algebra $(C^{*},\Delta_{C}^{*},\epsilon_{C})$ dual to $C$, with respect to the {\it convolution bracket} 
$$
[F,G]:=[\cdot,\cdot]\circ (F\otimes G)\circ \Delta_{C}\qquad\qquad F,G\in \Hom_{\bR}(C,\gg)\quad.
$$
The pair $(\Hom_{\bR}(C,\gg),[\cdot,\cdot])$ is the {\it convolution Lie superalgebra} produced from $\gg$ and $C$.
}\end{definition}
\subsection{Supermanifold}
\begin{definition}[\cite{Kt, L, V}]
\label{supermanifold}
\rm{A  {\bf supermanifold} $M$ of dimension $\dim M= m|n$ is a pair $(M_{0},\mathcal{A}_{M})$ where $M_{0}$ is a manifold of dimension $m$ (called the {\it body}) and 
$\mathcal{A}_{M}=(\mathcal{A}_{M})_{\ol{0}}\oplus(\mathcal{A}_{M})_{\ol{1}}$ is a sheaf of {\bf superfunctions}, \textit{i.e} a sheaf of superalgebras such that for all $p\in M_{0}$ there exists an open neighbourhood $U\owns p$ such that
$$
\mathcal{A}_{M}|_{U}\cong\mathcal{C}_{M_{0}}^{\infty}|_{U}\otimes \Lambda[\xi_{1},...,\xi_{n}]
$$
where $\mathcal{C}_{M_{0}}^{\infty}$ denotes the sheaf of smooth functions of $M_{0}$. Sections of $\mathcal{A}_{M}$ are called {\it superfunctions} of $M$. If coordinates $\{x^{i}\}$ on $U\subseteq M_{0}$ are given, the set
\be
\label{COORD}
\left\{\eta^{k}\right\}=\{x^{1},..,x^{m},\xi_{1},..,\xi_{n}\}
\ee
is called a set of {\it local coordinates} on $M$.
}\end{definition}
Locally any superfunction $f\in \mathcal{A}_{M}(U)$ is of the form 
\be
\label{local}
f=\sum_{\alpha}f_{\alpha}(x^{1},..,x^{m})\xi^{\alpha}
\ee
where $f_{\alpha}(x^{1},....,x^{m})\in \mathcal{C}^{\infty}_{M_{0}}(U)$ and $\xi^{\alpha}:=\xi_{1}^{\alpha_{1}}\wedge\cdot\cdot\cdot\wedge\xi_{n}^{\alpha_{n}}$, $\alpha=(\alpha_{1},...,\alpha_{n})\in\bZ_{2}^{n}$. The decomposition $f=f_{\ol{0}}+f_{\ol{1}}\in\mathcal{A}_{M}(U)_{\ol{0}}\oplus\mathcal{A}_{M}(U)_{\ol{1}}$ is locally given by
$$
f_{\ol{0}}=\sum_{|\alpha|=\0}f_{\alpha}(x^{1},..,x^{m})\xi^{\alpha}\quad,\quad f_{\ol{1}}=\sum_{|\alpha|=\ou}f_{\alpha}(x^{1},..,x^{m})\xi^{\alpha}\quad(|\alpha|:=\sum_{i=1}^{n}\alpha_{i}\in\bZ_{2})
$$
Definition \ref{supermanifold} implies the existence of a canonical epimorphism, called the {\bf evaluation map}, denoted by
$$
\ev:\mathcal{A}_{M}\rightarrow \mathcal{C}^{\infty}_{M_{0}}
$$ 
whose action on a superfunction (\ref{local}) is given by
$\tilde{f}:=\ev(f)=f_{0,...,0}(x^{1},...,x^{n})$.
\begin{example}
\rm{The supermanifold $\bR^{m,n}$ is the pair $(\bR^{m},\mathcal{A}_{\bR^{m,n}})$ where
$$
\mathcal{A}_{\bR^{m,n}}:=\mathcal{C}_{\bR^{m}}^{\infty}\otimes\Lambda[\xi_{1},...,\xi_{n}]\quad.
$$
}\end{example}
\begin{example}[\cite{Bat}]
\rm{Let $M_{0}$ be a real smooth manifold of dimension $m$, $E$ a rank $n$ vector bundle over $M_{0}$ and $\Lambda(E)$ the associated exterior bundle. The sheaf 
$$
\label{bat}
U\longrightarrow\mathcal{A}_{M}(U):=\Gamma(U,\Lambda(E))
$$
defines a supermanifold $M=(M_{0},\mathcal{A}_{M})$ of dimension $m|n$. Supermanifolds of this form are called {\it split} and are $\bZ$-graded. Every real smooth supermanifold is non canonically isomorphic to a split supermanifold. 
}\end{example}
\subsection{Morphism of supermanifolds and tensor sheaf}\hspace{0 cm}\newline\\
Let $M=(M_{0},\mathcal{A}_{M})$ and $N=(N_{0},\mathcal{A}_{N})$ be two supermanifolds.
\begin{definition}[\cite{Kt}]
\rm{A {\it morphism} $\phi\in\Mor(N,M)$ is a pair  
$\phi=(\phi_{0},\phi^{*})$ where $\phi_{0}:N_{0}\rightarrow M_{0}$ is a smooth map and $\phi^{*}:\mathcal{A}_{M}\rightarrow (\phi_{0})_{*}\mathcal{A}_{N}$ is a morphism of sheaves
of superalgebras over $M_{0}$ called the {\bf pull-back}. A {\it diffeomorphism} is a morphism $\phi=(\phi_{0},\phi^{*})$ such that $\phi_{0}$ is a diffeomorphism and 
$\phi^{*}$ is an isomorphism.
}\end{definition}
Every morphism $\phi=(\phi_{0},\phi^{*})\in\Mor(N,M)$ is uniquely defined by the pull-back $\phi^{*}: \mathcal{A}(M)\longrightarrow \mathcal{A}(N)$ on global sections of the structure sheaves (\cite{Kt}).
\begin{example}
\label{point}
\rm{The evaluation map
$$\ev:\mathcal{A}(M)\rightarrow\mathcal{C}^{\infty}(M_{0})$$
$$\phantom{c}\,f\mapsto \tilde{f}$$
defines the canonical embedding $\phi:=(id_{M_{0}},\ev)$ of the body $(M_{0},\mathcal{C}^{\infty}_{M_{0}})$ inside the supermanifold $(M_{0},\mathcal{A}_{M})$ while, for every $p\in M_{0}$, the {\it evaluation at p} $$\ev_{p}:\mathcal{A}(M)\longrightarrow \bR\phantom{ccccc}$$
$$\phantom{ccccccccccccccccccc}\,f\mapsto\tilde{f}(p)=:f(p)\phantom{ccccc}$$
defines the canonical embedding $p=(p_{0},\ev_{p})\in\Mor(\bR^{0,0},M)$. Recall that 
a superfunction $f\in\mathcal{A}(M)$ is not determined by its values $f(p)$, $p\in M_{0}$. 
For every fixed supermanifold $N$, denote by $\hat{p}\in\Mor(N,M)$ the {\it constant map} $p\in M_{0}$, \textit{i.e.} the composition of $p\in\Mor(\bR^{0,0},M)$ with the unique element of $\Mor(N,\bR^{0,0})$. 
}\end{example}
\begin{example}[\cite{L}]
\rm{
A morphism $\phi=(\phi_{0},\phi^{*})\in\Mor(\bR^{m,n},\bR^{m,n})$ is uniquely defined by the formulae
$$
\begin{cases}
\phi^{*}(x^{i})=x^{i}\circ\phi_{0}+\sum_{|\alpha|=\0}\phi^{i}_{\alpha}\xi^{\alpha} & \left\{\phi_{\alpha}^{i}\right\}_{|\alpha|=\0}\subseteq \mathcal{C}^{\infty}(\bR^{m})\\
\phi^{*}(\xi_{j})=\sum_{|\alpha|=\ou}\phi^{j}_{\alpha}\xi^{\alpha} & \left\{\phi_{\alpha}^{j}\right\}_{|\alpha|=\ou}\subseteq \mathcal{C}^{\infty}(\bR^{m})
\end{cases}
$$
and is a {\it change of coordinates} whenever it is a diffeomorphism.
}\end{example}
\begin{example}[\cite{Gr, L}]
\rm{
The structure sheaf of $M\times N$ is the completion $\overline{\mathcal{A}_{M}\otimes\mathcal{A}_{N}}$ of $\mathcal{A}_{M}\otimes\mathcal{A}_{N}$ with respect to Grothendieck's $\pi$-topology. The embedding $\mathcal{A}(M)\subseteq\overline{\mathcal{A}(M)\times\mathcal{A}(N)}$ defines the projection of $M\times N$ onto $M$.
}\end{example}
The sheaf $Der\mathcal{A}_{M}$ of derivations of $\mathcal{A}_{M}$ over $\bR$ is a sheaf of left $\mathcal{A}_{M}$-supermodules: $Der\mathcal{A}_{M}=(Der \mathcal{A}_{M})_{\ol{0}}\oplus(Der \mathcal{A}_{M})_{\ol{1}}$, where for $\alpha\in\bZ_{2}$
$$
(Der \mathcal{A}_{M})_{\alpha}:=\left\{ X\in End_{\bR}(\mathcal{A}_{M},\mathcal{A}_{M})_{\alpha}\,|\, X(fg)=X(f)g+(-1)^{\alpha|f|}fX(g)\right\}\quad.
$$
\begin{definition}[\cite{V}]
\label{tensorsheaf}
\rm{
The sheaf $\mathcal{T}_{M}:=Der\mathcal{A}_{M}$ is called the {\bf tangent sheaf} of the supermanifold $M=(M_{0},\mathcal{A}_{M})$. The sections of $\mathcal{T}_{M}$ are called {\it vector fields}. The {\it cotangent sheaf} is the sheaf $\mathcal{T}_{M}^{*}:=Hom_{\mathcal{A}_{M}}(\mathcal{T}_{M},\mathcal{A}_{M})$. The {\it tensor sheaf} $\oplus_{r,s}(\mathcal{T}_{M})^{r}_{s}$ is the sheaf of graded algebras generated by tensor products (graded over $\mathcal{A}_{M}$) of $\mathcal{T}_{M}$ and $\mathcal{T}_{M}^{*}$. Its sections are called {\it tensor fields}.
}\end{definition}
\begin{example}
\rm{A {\it metric} on $M$ is a non-degenerate symmetric even field $g:\mathcal{T}_{M}\otimes\mathcal{T}_{M}\rightarrow\mathcal{A}_{M}$ \textit{i.e.} a tensor field of type $(0,2)$ satisfying 
$$
{\bf (Non\phantom{c} degeneracy)\,\,\phantom{c}} 
{\rm The\phantom{c}map\phantom{c}} X\rightarrow g(X,\cdot)\phantom{c} {\rm is\phantom{c} an\phantom{c} isomorphism\phantom{c}} \mathcal{T}_{M}\rightarrow\mathcal{T}_{M}^{*},
$$
$$
{\bf (Symmetry)\phantom{ccccccccccccccccccc}} 
g(X,Y)=(-1)^{|X||Y|}g(Y,X),\phantom{cccccccccccccc}
$$
$$
{\bf (Parity)\phantom{ccccccccccccccccccccccc}} 
|g(X,Y)|=|X|+|Y|\quad.\phantom{cccccccccccccccccccc}
$$
}\end{example}
The tangent sheaf $\mathcal{T}_{M}$ is a sheaf of Lie superalgebras with respect to the bracket
\be
\label{bracket}
[X,Y]:=X\circ Y-(-1)^{|X||Y|}Y\circ X\qquad X,Y\in\mathcal{T}_{M}(U)
\ee
Any coordinate system (\ref{COORD}) gives rise to basic even $\frac{\partial}{\partial x^{i}}$ and odd $\frac{\partial}{\partial \xi_{j}}$ vector fields whose action on (\ref{local}) is 
\begin{displaymath}
\begin{array}{l}
\frac{\partial f}{\partial x^{i}}=\sum_{\alpha}\frac{\partial f_{\alpha}(x^{1},..,x^{m})}{\partial x^{i}}\xi^{\alpha},\\
\\
\frac{\partial f}{\partial \xi_{j}}=\sum_{\alpha}\alpha_{j}(-1)^{\alpha_{1}+\cdot\cdot+\alpha_{j-1}}f_{\alpha}(x^{1},..,x^{m})\xi_{1}^{\alpha_{1}}\wedge\cdot\cdot\wedge\xi_{j-1}^{\alpha_{j-1}}\wedge\xi_{j+1}^{\alpha_{j+1}}\wedge\cdot\cdot\wedge\xi_{n}^{\alpha_{n}}\quad.
\end{array}
\end{displaymath}
\subsection{
Even-value and odd-value}\hspace{0 cm}\newline\\
Let $\mathcal{A}_{M,p}$ denote the stalk of the sheaf $\mathcal{A}_{M}$ at $p\in M_{0}$.
\begin{definition}[\cite{L}]
\rm{The {\it tangent space} $T_{p}M=(T_{p}M)_{\ol{0}}\oplus (T_{p}M)_{\ol{1}}$ at $p\in M_{0}$ is the supervector space of dimension $\dim M=m|n$ defined by
$$
T_{p}M:=\left\{ X\in\Hom_{\bR}(\mathcal{A}_{M,p},\bR)\,|\,X(fg)=X(f)g(p)+f(p)X(g)\right\}\quad.
$$
Elements of $T_{p}M$ are called {\it tangent vectors} and $$TM:=\cup_{p\in M_{0}}T_{p}M=\cup_{p\in M_{0}}(T_{p}M)_{\ol{0}}+\cup_{p\in M_{0}}(T_{p}M)_{\ol{1}}=:(TM)_{\ol{0}}+(TM)_{\ol{1}}$$ is a graded vector bundle on $M_{0}$ called the {\bf tangent bundle} of $M$. The canonical isomorphism
$
(T_{p}M)_{\ol{0}}\cong T_{p}M_{0}
$
gives the canonical identification
$
(TM)_{\ol{0}}\cong TM_{0}
$. 
}\end{definition}
\begin{definition}
\label{valore}
\rm{
Let $X\in\mathcal{T}_{M}(U)$ be a vector field on $U\ni p$. The tangent vector $X|_{p}:=\ev_{p}\circ X:\mathcal{A}_{M,p}\rightarrow\bR$ is called the {\bf value}
of $X$ at the point $p$. The assignment $U\ni p\mapsto (X|_{p})_{\ol{0}}\in (T_{p}M)_{\ol{0}}\cong T_{p}M_{0}$ is a well-defined vector field $\tilde{X}_{\ol{0}}\in\mathcal{T}_{M_{0}}(U)$ of $M_{0}$ called the {\bf even-value} of $X$. The assignment $U\ni p\mapsto (X|_{p})_{\ol{1}}\in (T_{p}M)_{\ol{1}}$ is a well-defined section $\tilde{X}_{\ol{1}}\in\Gamma(U,(TM)_{\ol{1}})$ of $(TM)_{\ol{1}}$ called the {\bf odd-value} of $X$. 
}\end{definition}
Denote by
$$
\ev_{\ol{0}}:\mathcal{T}_{M}(U)\rightarrow\Gamma(U,TM_{0})
$$
\be
\label{even}
\phantom{ccccccccccc}X\mapsto \ev_{\ol{0}}X:=\tilde{X}_{\ol{0}}
\ee
$$
\phantom{cc}\ev_{\ol{1}}:\mathcal{T}_{M}(U)\rightarrow\Gamma(U,(TM)_{\ol{1}})
$$
\be
\label{odd}
\phantom{ccccccccccc}X\mapsto \ev_{\ol{1}}X:=\tilde{X}_{\ol{1}}
\ee
the operators which send a vector field $X\in\mathcal{T}_{M}(U)$ to its even and odd values. Note that (\ref{even}) and (\ref{odd}) do not uniquely determine $X$ and that (\ref{even}) is {\bf not} a Lie superalgebra morphism, unless $\dim M=m|0$. The above definitions are naturally extended to arbitrary tensor fields $T\in\mathcal{T}^{r}_{s}(M)$. The tensor space at a point $p\in M_{0}$ is denoted by $\oplus_{r,s}T_{p}M^{r}_{s}\ni T|_{p}$ and
$T(X_{1},...,X_{s})|_{p}=T|_{p}(X_{1}|_{p},...,X_{s}|_{p})$ for all $X_{1},...,X_{s}\in\mathcal{T}(M)$.

\subsection{\texorpdfstring{$\phi$}{Phi}-vector field 
}\hfill\newline\\
Let $\phi=(\phi_{0},\phi^{*})\in\Mor(N,M)$ be a morphism of supermanifolds.
\begin{definition}[\cite{CF}]
\rm{A {\bf $\phi$-vector field} on $U\subseteq M_{0}$ is a linear map 
\be
\label{alongv}
X:\mathcal{A}_{M}(U)\rightarrow \mathcal{A}_{N}(\phi^{-1}_{0}U)
\ee
such that its homogeneous components satisfy
$$
X(fg)=(Xf)\phi^{*}(g)+(-1)^{|X||f|}\phi^{*}(f)(Xg)
$$
for all $f,g\in\cA_{M}(U)$. 
The associated sheaf $\mathcal{T}_{\phi}$ over $M_{0}$ is a locally free sheaf of left $(\phi_{0})_{*}\mathcal{A}_{N}$-supermodules of rank $\dim M$. Any $X\in\mathcal{T}_{\phi}(U)$ can be uniquely written in coordinates (\ref{COORD}) as 
\be
\label{along}
X=\sum_{k} f^{k}\cdot(\phi^{*}\circ \frac{\partial}{\partial \eta^{k}})\qquad\qquad f^{k}\in\mathcal{A}_{N}(\phi_{0}^{-1}U)
\ee
}\end{definition}
\begin{definition}
\label{differential}
\rm{
The {\it differential} $\phi_{*}:\mathcal{T}(N)\rightarrow\mathcal{T}_{\phi}(M)$ 
is defined by 
$$\mathcal{T}(N)\ni Y\mapsto\phi_{*}Y:=Y\circ\phi^{*}\in\mathcal{T}_{\phi}(M)\quad.$$ 
Two vector fields $X\in\mathcal{T}(M)$, $Y\in\mathcal{T}(N)$ are {\it $\phi$-related} if $\phi_{*}Y=\phi^{*}\circ X\in\mathcal{T}_{\phi}(M)$. The {\it differential $\phi_{*,p}:T_{p}N\rightarrow T_{\phi_{0}(p)}M$ at the point $p\in N_{0}$} 
is defined by
$$
T_{p}N\ni v\mapsto v\circ \phi^{*}\in T_{\phi_{0}(p)}M
\quad.$$ 
The morphism $\phi$ is an {\it immersion} (resp. {\it submersion}) at $p\in N_{0}$ if $\phi_{*,p}$ is injective (resp. surjective). For the local description of these see \cite{L, V}.
}\end{definition}
If $\phi$ is a diffeomorphism the notation
$\phi_{*}Y:=(\phi^{-1})^{*}\circ Y\circ\phi^{*}\in\mathcal{T}(M)$ is used. In particular, there exists a natural action of the group of diffeomorphisms $\Aut(M)$ on $\mathcal{T}^{r}_{s}(M)$.
The {\it pull-back} $\phi^{*}T\in\mathcal{T}_{s}(N)$ of a covariant tensor field $T\in\mathcal{T}_{s}(M)$ under $\phi=(\phi_{0},\phi^{*})\in\Mor(N,M)$ is easily defined using Definition \ref{differential} and (\ref{along}). The usual functorial property holds.
\begin{lemma}
Let $\psi=(\psi_{0},\psi^{*})\in\Mor(L,N)$ and $\phi=(\phi_{0},\phi^{*})\in\Mor(N,M)$ be morphisms of supermanifolds. For all families of vector fields $\left\{Y_{i}\right\}\subseteq\mathcal{T}(N)$ and $\left\{Z_{i}\right\}\subseteq\mathcal{T}(L)$ satisfying
$
(\phi\circ\psi)_{*}(Z_{i})
=\psi^{*}\circ Y_{i}\circ \phi^{*}
$, we have that
$$
\label{pull-back2}
\psi^{*}((\phi^{*}T)(Y_{1},...,Y_{s}))=((\phi\circ\psi)^{*}T)(Z_{1},...,Z_{s})
$$
\end{lemma}
\subsection{Lie derivative}\hspace{0 cm}\newline\\
\mbox{For every vector field $Y\in\mathcal{T}(N)$ there exists a unique derivation $\mathcal{L}_{Y}$ of the tensor} sheaf $\oplus_{r,s}(\mathcal{T}_{N})_{s}^{r}$ commuting with contractions and such that
$$
\mathcal{L}_{Y}(f)=Yf\qquad,\qquad \mathcal{L}_{Y}X=[Y,X]
$$
for every $f\in\mathcal{A}(N)$ and $X\in\mathcal{T}(N)$.  
\begin{definition}\rm{
The representation $\mathcal{L}:\mathcal{T}_{N}\rightarrow\ggl_{\bR}((\mathcal{T}_{N})_{s}^{r})$ is called {\it Lie derivative}. }
\end{definition}
Let $\left\{Y_{i}\right\}_{i=1}^{j\in\bN}$ be a set of vector fields on a supermanifold $N=(N_{0},\cA_{N})$ of dimension $\dim N=m|n$ such that 
$
\Span_{\bR}[Y_{1}|_{p},..,Y_{j}|_{p}]=(T_{p}N)_{\ou}
$
for every $p\in N_{0}$. 
\begin{lemma}
\label{appendixA2}
A tensor field $T\in\mathcal{T}(N)^{r}_{s}$ on $N$ satisfying
\begin{itemize}
\item[i)] $T|p=0$ for every $p\in N_{0}$,
\item[ii)] $(\mathcal{L}_{Y_{i_{k}}}\cdot\cdot\cdot\mathcal{L}_{Y_{i_{1}}}(T))|_{p}=0$ for every $p\in N_{0}$ and $1\leq k \leq n$
\end{itemize}
is zero.
\end{lemma}
\begin{pf}
See the Appendix. 
\end{pf}
If $Y\circ\phi^{*}=0$, the {\it Lie derivative} of a $\phi$-vector field $X\in\mathcal{T}_{\phi}(M)$ is defined by 
$$\mathcal{L}_{Y}(X):=Y\circ X\in\mathcal{T}_{\phi}(M)\quad.$$
\section{Lie-Kostant-Koszul theory of supergroups}
\setcounter{equation}{0}
There are various approaches to the theory of Lie supergroups. Subsection \ref{Apparizione} deals with the categorical approach, \textit{i.e.} a Lie  supergroup is a group object $G=(G_{0},\cA_{G})$ in the category of supermanifolds. Subsection \ref{harishharish} describes a more algebraic approach. The notion of super Harish-Chandra pair $(G_{0},\gg)$ is introduced and, following the ideas of \cite{BS, Kz}, it is recalled how to explicitly reconstruct the structure sheaf of a Lie supergroup, together with its Hopf superalgebra structure. Lie supergroups are thus equivalent to a global even part together with an infinitesimal odd part and, for this reason, most of the results in Lie supergroup theory find their natural setting here. This approach is refined in subsection \ref{accanto} where an important, more geometric picture is described. Here the behaviour of odd infinitesimal symmetries is better understood and some new formulas related to it are obtained.
\subsection{Lie supergroup}
\label{Apparizione}
\begin{definition}[\cite{V}]
\label{Lie}
\rm{A supermanifold $G=(G_{0},\mathcal{A}_{G})$ is called a {\it Lie supergroup} if there exist morphisms of supermanifolds
$$m=(m_{0},m^{*}):G\times G\rightarrow G\quad,\quad
i=(i_{0},i^{*}):G\rightarrow G\quad,\quad
e=(e_{0},\ev_{e}):\bR^{0,0}\rightarrow G$$
satisfying the usual axioms of associativity, identity and taking inverses:
\begin{itemize}
\item[$i)$]
$m\circ(id_{G}\times m)=m\circ(m\times id_{G}):G\times G\times G\longrightarrow G$,
\item[$ii)$]
$m\circ\left\langle id_{G},\hat{e}\right\rangle=m\circ\left\langle \hat{e},id_{G}\right\rangle=id_{G}:G\longrightarrow G$,
\item[$iii)$]
$m\circ\left\langle id_{G},i\right\rangle=m\circ\left\langle i,id_{G}\right\rangle=\hat{e}:G\longrightarrow G.$
\end{itemize}
}\end{definition}
An equivalent definition is that the infinite-dimensional superalgebra $\mathcal{A}(G)$ has a structure of (completed) Hopf superalgebra. 
\begin{definition}[\cite{Kt}]
\label{supergrouphopf}
\rm{A supermanifold $G=(G_{0},\mathcal{A}_{G})$ is called a {\it Lie supergroup} if there exist morphisms of superalgebras 
\begin{itemize}
\item[]$\textbf{Comultiplication}$
$m^{*}:\mathcal{A}(G)\longrightarrow\overline{\mathcal{A}(G)\otimes\mathcal{A}(G)}$ ,
\item[]$\textbf{Counity}$
$\ev_{e}:\mathcal{A}(G)\longrightarrow \bR$,
\item[]$\textbf{Antipode}$
$i^{*}:\mathcal{A}(G)\longrightarrow\mathcal{A}(G)$,
\end{itemize}
satisfying the following axioms
$$
(\Id\otimes m^{*})\circ m^{*}=(m^{*}\otimes \Id)\circ m^{*}:\mathcal{A}(G)\rightarrow\overline{\mathcal{A}(G)\otimes\mathcal{A}(G)\otimes\mathcal{A}(G)},
$$
$$
(\Id\otimes \ev_{e})\circ m^{*}=(\ev_{e}\otimes \Id)\circ m^{*}=\Id:\mathcal{A}(G)\longrightarrow \mathcal{A}(G),
$$
$$
m_{\cA_{G}}\circ(\Id\otimes i^{*})\circ m^{*}=m_{\cA_{G}}\circ(i^{*}\otimes\Id)\circ m^{*}=\ev_{e}:\mathcal{A}(G)\longrightarrow \bR\subseteq \mathcal{A}(G).
$$
}\end{definition}
Note that Definition \ref{Lie} and Definition \ref{supergrouphopf} reduce to the definition of Lie group if $\dim M=m|0$. In general the body $(G_{0},m_{0},i_{0},e_{0})$ is a Lie group.\\
For every $g\in G_{0}$ {\bf left/right translations} by $g$ are the diffeomorphisms
$$
L_{g}:G\cong\{g\}\times G\hookrightarrow G\times G\stackrel{m}{\rightarrow}G\qquad,\qquad
R_{g}:G\cong G\times \{g\}\hookrightarrow G\times G\stackrel{m}{\rightarrow}G
$$
whose pull-backs on global superfunctions are
$$
L_{g}^{*}:\mathcal{A}(G)\longrightarrow \mathcal{A}(G)\qquad,\qquad R_{g}^{*}:\mathcal{A}(G)\longrightarrow \mathcal{A}(G)
$$
$$
f\mapsto (\ev_{g}\otimes \Id)(m^{*}f)\qquad\qquad f\mapsto (\Id\otimes \ev_{g})(m^{*}f).
$$
The map $R_{G_{0}}:G_{0}\rightarrow \Aut(\mathcal{A}(G))$ (resp. $L_{G_{0}}:G_{0}\rightarrow \Aut(\mathcal{A}(G))$) defined by $R_{G_{0}}(g):=R_{g}^{*}$
(resp. $L_{G_{0}}(g):=L_{g}^{*}$) is a group (resp. anti-) homomorphism.
\begin{definition}[\cite{V}]
\label{vara}
\rm{A vector field $A\in\mathcal{T}(G)$ on $G$ is said {\it left} (resp {\it right})-{\it invariant} if
\be
\label{lefta}
(\Id\otimes A)\circ m^{*}=m^{*}\circ A\qquad\qquad ({\rm resp.}\phantom{c} (A\otimes \Id)\circ m^{*}=m^{*}\circ A)
\ee
The {\it Lie superalgebra} $\gg=\gg_{\0}+\gg_{\ou}$ (resp. $\hat{\gg}=\hat{\gg}_{\0}+\hat{\gg}_{\ou}$) is the supervector space of all left (resp. right)-invariant vector fields on $G$ with bracket (\ref{bracket}).
}\end{definition}
\begin{lemma}[\cite{V}]
Let $G=(G_{0},\mathcal{A}_{G})$ be a Lie supergroup with Lie superalgebra $\gg=\gg_{\0}+\gg_{\ou}$. The evaluation of a left-invariant vector field at the point $e\in G_{0}$
$$
\gg\longrightarrow T_{e}G\quad,\quad
A\mapsto A|_{e}
$$
is an isomorphism of supervector spaces. The inverse map is given by
$$
T_{e}G\longrightarrow\gg\quad,\quad
v \mapsto A_{v}:=(\Id\otimes v)\circ m^{*}
$$
A similar result holds for the right-invariant vector field $_{v}A\in\hat{\gg}$.
\end{lemma}
If $A\in\gg_{\0}$, the even value $\ev_{\0}\circ A$ is a left-invariant vector field on $G_{0}$. 
Note that invariance by left translations with respect to $G_{0}$ is not equivalent to (\ref{lefta}):
\begin{lemma}
\label{allora3}A vector field $A\in\mathcal{T}(G)$ on $G$ is left-invariant if and only if 
\beq
\label{prima}
L_{g}^{*}\circ A=A\circ L_{g}^{*}\qquad,\qquad
\mathcal{L}_{B}A=0
\eeq
for all $g\in G_{0}$ and (odd) $B\in\hat{\gg}$.
\end{lemma}
\begin{pf}
Let $A$ be a left-invariant vector field. It is easy 
to show that it satisfies equations (\ref{prima}).
Conversely suppose $A\in\mathcal{T}(G)$ satisfies (\ref{prima}) and consider the left-invariant vector field $A_{v}$ associated with $v:=A|_{e}\in T_{e}G$. The vector field $A-A_{v}$ satisfies the hypothesis of Lemma \ref{appendixA2} and then $A=A_{v}$.
\end{pf}
The sheaf of derivations of a Lie supergroup is trivial.
\begin{lemma}
\label{trivial}
The sheaf of derivations $\mathcal{T}_{G}$ of a Lie supergroup $G$ is globally trivial; more precisely
$
\mathcal{A}_{G}\otimes\gg\cong\mathcal{T}_{G}
$.
\end{lemma}
\begin{pf}
It follows directly from Nakayama's Lemma (\cite{V}).
\end{pf}
\subsection{Super Harish-Chandra pair}\hspace{0 cm}\newline\\
\label{harishharish}
An important equivalent way of defining a Lie supergroup is the following.
\begin{definition}[\cite{Kt}]
\label{superharish}
\rm{A pair $(G_{0},\gg)$ consisting of a Lie group $G_{0}$ and a Lie superalgebra $\gg=\gg_{\ol{0}}+\gg_{\ol{1}}$ is a \textbf{super Harish-Chandra} (shortly {\bf sHC}) {\bf pair} if $\gg_{\0}=\gg_{0}$ and if there exists an {\bf adjoint action}, \textit{i.e.} a morphism of Lie groups
\be
\label{grazie}
\Ad:G_{0}\longrightarrow\Aut(\gg)
\ee
such that
$
\Ad:G_{0}\longrightarrow\Aut(\gg_{\ol{0}})
$
is the usual adjoint action and 
$$
\ad_{B}A
=\frac{d}{dt}|_{t=0}\Ad_{\exp(tB)}A
$$
for every $A\in\gg$ and $B\in\gg_{\0}$.
}\end{definition}
\begin{theorem}[\cite{Kt}]
\label{kostantone}
Any Lie supergroup $G=(G_{0},\mathcal{A}_{G})$ defines a sHC pair $(G_{0},\gg)$ where
$$
\Ad:G_{0}\longrightarrow\Aut(\gg)
\qquad,\qquad g\mapsto \Ad_{g}:=(A\rightarrow R_{g}^{*}\circ A \circ R_{g^{-1}}^{*})\quad.
$$
The associated correspondence 
$
(G_{0},\mathcal{A}_{G})\longrightarrow (G_{0},\gg)
$
is a bijection between the sets of Lie supergroups and of sHC pairs.
\end{theorem}
Due to Theorem \ref{kostantone}, Lie supergroups and sHC pairs are synonymous for us.
\begin{example}
\label{esempiopoincare}
\rm{
Denote by $\bR^{r,s}$ the vector space $\bR^{r+s}$ endowed with the standard inner product $\left\langle \cdot,\cdot\right\rangle
$ of signature $(r,s)$ and by $$\mathrm{Spin}_{r,s}^{0}\subseteq\mathrm{Spin}_{r,s}\subseteq\mathrm{Cl}^{0}_{r,s}\subseteq(\mathrm{Cl}_{r,s}, \cdot)\qquad,\qquad\bR^{r,s}\subseteq\mathrm{Cl}_{r,s}$$
the inclusions of the (connected) Spin group inside (the even part of) the Clifford algebra $(\mathrm{Cl}_{r,s}, \cdot)$ and of $\bR^{r,s}$ inside $\mathrm{Cl}_{r,s}$. The 2-fold cover 
$
\xi:\mathrm{Spin}_{r,s}\rightarrow\mathrm{SO}_{r,s}
$
$$
\xi(g):\bR^{r,s}\rightarrow\bR^{r,s}
$$
$$
\phantom{ccccccccccccccccccccccccccccccccccccccc}v\mapsto g\cdot v\cdot g^{-1}\qquad g\in\mathrm{Spin}_{r,s},\phantom{c}v\in\bR^{r,s}
$$
is the {\it vector representation} of $\mathrm{Spin}_{r,s}$. It induces an isomorphism of Lie algebras 
$$
\xi_{*}:\spin_{r,s}\rightarrow\so_{r,s}\qquad,\qquad
\frac{1}{4}[v_{1},v_{2}]\mapsto v_{1}\wedge v_{2}
$$
where $v_{1},v_{2}\in\bR^{r,s}$. The {\it real spin representation} is the restriction to the spin group $\Spin_{r,s}\subseteq\mathrm{Cl}_{r,s}$ of an irreducible real representation of the Clifford algebra $\mathrm{Cl}_{r,s}$
$$
\Delta_{r,s}:\mathrm{Cl}_{r,s}\rightarrow\End_{\bR}(S)\quad.
$$
The spin representation $\Delta|_{\Spin_{r,s}}:\Spin_{r,s}\rightarrow\End_{\bR}(S)$ is either irreducible or it decomposes into a sum of two irreducible representations $S^{\pm}$ depending on $r-s$ mod $8$ (see \cite{LM} for details). 
The {\it Poincare' Lie superalgebra} $\gg(\Gamma)=\gg_{\0}+\gg_{\ou}$ associated with a $\spin_{r,s}$-invariant symmetric bilinear map
\be
\label{mah}
\Gamma:S\vee S\rightarrow \bR^{r,s}
\ee
is given by the Poincare' algebra
$$
\phantom{ccccccc}\gg_{\ol{0}}=\spin_{r,s}\inplus\bR^{r,s}
$$
as the even part, the odd part
$
\gg_{\ol{1}}=S
$
considered as a $\gg_{\0}$-module via
$$
[A,s]:=\Delta_{r,s}(A)s\qquad,\qquad[\bR^{r,s},S]:=0\qquad\qquad A\in\spin_{r,s},\phantom{c} s\in S
$$
and a bracket $[\cdot,\cdot]:S\vee S\rightarrow \bR^{r,s}$ given by (\ref{mah}). 
The associated {\it Poincare' Lie supergroup} is a sHC pair $(G_{0},\gg)$, where $$G_{0}=\Spin_{r,s}^{0}\ltimes\bR^{r,s}\qquad,\qquad\gg=\gg(\Gamma)=(\spin_{r,s}\inplus\bR^{r,s})+S$$ and the morphism
$
\Ad:G_{0}\longrightarrow\Aut(\gg)
$
is defined, for $g\in\Spin_{r,s}^{0}$ and $v\in\bR^{r,s}$ by
$$
\Ad_{g}(A+w+s):=(g\cdot A\cdot g^{-1}+g\cdot w\cdot g^{-1}+\Delta_{r,s}(g)s)
$$
$$
\Ad_{v}(A+w+s):=(A-\xi_{*}(A)v+w+s)\phantom{cccccccccccccc}
$$
where $A\in\spin_{r,s}, w\in\bR^{r,s}, s\in S$.
}\end{example}
\cite{BS, Kz} show how to reconstruct, from a given sHC pair $(G_{0},\gg)$, the structure sheaf of the associated Lie supergroup, together with its Hopf algebra structure. Note that the representation of $\gg_{\0}$ as left-invariant vector fields of the Lie group $G_{0}$
\be
\label{omomorfismo}
\varphi_{0}:\gg_{\ol{0}}\rightarrow \mathcal{T}(G_{0})
\ee
defines a structure of sheaf of left $\mathcal{U}(\gg_{\ol{0}})$-modules on  $\cC^{\infty}:=(\mathcal{C}^{\infty}_{G_{0}},m_{\cC^{\infty}_{G_{0}}},1_{\cC^{\infty}_{G_{0}}})$
and that the universal enveloping algebra $(\mathcal{U}(\gg),m_{\cU(\gg)},1_{\cU(\gg)},\Delta_{\cU(\gg)},\delta_{\cU(\gg)})$ of $\gg$ has the natural structure of left $\mathcal{U}(\gg_{\ol{0}})$-supermodule given by the embedding $\mathcal{U}(\gg_{\0})\subseteq\mathcal{U}(\gg)$. 
\begin{theorem}[\cite{BS, Kz}]
\label{harish}
The structure sheaf $\mathcal{A}_{G}$ of the Lie supergroup $(G_{0},\mathcal{A}_{G})$ associated with a sHC pair $(G_{0},\gg)$ is the sheaf of convolution algebras given by
$$
U\longrightarrow \mathcal{A}_{G}(U):=\Hom_{\mathcal{U}(\gg_{\ol{0}})}(\mathcal{U}(\gg),\mathcal{C}^{\infty}(U))\phantom{ccc},\phantom{ccc}
F_{1}* F_{2}:=m_{\mathcal{C}^{\infty}}\circ(F_{1}\otimes F_{2})\circ\Delta_{\mathcal{U}(\gg)}
$$
where $F_{1}$, $F_{2}\in\mathcal{A}_{G}(U)$. The embedding
$\mathcal{U}(\gg_{\ol{0}})\subseteq\mathcal{U}(\gg)$ induces the evaluation map
$$
\mathcal{A}_{G}(U)=\Hom_{\mathcal{U}(\gg_{\ol{0}})}(\mathcal{U}(\gg),\mathcal{C}^{\infty}(U))\rightarrow \Hom_{\mathcal{U}(\gg_{\ol{0}})}(\mathcal{U}(\gg_{\ol{0}}),\mathcal{C}^{\infty}(U))\cong\mathcal{C}^{\infty}(U)
$$
$$
\phantom{cccccc}\qquad F\mapsto F(1_{\cU(\gg)})\qquad.
$$
The structure of Hopf superalgebra of $\cA(G)$ is given by the algebra morphisms:
\begin{itemize}
\item[1)]
{\bf Comultiplication
$
m^{*}:\mathcal{A}(G)\rightarrow \mathcal{A}(G\times G)\quad,\quad
F\mapsto m^{*}F$}
$$
(m^{*}F)(u\otimes v)(g,h):=F((\Ad_{h^{-1}}u)\cdot v)(gh)
$$
\item[2)]
{\bf Counity
$
\ev_{e}:\mathcal{A}(G)\rightarrow \bR\quad,\quad
F\mapsto\ev_{e}F
$}
$$
\ev_{e}F:=F(1_{\cU(\gg)})(e)
$$
\item[3)]
{\bf Antipode $
i^{*}:\mathcal{A}(G)\rightarrow\mathcal{A}(G)\quad,\quad
F\mapsto i^{*}F:=F^{*}$}
$$
F^{*}(u)(g):=F(\Ad_{g}\overline{u})(g^{-1}) 
$$
\end{itemize}
where $\Ad:G_{0}\rightarrow \Aut(\mathcal{U}(\gg))$ is the unique extension of (\ref{grazie}), $e\in G_{0}$ is the unity of the Lie group $G_{0}$ and $u,v\in \mathcal{U}(\gg)$, $g,h\in G_{0}$.
The Lie superalgebra morphism
\be
\label{representation1}
\varphi:\gg\longrightarrow \mathcal{T}(G)\quad,\quad
\varphi(a)F:=(u\mapsto (-1)^{|a|}F(ua))
\ee
is the representation of $\gg$ as
left-invariant vector fields of the Lie supergroup $G$. The right-invariant vector field associated with $a\in\gg$ is given by
$$
\mathcal{A}(G)\ni
F\mapsto (\varphi(\overline{a})F^{*})^{*}\in\mathcal{A}(G)
$$
so that the Lie superalgebra anti-homomorphism
\be\hat{\varphi}:\gg\rightarrow\mathcal{T}(G)\quad,\quad
\label{representation2}
\hat{\varphi}(a)(F)(u)(g)=(-1)^{|F||a|}F((\Ad_{g^{-1}}a)\cdot u)(g)
\ee
is the representation of $\gg$ as
right-invariant vector fields of the Lie supergroup $G$.
\end{theorem}
\begin{remark}[\cite{Kz}]\rm{The previous construction can be generalized to associate with any $\gg_{\0}$-manifold (\textit{i.e.} a manifold $M_{0}$ with a representation $\varphi_{0}:\gg_{\0}\rightarrow\mathcal{T}(M_{0})$) a supermanifold $M=(M_{0},\cA_{M})$ with a representation $\varphi:\gg\rightarrow\mathcal{T}(M)$ such that $\ev_{\0}\circ\varphi|_{\gg_{\0}}=\varphi_{0}$ and $\ev_{p}\circ\varphi|_{\gg_{\ol{1}}}:\gg_{\ol{1}}\rightarrow (T_{p}M)_{\ol{1}}$ is an isomorphism for every $p\in M_{0}$.
}
\end{remark}
For the sake of completeness, here is how to reconstruct a morphism of Lie supergroups from the associated morphism of sHC pairs.
Let $(G_{0},\gg)$ and $(H_{0},\gh)$ be two sHC pairs and, for simplicity, assume that $G_{0}$ is connected. A morphism $\phi:\gg\rightarrow\gh$ of Lie superalgebras whose even part integrates to a morphism of Lie groups
$
\phi_{0}:G_{0}\rightarrow H_{0}
$
defines a morphism of Lie supergroups
$$
\phi^{*}:\mathcal{A}(H)\rightarrow \mathcal{A}(G)
$$
given by
$$
(\phi^{*}F)(u)(g):=F(\phi u)(\phi_{0}g)
$$
where $F\in \mathcal{A}(H)$, $u\in \mathcal{U}(\gg)$, $g\in G_{0}$.
\subsection{Koszul theory}\hfill\newline\\
\label{accanto}
Theorem \ref{harish} gives a very elegant way to reconstruct the Hopf superalgebra structure of a Lie supergroup $G$ from its sHC pair $(G_{0},\gg)$. On the other hand the structure sheaf of a Lie supergroup is trivial and thus, in particular, $\bZ$-graded.
\begin{theorem}[\cite{Kz}]
\label{split}
Let $G$ be a Lie supergroup and $\gg=\gg_{\ol{0}}+\gg_{\ol{1}}$ its Lie superalgebra. There exists a canonical isomorphism
$
\mathcal{A}_{G}\cong \mathcal{C}^{\infty}_{G_{0}}\otimes\Lambda(\gg_{\ol{1}}^{*})
$
\end{theorem}
\begin{pf}
It is important that for any $a_{1},...,a_{p}\in\gg_{\ol{1}}$
$$
\gamma:\Lambda(\gg_{\ol{1}})\hookrightarrow \mathcal{U}(\gg)\qquad,\qquad a_{1}\wedge\cdot\cdot\wedge a_{p}\mapsto\frac{1}{p!}\sum_{\sigma\in \sum_{p}}\sgn(\sigma)a_{\sigma(1)}\cdot\cdot\cdot a_{\sigma(p)}
$$
is a homomorphism of coalgebras. This implies that
\be
\label{symmetrization2}
\underline{\gamma}:\mathcal{U}(\gg_{\ol{0}})\otimes\Lambda(\gg_{\ol{1}})\longrightarrow \mathcal{U}(\gg)\qquad,
\qquad u\otimes v\mapsto u\cdot\gamma(v)
\ee
is an isomorphism of coalgebras and of left $\mathcal{U}(\gg_{\ol{0}})$-supermodules (see Theorem \ref{petracciona}). Thus
$$
\cA_{G}(U)=\Hom_{\mathcal{U}(\gg_{\ol{0}})}(\mathcal{U}(\gg),\mathcal{C}^{\infty}(U))\rightarrow \Hom_{\bR}(\Lambda(\gg_{\ol{1}}),\mathcal{C}^{\infty}(U))
$$
$$
\phantom{cccccccccccccccccccccccccccccccccc}F\mapsto F\circ\underline{\gamma}|_{\Lambda(\gg_{\ou})}=F\circ \gamma=:f
$$
is an isomorphism of convolution superalgebras. The superalgebra identification
$$
\Hom_{\bR}(\Lambda(\gg_{\ol{1}}),\mathcal{C}^{\infty}(U))\rightarrow\mathcal{C}^{\infty}(U)\otimes\Lambda(\gg_{\ol{1}})^{*}\rightarrow\mathcal{C}^{\infty}(U)\otimes\Lambda(\gg_{\ol{1}}^{*})
$$
shows that $G=(G_{0},\cA_{G})$ is a split supermanifold with respect to the (trivial) vector bundle $G_{0}\times\gg_{\ol{1}}^{*}$ over $G_{0}$. Recall that, due to the rule of signs of supergeometry, the last isomorphism
$
\Lambda(\gg_{\ol{1}})^{*}\cong \Lambda(\gg_{\ol{1}}^{*})
$
is given by
$$
(a_{1}\wedge\cdot\cdot\wedge a_{p})^{*}\mapsto (-1)^{(p-1)p/2}a_{1}^{*}\wedge\cdot\cdot\wedge a_{p}^{*}
$$
for $a_{1},...,a_{p}\in\gg_{\ou}$. 
\end{pf}
The algebra isomorphism  
\be
\label{koskos}
\mathcal{A}_{G}(U)\cong\mathcal{C}^{\infty}(U)\otimes\Lambda(\gg_{\ol{1}}^{*})
\ee
can be used to construct a canonical embedding of superalgebras
\be
\label{embed1}
\mathcal{C}^{\infty}(U)\hookrightarrow\mathcal{C}^{\infty}(U)\otimes\Lambda(\gg_{\ol{1}}^{*})\cong\mathcal{A}_{G}(U)
\ee 
Given a smooth function $f\in \mathcal{C}^{\infty}(U)$, the corresponding superfunction $f\in\mathcal{A}_{G}(U)$ (denoted by the same symbol, by abuse of notation) is uniquely determined by 
\be
\label{funzioni}
f(1_{\cU(\gg)})=f\in\mathcal{C}^{\infty}(U)\qquad,\qquad f(\Im\gamma-\bR 1_{\cU(\gg)})=0
\ee
In particular, for every $a_{1},a_{2}\in\gg_{\ol{1}}$
$$
f(a_{1})=0\quad,\quad
f(2a_{1}\cdot a_{2})=f([a_{1},a_{2}])=\varphi_{0}([a_{1},a_{2}])f\in\mathcal{C}^{\infty}(U)
$$
and for every $G\in\mathcal{A}_{G}(U)$ and $a_{1},...,a_{p}\in\Lambda(\gg_{\ol{1}})$
$$
(f*G)(\gamma(a_{1}\wedge\cdot\cdot\wedge a_{p}))=f\cdot G(\gamma(a_{1}\wedge\cdot\cdot\wedge a_{p}))=f\cdot g(a_{1}\wedge\cdot\cdot\wedge a_{p})\quad.
$$
\begin{corollary}
\label{preserv}
For every $g,h\in G_{0}$, $u\in \mathcal{U}(\gg)$ and $F\in \mathcal{A}(G)$
$$
(L_{h}^{*}F)(u)(g)=F(u)(hg)\quad,
$$
$$
\phantom{cccccci}(R_{h}^{*}F)(u)(g)=F(\Ad_{h^{-1}}u)(gh)\quad.
$$
The antipode, left and right translations with respect to $G_{0}$ preserve the $\bZ$-grading of $\cA_{G}$ and, in particular, the embedding (\ref{embed1}). 
\end{corollary}\begin{pf}
The first part of the Corollary is a straightforward consequence of Theorem \ref{harish}. 
The last assertions follow from the equations
$$
\Ad_{h}\circ\gamma=\gamma\circ \Ad_{h}\qquad,\qquad \overline{\gamma(a_{1}\wedge\cdot\cdot\wedge a_{p})}=(-1)^{p}\gamma(a_{1}\wedge\cdot\cdot\wedge a_{p})
$$
where $h\in G_{0}$, $a_{1},..,a_{p}\in\gg_{\ol{1}}$.
\end{pf}
We want to understand how the representations (\ref{representation1}) and (\ref{representation2}) are translated with respect to the canonical isomorphism $\mathcal{A}_{G}\cong\cC^{\infty}_{G_{0}}\otimes\Lambda(\gg_{\ou}^{*})$; the results have a geometric role in some applications discussed in section \ref{sezione5}. An explicit formula for left-invariant vector fields is recalled (\cite{DP, Kz}). An analogous formula for the right invariant ones, which to the best of our knowledge has not been described in the literature, is obtained using the results of \cite{DP, P}. The restrictions of the representations (\ref{representation1}) and (\ref{representation2}) to $\gg_{\0}$ are compared to "trivial extensions" in subsection \ref{evensymmetries} while the restrictions to $\gg_{\ou}$ are more complicated to handle and they are studied in subsection \ref{oddsymmetries}.
\\
{\bf Notation}: The representation of $\gg_{\ol{0}}$ as right-invariant vector fields of the Lie group $G_{0}$ is denoted by 
\be
\label{firenze}
\hat{\varphi}_{0}:\gg_{\ol{0}}\rightarrow \mathcal{T}(G_{0})
\ee
\subsection{Even symmetries}\hspace{0 cm}\newline\\
\label{evensymmetries}
The algebra isomorphism (\ref{koskos}) can be used to construct a canonical embedding of Lie superalgebras
\be
\label{embed2}
\mathcal{T}_{G_{0}}(U)\hookrightarrow\Der_{\Lambda(\gg_{\ol{1}}^{*})}(\mathcal{C}^{\infty}(U)\otimes\Lambda(\gg_{\ol{1}}^{*}))\subseteq\Der_{\bR}(\mathcal{C}^{\infty}(U)\otimes\Lambda(\gg_{\ol{1}}^{*}))\cong\mathcal{T}_{G}(U)
\ee
Given a vector field $X\in\mathcal{T}_{G_{0}}(U)$, the corresponding vector field $X\in\mathcal{T}_{G}(U)$ (denoted by the same symbol by abuse of notation) acts trivially on $\gg_{\ol{1}}^{*}\subseteq\mathcal{C}^{\infty}(M_{0})\otimes\Lambda(\gg_{\ol{1}}^{*})$.
Composing (\ref{embed2}) with (\ref{omomorfismo}) and (\ref{firenze}) construct two representations
$$
\varphi_{0}:\gg_{\ol{0}}\hookrightarrow\Der_{\Lambda(\gg_{\ol{1}}^{*})}(\mathcal{C}^{\infty}(U)\otimes\Lambda(\gg_{\ol{1}}^{*}))\quad,\quad\hat{\varphi}_{0}:\gg_{\ol{0}}\hookrightarrow\Der_{\Lambda(\gg_{\ol{1}}^{*})}(\mathcal{C}^{\infty}(U)\otimes\Lambda(\gg_{\ol{1}}^{*}))\quad.
$$
It is not difficult to see that 
$$
\hat{\varphi}|_{\gg_{\ol{0}}}=\hat{\varphi}_{0}:\gg_{\ol{0}}\rightarrow\mathcal{T}(G)\qquad,
$$ 
while, 
for every $a\in\gg_{\ol{0}}$
$$
\varphi(a)|_{\mathcal{C}^{\infty}}=\varphi_{0}(a)\quad,\quad
\varphi(a)(\gg_{\ol{1}}^{*})=0\Longleftrightarrow \varphi(a)=\varphi_{0}(a)  \Longleftrightarrow [a,\gg_{\ol{1}}]=0\quad.
$$
This behaviour is described in more detail in the following two Propositions.
\begin{proposition}
\label{parikill}
The representation $\hat{\varphi}|_{\gg_{\0}}:\gg_{\0}\rightarrow\Der_{\bR}(\mathcal{C}^{\infty}(U)\otimes\Lambda(\gg_{\ou}^{*}))$ satisfies
$$
\hat{\varphi}(a)=\hat{\varphi}_{0}(a)
$$
for every $a\in\gg_{\0}$. In particular it preserves the $\bZ$-gradation of $\mathcal{C}^{\infty}(U)\otimes\Lambda(\gg_{\ou}^{*})$. 
\end{proposition}
\begin{pf}
Straightforward consequence of (\ref{representation2}) and the definition of the sheaf $\cA_{G}$.
\end{pf}
\begin{proposition}[\cite{Kz}]
\label{koszulpari}
The difference between $\varphi|_{\gg_{\ol{0}}}:\gg_{\ol{0}}\rightarrow\Der_{\bR}(\mathcal{C}^{\infty}(U)\otimes\Lambda(\gg_{\ol{1}}^{*}))$ and $\varphi_{0}:\gg_{\ol{0}}\hookrightarrow\Der_{\Lambda(\gg_{\ol{1}}^{*})}(\mathcal{C}^{\infty}(U)\otimes\Lambda(\gg_{\ol{1}}^{*}))$ is given by the formula 
$$
[\varphi(a)(f)-\varphi_{0}(a)(f)](a_{1}\wedge\cdot\cdot\wedge a_{p})=-\sum_{i=1}^{p}f(a_{1}\wedge\cdot\cdot\wedge [a,a_{i}]\wedge\cdot\cdot\wedge a_{p})
$$
where $a\in\gg_{\ol{0}}$, $a_{1},..,a_{p}\in\gg_{\ol{1}}$ and $f\in\Hom_{\bR}(\Lambda(\gg_{\ol{1}}),\mathcal{C}^{\infty}(U))$. The representation $\varphi|_{\gg_{\ol{0}}}:\gg_{\ol{0}}\rightarrow\Der_{\bR}(\mathcal{C}^{\infty}(U)\otimes\Lambda(\gg_{\ol{1}}^{*}))$ preserves the $\bZ$-gradation of $\mathcal{C}^{\infty}(U)\otimes\Lambda(\gg_{\ou}^{*})$. 
\end{proposition}
The following example clarifies the situation. 
\begin{example}
\label{examplepoincarè}
\rm{We reconstruct the Poincare' Lie supergroup $G$ associated with the sHC pair $(G_{0},\gg)$ of Example \ref{esempiopoincare} in signature $(1,2)$. Recall that
$$
\Spin^{0}_{1,2}\cong \SL(2,\mathrm{\bR})
$$
and that the spin module $S$ is a $2$-dimensional real vector space $S=[[s_{0},s_{1}]]$. The vector representation of $\SL(2,\bR)$ on $\bR^{1,2}=[[e_{0},e_{1},e_{2}]]$ is the adjoint representation via the identification of vector spaces endowed with a quadratic form
$$(\bR^{1,2},\left\langle \cdot,\cdot\right\rangle)\cong(\sl(2,\bR),\det)$$
given by
$$
e_{0}\cong\begin{pmatrix} 0 & 1 \\ -1 & 0 \end{pmatrix}\quad e_{1}\cong\begin{pmatrix} 0 & 1 \\ 1 & 0 \end{pmatrix}\quad e_{2}\cong\begin{pmatrix} 1 & 0 \\ 0 & -1 \end{pmatrix}\qquad.
$$
An element of
$
G_{0}=\SL(2,\bR)\ltimes\bR^{1,2}
$
is denoted by $(g,v)\in\SL(2,\bR)\ltimes\bR^{1,2}$.
Fix an $\sl(2,\bR)$-invariant symmetric bilinear map
$
\Gamma:S\vee S\rightarrow \bR^{1,2}
$
and denote by $\gg=\gg(\Gamma)=\gg_{\ol{0}}+\gg_{\ol{1}}=(\sl(2,\bR)\inplus\bR^{1,2})+S$ the Poincare' Lie superalgebra with
$$
[s_{\alpha},s_{\beta}]:=\Gamma_{\alpha\beta}^{k}e_{k}
$$
for $0\leq \alpha,\beta\leq 1$, $0\leq k\leq 2$, where we have used the Einstein convention, which is used throughout this example. 
An element of $\gg$ is denoted by $(A,w,s)\in(\sl(2,\bR)\inplus\bR^{1,2})+S$ where
$$
A=\begin{pmatrix} a_{0}^{0} & a_{1}^{0} \\ a_{0}^{1} & a_{1}^{1} \end{pmatrix}\in\sl(2,\bR)\quad.
$$ 
For every $0\leq i,j\leq 1$, denote by
\begin{displaymath}
\begin{array}{cc}
y^{i}_{j}:\SL(2,\bR)\rightarrow\bR \phantom{ccccccc}&\phantom{cccccc} x^{k}:\bR^{1,2}\rightarrow\bR \\
g=\begin{pmatrix} g_{0}^{0} & g_{1}^{0} \\ g_{0}^{1} & g_{1}^{1} \end{pmatrix}\rightarrow g_{j}^{i} \phantom{ccccccccc}&\phantom{cccccccc} v=v^{l}e_{l}\rightarrow v^{k}
\end{array}
\end{displaymath}
the coordinates of $\SL(2,\bR)$ and $\bR^{1,2}$. The following linear maps (recall (\ref{funzioni}))
\begin{align*}
& y^{i}_{j}:\Lambda(S)\rightarrow\mathcal{C}^{\infty}(\SL(2,\bR))\qquad\quad;\quad x^{k}:\Lambda(S)\rightarrow\mathcal{C}^{\infty}(\bR^{1,2})\qquad\quad;\quad s^{\alpha}:\Lambda(S)\rightarrow\bR\\\
&\phantom{cccccccc}1\rightarrow y^{i}_{j}\qquad\qquad\quad\qquad\quad\quad\phantom{ccccccccc}1\rightarrow x^{k}\phantom{cccccccccccccccccccccccc} s_{\beta}\rightarrow \delta^{\alpha}_{\beta}\\\
&\phantom{c}\Lambda^{l\neq 0}(S)\rightarrow 0\qquad\qquad\quad\qquad\quad\phantom{ccccc}\Lambda^{l\neq 0}(S)\rightarrow 0\phantom{cccccccccccccccccccc}\Lambda^{l\neq 1}(S)\rightarrow 0
\end{align*}
are coordinates of $G$. The $\mathcal{U}(\gg_{\ol{0}})$-linear maps from $\mathcal{U}(\gg_{\ol{0}})\otimes\Lambda(S)$ to $\mathcal{C}^{\infty}(G_{0})$ are denoted by the same letter $y^{i}_{j}$, $x^{k}$, $s^{\alpha}$, while the associated superfunctions by $$Y_{j}^{i}, X^{k}, S^{\alpha}\in\mathcal{A}(G)=\Hom_{\mathcal{U}(\gg_{\ol{0}})}(\mathcal{U}(\gg),\mathcal{C}^{\infty}(G_{0}))\quad.$$
Recall that $y_{j}^{i}=Y_{j}^{i}\circ \underline{\gamma}$, $x^{k}=X^{k}\circ\underline{\gamma}$, $s^{\alpha}=S^{\alpha}\circ\underline{\gamma}$. 
The representation $\varphi:\gg\rightarrow \mathcal{T}(G)$ of the Lie superalgebra $\gg=\sl(2,\bR)\inplus\bR^{1,2}+S$ as left-invariant vector fields of $G$ is 
$$
\varphi(A)=\varphi_{0}(A)-a^{\alpha}_{\beta}s^{\beta}\frac{\partial}{\partial s^{\alpha}}\quad,\quad
\varphi(w)=\varphi_{0}(w)
$$
$$
\varphi(s_{\alpha})=-\frac{\partial}{\partial s^{\alpha}}-\frac{1}{2}s^{\eta}\varphi_{0}(\Gamma_{\alpha\eta}^{k}e_{k})
$$
while that as right-invariant vector fields $\hat{\varphi}:\gg\rightarrow \mathcal{T}(G)$ is
$$
\hat{\varphi}(A)=\hat{\varphi}_{0}(A)\quad,\quad
\hat{\varphi}(w)=\hat{\varphi}_{0}(w)\quad,\quad
\hat{\varphi}(s_{\beta})=y^{\alpha}_{\beta}(g^{-1})[-\frac{\partial}{\partial s^{\alpha}}+
\frac{1}{2}s^{\eta}\varphi_{0}(\Gamma_{\alpha\eta}^{k}e_{k})] 
$$
where $0\leq\eta\leq 1$.
We calculate the comultiplication and the coinverse for the Lie subsupergroup of $G$ given by the sHC pair $(\bR^{1,2},\bR^{1,2}+S)$. By definition
$$
m^{*}(x^{k}):=m^{*}(X^{k})\circ(\gamma\otimes\gamma)\quad,
$$
\textit{i.e.}
$$
m^{*}(x^{k})(r\otimes t)(v_{0},v_{1})=X^{k}((\Ad_{-v_{1}}\gamma(r))\cdot \gamma(t))(v_{0}+v_{1})=
$$
$$
X^{k}(\gamma(r)\cdot\gamma(t))(v_{0}+v_{1})=
x^{k}(\underline{\gamma}^{-1}(\gamma(r)\cdot\gamma(t)))(v_{0}+v_{1})
$$
where $r,t\in\Lambda(S)$ and $v_{0},v_{1}\in \bR^{1,2}$. The computation of $\underline{\gamma}^{-1}(\gamma(r)\cdot\gamma(t))$ implies that
$$
m^{*}(x^{k})=m_{0}^{*}(x^{k})\otimes 1+\frac{1}{2}\sum_{0\leq \alpha,\beta\leq 1}\Gamma_{\alpha\beta}^{k}\otimes(s_{\alpha}\otimes s_{\beta})^{*}\in\mathcal{C}^{\infty}(\bR^{1,2}+\bR^{1,2})\otimes(\Lambda(S)\otimes\Lambda(S))^{*}\quad.
$$
Denoting by $\vartheta^{\alpha}\in\mathcal{C}^{\infty}(\bR^{1,2})\otimes\Lambda(S^{*})$ and $\theta^{\beta}\in\mathcal{C}^{\infty}(\bR^{1,2})\otimes\Lambda(S^{*})$ global odd coordinates on two copies of the Lie supergroup $(\bR^{1,2},\bR^{1,2}+S)$ we get the common forms of comultiplication and coinverse in the mathematical physics literature
$$
m^{*}(x^{k})=m_{0}^{*}(x^{k})-\frac{1}{2}
\Gamma_{\alpha\beta}^{k}\vartheta^{\alpha}\theta^{\beta}\qquad,\qquad
m^{*}(s^{\alpha})=\vartheta^{\alpha}+\theta^{\alpha}
$$
$$
i^{*}(x^{k})=
-x^{k}\qquad\qquad\qquad\quad\phantom{cccc},\qquad i^{*}(s^{\alpha})=-s^{\alpha}\quad.$$
}\end{example} 
\subsection{Odd symmetries}
\label{oddsymmetries}
\subsubsection{Introduction}\hspace{0 cm}\newline\\
\cite{DP, P} give a proof of Koszul's Theorem \ref{Koszulsx}. The arguments involved in the proof are briefly recalled and then Theorem \ref{Andrea} proves an analogous formula. 
Recall that the symmetric algebra $S(\gg)$ of $\gg$ has a natural structure of a $\bZ$-graded cocommutative Hopf superalgebra $(S(\gg),m_{S(\gg)},1_{S(\gg)},\Delta_{S(\gg)},\epsilon_{S(\gg)},\delta_{S(\gg)})$.
\begin{definition}[\cite{DP, P}]
\rm{Let $(A,m_{A},1_{A})$ be a graded algebra. The convolution algebra $\Hom_{\bR}(S(\gg),A)$ is the space of {\it formal functions} on $\gg$ with values in $A$. The algebra $A$ is embedded as the subalgebra of $\Hom_{\bR}(S(\gg),A)$ which consists of {\it constant functions}
$$
a:1\mapsto a
$$
$$
\phantom{ccccccccccccccccc}S^{n}(\gg)\mapsto 0\qquad\qquad {\rm if}\quad n\gneq 0
$$
where $a\in A$. 
}\end{definition}
An analogous definition of constant formal vector fields can be given for the convolution Lie superalgebras $\gg_{x}:=\Hom_{\bR}(S(\gg),\gg)$ and $\gg_{y}:=\Hom_{\bR}(\L(\gg_{\ou}),\gg)$.
\begin{definition}[\cite{DP, P}]
\rm{
The map $x\in\gg_{x}$ defined by
$$
S^{1}(\gg)=\gg\stackrel{\Id_{\gg}}{\rightarrow}\gg\quad,\qquad S^{n}(\gg)\rightarrow 0\quad (n\neq 1)
$$
is the {\it generic point} of $\gg$. Its restriction $y\in\gg_{y}$ to $\Lambda(\gg_{\ou})$
$$
\L^{1}(\gg_{\ou})=\gg_{\ou}\hookrightarrow\gg\quad,\quad \L^{n}(\gg_{\ou})\rightarrow 0\quad (n\neq 1)
$$
is the {\it generic point} of $\gg_{\ou}$.
}\end{definition}
For every $n\in\bN$ the map $(\ad x)^{n}:\gg_{x}\rightarrow\gg_{x}$ is a $S(\gg)^{*}$-linear morphism of $\gg_{x}$. In particular, for every $a\in\gg\subseteq \gg_{x}$ the map $(\ad x)^{n}(a):S(\gg)\rightarrow\gg$ is given by
\be
\label{generic}
S^{p}(\gg)\ni a_{1}\cdot\cdot\cdot a_{p}\mapsto
\begin{cases}
0 & p\neq n\\
\pm\sum_{s\in \sum_{p}}\ad a_{s(1)}\circ\cdot\cdot\circ \ad a_{s(p)}(a) & p=n
\end{cases}
\ee
where $\pm$ is given by the rule of signs in supergeometry. An analogous formula holds for the generic point of $\gg_{\ou}$.
A formal power series in the even variable t is denoted by $$f=\sum_{i=0}^{+\infty}f_{i}t^{i}\in\bR[[t]]\quad.$$ 
The following definition makes sense due to equation (\ref{generic}).
\begin{definition}[\cite{DP, P}]
\rm{We denote by $f(\ad x)(a)\in\gg_{x}$ the formal vector field
$$
f(\ad x)(a):=f_{0}a+f_{1}(\ad x)(a)+f_{2}(\ad x)^{2}(a)+\cdot\cdot\cdot\cdot
$$
where $a\in\gg$. The associated coderivations of $S(\gg)$ are denoted by
$$
c_{f(\ad x)(a)}:=\Id_{S(\gg)}*\phantom{c}f(\ad x)(a)\quad\phantom{c},\quad
\phantom{a}_{f(\ad x)(a)}c:=f(\ad x)(a)*\Id_{S(\gg)}\quad.
$$
}\end{definition}
\cite{P} investigates for which power series $f\in\bR[[t]]$ the set made up of the maps
$$
\phantom{cccc}\,\gg\rightarrow \CoDer(S(\gg))
$$
\be
\label{universal}
a\mapsto c_{f(\ad x)(a)}
\ee 
is a universal (in a functorial sense, see \cite{P} for the precise meaning of this assertion) representation by coderivations in the category of Lie superalgebras over $\bR$. 
\begin{definition}[\cite{DP, P}]
\rm{Let $c\in\bR$ be a real number. The power series $f_{c}\in\bR[[t]]$ is defined as  
$$
\phantom{ccc}f_{c}(t):=\frac{t}{e^{\frac{t}{c}}-1}\in\bR[[t]]
$$
when $c\in\bR-\left\{0\right\}$ and
$$
f_{c}(t)=f_{0}(t):=-t\in\bR[[t]]
$$
when $c=0\in\bR$.
}\end{definition}
Any power series $f_{c}\in\bR[[t]]$ gives rise to a universal representation (\ref{universal}). Denote by 
$$
\phantom{c}\Phi_{c}:\mathcal{U}(\gg)\rightarrow \End_{\bR}(S(\gg))\quad,\quad
u\mapsto\Phi_{c}^{u}
$$
the unique extension of this representation to an algebra morphism from $\cU(\gg)$. Denote by the same symbol the associated map
$\Phi_{c}:\cU(\gg)\otimes S(\gg)\rightarrow S(\gg)$.
\begin{definition}[\cite{DP, P}]
\rm{
The map 
$$
\Phi_{c}(1):\cU(\gg)\rightarrow S(\gg)\quad,\quad
u\mapsto\Phi_{c}^{u}(1)
$$
is called the {\it symbol map} associated with the power series $f_{c}\in\bR[[t]]$.
}\end{definition}
\begin{theorem}[\cite{P}]
\label{petracciona}
The inverse of the symbol map $\Phi_{1}(1)$ is the symmetrization $\gamma:S(\gg)\rightarrow \cU(\gg)$
$$
\gamma(a_{1}\cdot\cdot\cdot a_{p})=\pm\frac{1}{p!}\sum_{s\in \sum_{p}}a_{s(1)}\cdot\cdot\cdot a_{s(p)}\qquad\qquad a_{1}\cdot\cdot\cdot a_{p}\in S^{p}(\gg)
$$
where $\pm$ is given by the rule of signs in supergeometry. The map $\gamma$ commutes with endomorphisms and derivations of $\gg$ and it is an isomorphism of coalgebras.
\end{theorem}
\begin{corollary}[\cite{P}]
\label{right1}
For all $a\in\gg$
\begin{itemize}
\item[1)] $\gamma^{-1}\circ \ad(a) \circ \gamma=\Phi_{0}^{a}$
\item[2)] $\gamma^{-1}\circ L_{a} \circ \gamma=\Phi_{1}^{a}$
\item[3)] $\gamma^{-1}\circ R_{a} \circ \gamma=-\Phi_{-1}^{a}$
\end{itemize}
where $\ad(a)$, $L_{a}$, $R_{a}\in\End_{\bR}(\cU(\gg))$.
\end{corollary}
\begin{pf}\\
1) The derivations $\ad(a)$ of $\cU(\gg)$ and $\Phi_{0}^{a}$ of $S(\gg)$ are extensions of the derivation $\ad(a)\in\Der(\gg)$ of $\gg$. \\
2) For every $b\in S(\gg)$, $\gamma^{-1}(a\cdot\gamma(b))=\Phi_{1}^{a\cdot\gamma(b)}(1)=\Phi_{1}^{a}\circ\Phi_{1}^{\gamma(b)}(1)=\Phi_{1}^{a}(b)$.\\
3) The previous cases give $\gamma^{-1}\circ R_{a}\circ\gamma=\Phi_{1}^{a}-\Phi_{0}^{a}=-\Phi_{-1}^{a}$.
\end{pf}
It follows that (\ref{symmetrization2})
is a coalgebra isomorphism, where $\mathcal{U}(\gg_{\0})\otimes \Lambda(\gg_{\ou})$ has the tensor product Hopf-algebra structure. 
\\
For every $c\in\bR-\left\{0\right\}$, denote the formal power series of $t\cdot coth(\frac{t}{c})$ and $th(\frac{t}{2c})$ by 
$$
\phantom{ccc}p_{c}(t):=t\cdot coth(\frac{t}{c})=c+\sum_{n=1}^{+\infty}b_{2n}\frac{2^{2n}t^{2n}}{c^{2n-1}(2n)!}\in\bR[[t]]
$$
$$
q_{c}(t):=-th(\frac{t}{2c})=-\sum_{n=1}^{+\infty}b_{2n}\frac{(2^{2n+1}-2)t^{2n-1}}{c^{2n-1}(2n)!}\in\bR[[t]]
$$
where $\left\{b_{p}\right\}_{p\in\bN}$ is the sequence of Bernoulli numbers. For $a\in\gg_{\ou}$, denote by
$$
\alpha^{a}_{c}:=p_{c}(\ad y)(a)\in \Hom_{\bR}(\Lambda(\gg_{\ou}),\gg_{\ou})_{\ou}
$$
and
$$
\theta^{a}_{c}:=q_{c}(\ad y)(a)\in \Hom_{\bR}(\Lambda(\gg_{\ou}),\gg_{\0})_{\ou}
$$
the formal vector fields given by
$$
\alpha^{a}_{1}:a_{1}\wedge\cdot\cdot\cdot\wedge a_{2p}\mapsto\sum_{s\in\sum_{2p}}\frac{b_{2p}2^{2p}}{(2p)!}\sgn(s)\ad a_{s(1)}\circ\cdot\cdot\circ \ad a_{s(2p)}(a)
$$
and
$$
\theta^{a}_{1}:a_{1}\wedge\cdot\cdot\cdot\wedge a_{2p-1}\mapsto\sum_{s\in\sum_{2p-1}}\frac{b_{2p}(2^{2p+1}-2)}{(2p)!}\sgn(s)\ad a_{s(1)}\circ\cdot\cdot\circ \ad a_{s(2p-1)}(a)\quad.
$$
\subsubsection{Koszul's theorems}\hspace{0 cm}\newline\\
\cite{DP} proves the following
\begin{theorem}[\cite{DP, Kz}]
\label{teoremone2}
For every $a\in\gg_{\ou}$, the coderivation of $\cU(\gg_{\0})\otimes \Lambda(\gg_{\ou})$ induced by $R_{a}\in\CoDer_{\bR}(\mathcal{U}(g))$ via the isomorphism (\ref{symmetrization2}) is given by
$$
\underline{\gamma}^{-1}\circ R_{a}\circ\underline{\gamma}=
$$
$$
[-(m_{\cU(\gg_{\0})}\otimes \Id)\circ(\Id\otimes \theta^{a}_{1}\otimes \Id)+(\Id\otimes m_{\L(\gg_{\ou})})\circ(\Id\otimes \alpha^{a}_{1}\otimes \Id)]\circ (\Id_{\cU(\gg_{\0})}\otimes\Delta_{\Lambda(\gg_{\ou})})
$$
\end{theorem}
The aim is to prove an analogous result for the coderivation $L_{a}\in\CoDer_{\bR}(\cU(\gg))$. First, we prove the following structural result.
\begin{proposition}
\label{shape}
For every $a\in\gg_{\ou}$, the coderivation of $\cU(\gg_{\0})\otimes\Lambda(\gg_{\ou})$ induced by $L_{a}\in\CoDer(\cU(\gg))$ via the isomorphism (\ref{symmetrization2}) satisfies
\be
\label{shape2}
\underline{\gamma}^{-1}\circ L_{a}\circ\underline{\gamma}_{|\Lambda(\gg_{\ou})}=(\theta^{a}\otimes\Id)\circ\Delta_{\Lambda(\gg_{\ou})}+
m_{\Lambda(\gg_{\ou})}\circ 
(\epsilon^{a}\otimes\Id)\circ\Delta_{\Lambda(\gg_{\ou})}
\ee
for some $\theta^{a}\in \Hom_{\bR}(\L(\gg_{\ou}),\gg_{\0})$ and $\epsilon^{a}\in \Hom_{\bR}(\L(\gg_{\ou}),\gg_{\ou})$. 
\end{proposition}
\begin{pf}
For every $a\in\gg_{\ou}$, Theorem \ref{R} applied to the coderivation $\underline{\gamma}^{-1}\circ L_{a}\circ\underline{\gamma}\in\CoDer(\mathcal{U}(\gg_{\0})\otimes \Lambda(\gg_{\ou}))$ implies that
$$
\underline{\gamma}^{-1}\circ L_{a}\circ\underline{\gamma}=m_{\cU(\gg_{\0})\otimes \Lambda(\gg_{\ou})}\circ(\theta^{a}\otimes\Id+\epsilon^{a}\otimes\Id)\circ\Delta_{\cU(\gg_{\0})\otimes \Lambda(\gg_{\ou})}
$$
for some $\epsilon^{a}\in \Hom_{\bR}(\cU(\gg_{\0})\otimes \Lambda(\gg_{\ou}),\gg_{\ou})$ and $\theta^{a}\in \Hom_{\bR}(\cU(\gg_{\0})\otimes \Lambda(\gg_{\ou}),\gg_{\0})$. It is then sufficient to restrict the equality to $\Lambda(\gg_{\ou})$.
\end{pf}
Now we determine $\theta^{a}\in \Hom_{\bR}(\L(\gg_{\ou}),\gg_{\0})$ and $\epsilon^{a}\in \Hom_{\bR}(\L(\gg_{\ou}),\gg_{\ou})$.
A naive idea could be to obtain results analogous to part 1) of Corollary \ref{right1}. Unfortunately, the map (\ref{symmetrization2}) {\bf does not} commute with the odd derivations of $\gg$, and in particular with the adjoint action of an element $a\in\gg_{\ou}$. It turns out that handling the antipodes of $\mathcal{U}(\gg)$, $\mathcal{U}(\gg_{\0}$) and $\Lambda(\gg_{\ou})$ is useful to obtain what we are looking for.
\begin{lemma}
\label{andreone2}
\label{andreone}
For every $a\in\gg_{\ou}$
$$-\delta_{\cU(\gg)}\circ R_{a}\circ \delta_{\cU(\gg)}=L_{a}\quad.$$
The map induced by the antipode $\delta_{\cU(\gg)}$ of $\cU(\gg)$ via the isomorphism (\ref{symmetrization2}) is given by
$$
\delta_{\cU(\gg)}\circ\underline{\gamma}=\underline{\gamma}\circ(\Id\otimes\Phi_{0})\circ(\delta_{\cU(\gg_{\0})}\otimes\Id \otimes \delta_{\L(\gg_{\ou})})\circ (\Delta_{\cU(\gg_{\0})}\otimes\Id)
$$
\end{lemma}
\begin{pf}
The proof of the first assertion is straightforward.
The proof of the second assertion consists of repeated applications of 
$$
\gamma(a_{1}\cdot\cdot\cdot a_{p})u_{i}=u_{i}\gamma(a_{1}\cdot\cdot\cdot a_{p})-\gamma(\Phi_{0}^{u_{i}}(a_{1}\cdot\cdot\cdot a_{p}))
$$
to the equation
$$
\delta_{\cU(\gg)}\circ\underline{\gamma}(u_{1}\cdot\cdot\cdot u_{n}\otimes a_{1}\cdot\cdot\cdot a_{p})=(-1)^{p+n}\gamma(a_{1}\cdot\cdot\cdot a_{p})u_{n}\cdot\cdot\cdot u_{1}
$$
where $u_{i}\in\gg_{\0}$ for $i=1,..,n$ and $a_{j}\in\gg_{\ou}$ for $j=1,...,p$.
\end{pf}
Denote the formal power series of $p_{1}(t)+tq_{1}(t)=t\cdot(coth(t)-th(\frac{t}{2}))$ by
$$
e(t):=t\cdot(coth(t)-th(\frac{t}{2}))=\sum_{n=0}^{+\infty}b_{2n}\frac{(-2^{2n+1}+2^{2n}+2)t^{2n}}{(2n)!}\in\bR[[t]]\quad.
$$
For $a\in\gg_{\ou}$, denote by
$$
\epsilon^{a}:=e(\ad y)(a)\in \Hom_{\bR}(\Lambda(\gg_{\ou}),\gg_{\ou})_{\ou}
$$
the formal vector field given by
\be
\label{macheneso}
\epsilon^{a}:a_{1}\wedge\cdot\cdot\cdot\wedge a_{2p}\mapsto\sum_{\sigma\in\sum_{2p}}\frac{b_{2p}(-2^{2p+1}+2^{2p}+2)}{(2p)!}\sgn(\sigma)\ad a_{\sigma(1)}\circ\cdot\cdot\circ \ad a_{\sigma(2p)}(a)
\ee
\begin{theorem}
\label{Teoremone}
For every $a\in\gg_{\ou}$, the restriction to $\L(\gg_{\ou})$ of the coderivation of $\cU(\gg_{\0})\otimes \Lambda(\gg_{\ou})$ induced by $L_{a}\in\CoDer(\cU(\gg))$ via the isomorphism (\ref{symmetrization2}) is given by
$$
\underline{\gamma}^{-1}\circ L_{a}\circ\underline{\gamma}_{|\Lambda(\gg_{\ou})}=(\theta_{1}^{a}\otimes\Id)\circ\Delta_{\Lambda(\gg_{\ou})}+
m_{\Lambda(\gg_{\ou})}\circ 
(\epsilon^{a}\otimes\Id)\circ\Delta_{\Lambda(\gg_{\ou})}
$$
where $\epsilon^{a}\in\Hom_{\bR}(\Lambda(\gg_{\ou}),\gg_{\ou})_{\ou}$ is given by (\ref{macheneso}).
\end{theorem}
\begin{example}
For every $a,a_{1},a_{2}\in\gg_{\ou}$
$$
a\gamma(a_{1}\cdot a_{2})=-\frac{1}{6}[a_{1},[a_{2},a]]+\frac{1}{6}[a_{2},[a_{1},a]]+\gamma(a\cdot a_{1}\cdot a_{2})+\frac{1}{2}[a_{1},a]a_{2}-\frac{1}{2}[a_{2},a]a_{1}
\quad.$$
\end{example}
\begin{pf}
From Theorem \ref{teoremone2} and Lemma \ref{andreone2}, the l.h.s. of (\ref{shape2}) 
equals
$$
-(\Id\otimes\Phi_{0})\circ(\delta_{\cU(\gg_{\0})}\otimes \Id \otimes \delta_{\L(\gg_{\ou})})\circ (\Delta_{\cU(\gg_{\0})}\otimes\Id)\circ
$$
$$
[m_{\L(\gg_{\ou})}\circ(\alpha^{a}_{1}\otimes \Id)-(\theta^{a}_{1}\otimes \Id)]\circ \Delta_{\L(\gg_{\ou})}\circ\delta_{\L(\gg_{\ou})}\quad.
$$
Using the sigma notation of \cite{Sw}, the previous equation can be re-written as
$$
v\mapsto-(-1)^{|v_{(1)}|}\theta_{1}^{a}v_{(1)}\otimes v_{(2)}+(-1)^{|v_{(1)}|}1\otimes\Phi_{0}^{\theta_{1}^{a}v_{(1)}}v_{(2)}+(-1)^{|v_{(1)}|}1\otimes\alpha_{1}^{a}v_{(1)}\cdot v_{(2)}
$$
$$
=\theta_{1}^{a}v_{(1)}\otimes v_{(2)}-1\otimes\Phi_{0}^{\theta_{1}^{a}v_{(1)}}v_{(2)}+1\otimes\alpha_{1}^{a}v_{(1)}\cdot v_{(2)}
$$
where $v\in\Lambda(\gg_{\ou})$. The aim is to re-write the second term. Proposition \ref{shape} states that there exists $\lambda^{a}\in\Hom_{\bR}(\L(\gg_{\ou}),\gg_{\ou})$ such that
$$
-\Phi_{0}^{\theta_{1}^{a}v_{(1)}}v_{(2)}=m_{\L(\gg_{\ou})}\circ(\lambda^{a}\otimes \Id)\circ \Delta_{\L(\gg_{\ou})}\quad.
$$
The fact that the antipode $\delta_{\L(\gg_{\ou})}=\delta_{S(\gg)}|_{\Lambda(\gg_{\ou})}$ of $\L(\gg_{\ou})$ is the convolution inverse of $\Id_{\L(\gg_{\ou})}$ implies that $\lambda^{a}\in\Hom_{\bR}(\L(\gg_{\ou}),\gg_{\ou})$ equals  
$$
\Lambda(\gg_{\ou})\ni v\mapsto v_{(1)}\otimes v_{(2)}\mapsto-\Phi_{0}^{\theta_{1}^{a}v_{(1)}}v_{(2)}\cdot\delta_{S(\gg)}(v_{(3)})=
$$
$$
\phantom{cccccccc}t(\ad y)(\theta_{1}^{a}v_{(1)})(v_{(2)})
=-t\cdot th(\frac{t}{2})(\ad y)(a)(v)
$$ 
which is the assertion of the Theorem. The last equality is a consequence of the following equations.
Since $t(\ad y)(\theta_{1}^{a}v_{(1)})\in\gg_{y}$ acts non-trivially only on $\gg_{\ou}\subseteq \L(\gg_{\ou})$
\be
\label{sopra}
t(\ad y)(\theta_{1}^{a}v_{(1)})(v_{(2)})
=\sum_{j=1}^{2n}(-1)^{j}[v_{j},\theta_{1}^{a}(v_{1}\cdot\cdot \hat{v_{j}}\cdot\cdot v_{2n})]
\ee
where $v=v_{1}\cdot\cdot\cdot v_{2n}\in\L(\gg_{\ou})$, $v_{j}\in\gg_{\ou}$ for $j=1,...,2n$.
Moreover
$$
\theta_{1}^{a}(v_{1}\cdot\cdot \hat{v_{j}}\cdot\cdot v_{2n})=
\sum_{\sigma\in\sum_{2n-1}}\frac{b_{2n}(2^{2n+1}-2)}{(2n)!}\sgn(\sigma)\ad v_{\sigma(1)}\circ\cdot\cdot\circ\hat{\ad v_{\sigma(j)}}\circ\cdot\cdot\circ\ad v_{\sigma(2n)}(a)
$$
which implies that the last term of (\ref{sopra}) equals
$$
\sum_{j=1}^{2n}\sum_{\sigma\in\sum_{2n-1}}(-1)^{j}\frac{b_{2n}(2^{2n+1}-2)}{(2n)!}\sgn(\sigma)\ad v_{j}\circ\ad v_{\sigma(1)}\circ\cdot\cdot\circ\hat{\ad v_{\sigma(j)}}\circ\cdot\cdot\circ\ad v_{\sigma(2n)}(a)
$$
$
=-\sum_{\sigma\in\sum_{2n}}\frac{b_{2n}(2^{2n+1}-2)}{(2n)!}\sgn(\sigma)\ad v_{\sigma(1)}\circ\cdot\cdot\circ\ad v_{\sigma(2n)}(a)
$.
\end{pf}

\subsubsection{Left odd symmetries}\hspace{0 cm}\newline\\
The following Theorem of Koszul holds.
\begin{theorem}[\cite{DP, Kz, P}]
\label{Koszulsx}
For every $a\in\gg_{\ou}$ and $f\in \Hom_{\bR}(\Lambda(\gg_{\ou}),\mathcal{C}^{\infty}(U))$
\be
\label{koszulleft}
\varphi(a)(f)=-\theta_{1}^{a}* f+(-1)^{|f|}f\circ(\alpha_{1}^{a}* \Id_{\Lambda(g_{\ou})})
\ee
where the convolution product $\theta_{1}^{a}* f$ is defined via the map $\varphi_{0}:\gg_{\0}\rightarrow\mathcal{T}(G_{0})$.
\end{theorem}
\begin{pf}
Theorem \ref{teoremone2} implies the relation
$$
(-1)^{|v|}\gamma(v)a=\gamma(\alpha_{1}^{a}(v_{(1)})\cdot v_{(2)})-\theta_{1}^{a}(v_{(1)})\cdot \gamma(v_{(2)})
=\gamma\circ(\alpha_{1}^{a}*\Id_{\Lambda(\gg_{\ou})})-\theta_{1}^{a}*\gamma
$$
for $a\in\gg_{\ou}$, $v\in\Lambda(\gg_{\ou})$. It follows that
$$
(-1)^{|f|}\varphi(a)(f)(v)=(-1)^{|f|}\varphi(a)(F)(\gamma(v))=(-1)^{|f|+1}F(\gamma(v)a)
$$
$$
=(-1)^{|v|}F(\gamma(v)a)=[f\circ(\alpha_{1}^{a}*\Id_{\Lambda(\gg_{\ou})})-(-1)^{|f|}\theta_{1}^{a}*f](v)\quad.
$$
\end{pf}
For every $g\in G_{0}$  
\be
\label{tangentedispari}
\varphi(a)|_{g}=-\frac{\partial}{\partial a^{*}}|_{g}\cong a\in\gg_{\ou}
\ee
\textit{i.e.} $\ev_{\ou}(\varphi(a))$ equals the algebraic derivation given by contraction of $\gg_{\ou}$ and $\gg_{\ou}^{*}$. On the other hand, if $[\gg_{1},\gg_{1}]=0$, formula (\ref{koszulleft}) implies that 
$$\varphi(a)(f)(a_{1}\wedge\cdot\cdot\wedge a_{p})=(-1)^{|f|}f(a\wedge a_{1}\wedge\cdot\cdot\wedge a_{p})\qquad,$$
\textit{i.e.} that 
$$
\varphi(a)=-\frac{\partial}{\partial a^{*}}\cong a\in\gg_{\ou}
$$
is the algebraic derivation given by contraction of $\gg_{\ou}$ and $\gg_{\ou}^{*}$ {\bf globally}.
Informally, Proposition \ref{koszulpari} together with Theorem \ref{Koszulsx} show that the bracket $[\gg,\gg_{\ou}]$ measure how much the representation $\varphi:\gg\rightarrow\Der_{\bR}(\mathcal{C}^{\infty}(U)\otimes\Lambda(\gg_{\ou}^{*}))$ preserves the exterior algebra structure of (\ref{koskos}). 
\subsubsection{Right odd symmetries}
\begin{theorem}
\label{Andrea}
For every $a\in\gg_{\ou}$ and $f\in \Hom_{\bR}(\Lambda(\gg_{\ou}),\mathcal{C}^{\infty}(U))$
\be
\label{Andreadestra}
\hat{\varphi}(a)(f)=\theta_{1}^{\Ad_{g}^{-1}a}*f+(-1)^{|f|}f\circ(\epsilon^{\Ad_{g^{-1}}a}*\Id_{\Lambda(\gg_{\ou})})
\ee
where the convolution product $\theta_{1}^{\Ad_{g}^{-1}a}*f$ is defined via the map $\varphi_{0}:\gg_{\0}\rightarrow\mathcal{T}(G_{0})$.
\end{theorem}
\begin{pf}
Theorem \ref{Teoremone} implies the relation
$$
a\gamma(v)=\gamma(\epsilon^{a}(v_{(1)})\cdot v_{(2)})+\theta_{1}^{a}(v_{(1)})\cdot \gamma(v_{(2)})
=\gamma\circ(\epsilon^{a}*\Id_{\Lambda(\gg_{\ou})})+\theta_{1}^{a}*\gamma
$$
for $a\in\gg_{\ou}$, $v\in\Lambda(\gg_{\ou})$. It follows that
$$
(\hat{\varphi}(a)f)(v)(g)=(\hat{\varphi}(a)F)(\gamma(v))(g)=(-1)^{|f|}F((\Ad_{g^{-1}}a)\gamma(v))(g)
$$
$$
=[(-1)^{|f|}f\circ(\epsilon^{\Ad_{g^{-1}}a}*\Id_{\Lambda(\gg_{\ou})})+\theta_{1}^{\Ad_{g^{-1}}a}*f](v)(g)\quad.
$$
\end{pf}
Denote a basis of $\gg_{\ou}$ by $\left\{s_{\alpha}\right\}_{\alpha=1}^{\dim\gg_{\ou}}$ and its dual basis by $\left\{s^{\alpha}\right\}_{\alpha=1}^{\dim\gg_{\ou}}$. For every $g\in G_{0}$ and every natural number $1\leq\beta\leq\dim\gg_{\ou}$, the value of the right-invariant odd vector field $\hat{\varphi}(s_{\beta})$ is given by
\be
\label{andiamoeven}
\hat{\varphi}(s_{\beta})|_{g}=-\Ad^{\alpha}_{\beta}(g^{-1})\otimes\frac{\partial}{\partial s^{\alpha}}|_{g}\cong\Ad^{\alpha}_{\beta}(g^{-1})\otimes s_{\alpha}\cong\Ad_{g^{-1}}s_{\beta}
\ee 
where 
$$
\Ad_{g}|_{\gg_{\ou}}=\begin{pmatrix} \Ad_{1}^{1}(g) & \Ad_{2}^{1}(g) & \cdot\cdot & \Ad_{n}^{1}(g) \\ \Ad_{1}^{2}(g) & \Ad_{2}^{2}(g) & \cdot\cdot & \Ad_{n}^{2}(g) \\ \cdot & \cdot & \cdot\cdot & \cdot \\ \cdot & \cdot & \cdot\cdot & \cdot \\ \Ad_{1}^{n}(g) & \Ad_{2}^{n}(g) & \cdot\cdot & \Ad_{n}^{n}(g) \end{pmatrix}\in\End_{\bR}(\gg_{\ou})
$$
is the matrix representation of the adjoint action of $G_{0}$ on $\gg_{\ou}$. Note that $\ev_{\ou}(\hat{\varphi}(a))$ depends on the bracket $[\gg_{\0},\gg_{\ou}]$. Moreover, if $[\gg_{1},\gg_{1}]=0$, formula (\ref{Andreadestra}) implies that 
$$\hat{\varphi}(a)(f)(a_{1}\wedge\cdot\cdot\wedge a_{p})=(-1)^{|f|}f(\Ad_{g^{-1}}a\wedge a_{1}\wedge\cdot\cdot\wedge a_{p})\quad,$$
\textit{i.e.} that 
$$
\hat{\varphi}(a)=-\frac{\partial}{\partial (\Ad_{g^{-1}}a)^{*}}\cong \Ad_{g^{-1}}a\in\gg_{\ou}
$$
is the algebraic derivation given by contraction of $\gg_{\ou}$ and $\gg_{\ou}^{*}$ {\bf globally}.
\section{Homogeneous supermanifolds}
\label{section3}
\setcounter{equation}{0}
As in the classical theory, a homogeneous supermanifold can be understood either in terms of a transitive action of a Lie supergroup $G=(G_{0},\mathcal{A}_{G})$ on a supermanifold $M=(M_{0},\cA_{M})$ (see Definition \ref{transitive}) or by means of a quotient of $G=(G_{0},\cA_{G})$ with respect to a closed Lie subsupergroup $H=(H_{0},\cA_{H})$ (see Definition \ref{Fio}).
\subsection{Action of a Lie supergroup and of a super Harish-Chandra pair}
\begin{definition}[\cite{DM}]
\label{action}
\rm{An action of a Lie supergroup $G=(G_{0},\mathcal{A}_{G})$ on a supermanifold $M=(M_{0},\mathcal{A}_{M})$ is a morphism 
$
\rho=(\rho_{0},\rho^{*}):G\times M\rightarrow M
$
satisfying
\begin{itemize}
\item[$1)$]
$\rho\circ(id_{G}\times \rho)=\rho\circ(m\times id_{M}):G\times G\times M\longrightarrow M$ 
\item[$2)$]
$\rho\circ\left\langle \hat{e},id_{M}\right\rangle=id_{M}:M\longrightarrow M.$
\end{itemize}
The triple $(M_{0},\mathcal{A}_{M},\rho)$ is called a {\it $G$-supermanifold}.
}\end{definition}
\begin{definition}[\cite{DM}]
\label{actionSHCP}
\rm{An action of a sHC pair $(G_{0},\gg)$ on a supermanifold $M=(M_{0},\mathcal{A}_{M})$ is a pair $(\rho_{\ol{0}},\hat{\rho})$ consisting of a {\bf global} action, \textit{i.e.} a group homomorphism 
$$
\rho_{\ol{0}}:G_{0}\rightarrow\Aut(M)
$$ 
and a {\bf fundamental} action, \textit{i.e.} a Lie superalgebra anti-homomorphism
$$
\hat{\rho}:\gg\rightarrow\mathcal{T}(M)
$$ 
satisfying 
$$
\hat{\rho}(B)(f)=(\rho_{\0})_{*}(B)(f):=\frac{d}{dt}|_{t=0}(\rho_{\ol{0}}\circ\exp(tB))^{*}(f)\quad\quad f\in\mathcal{A}(M)
$$
for every $B\in\gg_{\ol{0}}$. The triple $(\mathcal{A}_{M},\rho_{\ol{0}},\hat{\rho})$ is called {\it a $(G_{0},\gg)$-supermanifold}.
}\end{definition}
If there is no danger of confusion, the global diffeomorphisms and the fundamental vector fields are denoted by $$\Aut(M)\ni\rho_{\ol{0}}(g):=g\qquad,\qquad \mathcal{T}(M)\ni\hat{\rho}(A):=\hat{A}$$
where $g\in G_{0}$ and $A\in\gg$.
\begin{theorem}[\cite{DM}]
Any action $(M_{0},\mathcal{A}_{M},\rho)$ of a Lie supergroup $G=(G_{0},\mathcal{A}_{G})$ defines an action $(\mathcal{A}_{M},\rho_{\ol{0}},\hat{\rho})$ of the corresponding sCH pair $(G_{0},\gg)$ where
$$
g^{*}:=(\ev_{g}\otimes \Id)\circ \rho^{*}:\mathcal{A}(M)\longrightarrow \mathcal{A}(M)
\qquad,\qquad
\hat{A}:=(A|_{e}\otimes \Id)\circ\rho^{*}$$
for all $g\in G_{0}$ and $A\in\gg$. The correspondence 
$
\rho\mapsto(\rho_{\ol{0}},\hat{\rho})
$
is a bijection between the sets of actions of a Lie supergroup and those of the associated sHC pair respectively.
\end{theorem}
\begin{example}
\rm{
Let $G=(G_{0},\gg)$ be a sHC pair and let $V=V_{\0}+V_{\ou}$ be a finite dimensional supervector space. A representation $\phi:G\rightarrow\GL_{\bR}(V)$ of $G$ on $V$ consists of
\begin{itemize}
\item[1)] a representation $\varphi_{\ol{0}}:G_{0}\rightarrow\GL_{\bR}(V)$ of $G_{0}$ on $V$,
\item[2)] a representation $\varphi:\gg\rightarrow\ggl_{\bR}(V)$ of $\gg$ on $V$,
\end{itemize}
such that $(\varphi_{\ol{0}})_{*}=\varphi|_{\gg_{\ol{0}}}:\gg_{\ol{0}}\rightarrow\ggl_{\bR}(V)_{\ol{0}}$. The representation (\ref{grazie}) of  Definition \ref{superharish} is called {\it adjoint representation}; for any closed Lie subsupergroup $(H_{0},\gh)\subseteq (G_{0},\gg)$ (\textit{i.e.} $H_{0}\subseteq G_{0}$ closed) it restricts to a well-defined map
$
\Ad:H\rightarrow \GL_{\bR}(\gg/\gh)
$.
}\end{example}
\begin{example}
\rm{
The adjoint action $a:G\times G\longrightarrow G$ of a Lie supergroup $G=(G_{0},\mathcal{A}_{G})$ on itself 
is defined as an action of the underlying $sHC$ pair $(G_{0},\gg)$:
$$
\rho_{\ol{0}}(g)^{*}:=a_{g}^{*}:=L_{g}^{*}\circ R_{g^{-1}}^{*}\qquad,\qquad
\hat{A}:=(A|_{e}\otimes \Id)\circ m^{*}-A
$$
for every $g\in G_{0}$ and $A\in\gg$.
}\end{example}
\begin{lemma}
\label{obvious}
Let $(M_{0},\mathcal{A}_{M},\rho)$ be a $G$-supermanifold, $A\in\gg$ 
and $g\in G_{0}$. Then
\begin{itemize}
\item[1)] $(A\otimes \Id)\circ\rho^{*}=(\Id\otimes\hat{A})\circ\rho^{*}\qquad,$
\item[2)] $(_{v}A\otimes \Id)\circ\rho^{*}=\rho^{*}\circ\hat{A}\qquad,$
\item[3)] $g_{*}\hat{A}=\hat{\Ad_{g}A}\qquad,$
\end{itemize}
where $v:=A|_{e}\in T_{e}(G)$.
\end{lemma}
\subsection{Orbit map and stability subgroup at a point 
}\hfill\newline\\
Every point $p\in M_{0}$ of a $G$-supermanifold $(M_{0},\mathcal{A}_{M},\rho)$ defines a morphism of supermanifolds
$$
\rho_{p}=((\rho_{p})_{0},\rho_{p}^{*}):=\rho\circ\left\langle id_{G},\hat{p}\right\rangle:G\longrightarrow M\qquad,
$$
called the {\it orbit map} of the point $p\in M_{0}$, which satisfies
$$
\rho_{p}\circ R_{g}=\rho_{g\cdot p}\qquad\qquad,\qquad\qquad g\circ\rho_{p}=\rho_{p}\circ L_{g}
$$
for every $g\in G_{0}$.
\begin{definition}[\cite{Go}]
\label{transitive}
\rm{The action $\rho$ is {\it transitive} if $\rho_{p}$ is a surjective submersion for every $p\in M_{0}$, \textit{i.e.}
the underlying action of $G_{0}$ on $M_{0}$ is transitive and the map
$$
\gg\longrightarrow T_{p}M\qquad,
\phantom{cccccccccccc}A\mapsto A|_{e}\circ \rho_{p}^{*}=\hat{A}|_{p}
$$
is onto. 
The triple $(M_{0},\mathcal{A}_{M},\rho)$ is called a {\it homogeneous $G$-supermanifold}.
}
\end{definition}
\begin{definition}[\cite{Go}]
\rm{The {\it stability subgroup} at $p\in M_{0}$ is the closed Lie subsupergroup $G_{p}$ of $G$ defined by the sHC pair
$
G_{p}:=((G_{0})_{p},\gg_{p})
$
where
$$
(G_{0})_{p}:=\left\{h\in G_{0}\,|\, h(p)=p\right\}
$$
is the stability subgroup of $G_{0}$ and 
$$
\qquad\gg_{p}:=\left\{B\in\gg\,|\, \hat{B}|_{p}=0\right\}
$$
consists of the fundamental vector fields with zero value at the point $p\in M_{0}$.
The action $\Ad:(G_{0})_{p}\longrightarrow\Aut(\gg_{p})$ is the restriction of $\Ad:G_{0}\longrightarrow\Aut(\gg)$.
}\end{definition}
The sheaf theoretic description of the stability subgroup is more complicated and can be found in \cite{BV}. Note that \cite{BV} defines the canonical closed embedding $G_{p}\hookrightarrow G$ in terms of commutativity of the following diagram
\begin{equation*}
\begin{CD}
G_{p} @>>> \bR^{0,0} \\
@VVV               @VV p V    \\
G     @>\rho_{p}>>    M 
\end{CD}
\end{equation*}
\subsection{Homogeneous supermanifold}\hspace{0 cm}\newline\\
Let $(H,\gh)\subseteq (G,\gg)$ be a closed Lie subsupergroup of a Lie supergroup $G=(G_{0},\mathcal{A}_{G})$. Denote by 
$$
\pi_{0}:G_{0}\rightarrow G_{0}/H_{0}\quad,\quad pr_{1}:G\times H\rightarrow G\quad,\quad R_{H}:=m_{|_{G\times H}}:G\times H\rightarrow G
$$
the canonical projections and the right action of $H$ on $G$. For every open set $U\subseteq G_{0}/H_{0}$, define three subalgebras of $\cA_{G}(\pi_{0}^{-1}U)$:
\begin{itemize}
\item[i)] The superalgebra of $R_{H}$-invariant superfunctions
$$
\mathcal{A}_{G/H}(U):=\left\{f\in \mathcal{A}_{G}(\pi_{0}^{-1}U)\,|\,R_{H}^{*}f=pr_{1}^{*}f\right\},
$$
\item[ii)]
The superalgebra of $R_{H_{0}}$-invariant superfunctions 
$$
\mathcal{A}_{H_{0}}(U):=\left\{f\in \mathcal{A}_{G}(\pi_{0}^{-1}U)\,|\,(R_{h})^{*}f=f,\forall h\in H_{0}\right\},
$$
\item[iii)]
The superalgebra of $\gh$-invariant superfunctions
$$
A_{\gh}(U):=\left\{f\in \mathcal{A}_{G}(\pi_{0}^{-1}U)\,|\, Bf=0,\forall B\in\gh\right\}.
$$
\end{itemize}
The associated sheaves over $G_{0}/H_{0}$, denoted by $\mathcal{A}_{G/H}$, $\mathcal{A}_{H_{0}}$ and $\mathcal{A}_{\gh}$ satisfy
\be
\label{homogeneous}
\mathcal{A}_{G/H}=\mathcal{A}_{H_{0}}\cap \mathcal{A}_{\gh}
\ee
\begin{theorem}[\cite{FLV, Kt}]
\label{Fio}
The pair $G/H=(G_{0}/H_{0},\mathcal{A}_{G/H})$
is a supermanifold. There exists a canonical projection
$
\pi:G\longrightarrow G/H
$
satisfying 
\begin{itemize}
\item[1)] $\Ker(\pi_{*})=\mathcal{A}(G)\otimes\gh$
\item[2)] $\pi_{*,e}:T_{e}G\longrightarrow T_{o}(G/H)$ is onto with kernel $\gh\subseteq\gg\cong T_{e}G$ $(o:=eH=\pi_{0}(e))$
\end{itemize}
and a canonical left action $\mu:G\times G/H\rightarrow G/H$ such that the following diagram is commutative
\begin{equation*}
\begin{CD}
G\times G @>m>> G \\
@V id_{G}\times \pi VV               @VV\pi V    \\
G\times G/H     @>\mu>>    G/H 
\end{CD}
\end{equation*}
For every $G$-supermanifold $(M_{0},\mathcal{A}_{M},\rho)$,
the orbit map 
\be
\label{identification}
\rho_{p}=((\rho_{p})_{0},\rho_{p}^{*}):G/G_{p}\rightarrow M
\ee
is well-defined. Whenever $(M_{0}, \mathcal{A}_{M}, \rho)$ is homogeneous, the orbit map (\ref{identification}) is a diffeomorphism with
$(\Id\otimes \rho_{p}^{*})\circ\rho^{*}=\mu^{*}\circ \rho_{p}^{*}$.
\end{theorem}
As in the classical case, there exist adapted coordinates for the canonical projection, \textit{i.e.} $\cA_{G}(\pi_{0}^{-1}U)\cong\cA_{G/H}(U)\otimes\cA(H)$ for a suitable open subset $U\subseteq G_{0}/H_{0}$.
A homogeneous supermanifold $G/H$ is {\it reductive} if the Lie superalgebra $\gg$ of $G$ can be decomposed into a subsupervector space direct sum 
\be
\label{reductived}
\gg=\gh+\gm
\ee
such that $\Ad_{h}\gm\subseteq\gm$ for every $h\in H_{0}$ and $[\gh,\gm]\subseteq\gm$.
If $[\gm,\gm]\subseteq\gh$, then $G/H$ is {\it symmetric}.
\begin{example}
\label{spst}
\rm{
Let $G$ be a Poincare' Lie supergroup and $H_{0}:=\Spin_{r,s}^{0}\subseteq G_{0}$. The {\it Poincare' superspacetime $M=G/H_{0}$} is a reductive homogeneous supermanifold.
}\end{example}
Lemma \ref{trivial} and Theorem \ref{Fio} imply the following useful description of the space of vector fields of a reductive homogeneous supermanifold.
\begin{lemma}
\label{utile}
Let $G/H$ be a homogeneous supermanifold with reductive decomposition (\ref{reductived}).
The map
$$
[\mathcal{A}(G)\otimes\gm]^{R_{H}}\rightarrow \mathcal{T}(G/H)\qquad,\qquad
X\mapsto X\circ\pi^{*}\quad.
$$
is an isomorphism.
\end{lemma}
The induced Lie superalgebra structure on $[\mathcal{A}(G)\otimes\gm]^{R_{H}}$ is given by
$$
[X,Y]=[\sum_{i}f^{i}\otimes A_{i},\sum_{j}g^{j}\otimes B_{j}]:=\sum_{i,j}(-1)^{|A_{i}||g^{j}|}f^{i}g^{j}\otimes[A_{i},B_{j}]_{\gm}\phantom{ccccccccccccccccccccccccc}
$$
$$
\phantom{cccccccccccccccccccccccccccc}+\sum_{j}(Xg^{j})\otimes B_{j}-\sum_{i}(-1)^{|X||Y|}(Yf^{i})\otimes A_{i}
$$
where $[A,B]:=[A,B]_{\gh}+[A,B]_{\gm}\in\gh+\gm=\gg$ for every $A,B\in\gm$.
\subsection{Isotropy representation}
\begin{definition}[\cite{Go}]
\rm{
The {\it linear isotropy representation} $\phi=(\varphi_{\0},\varphi):G_{p}\rightarrow \GL_{\bR}(T_{p}M)$ of the stability subgroup $G_{p}:=((G_{0})_{p},\gg_{p})$ on $T_{p}M$
is defined by
$$
\begin{cases}
\varphi_{\0}:(G_{0})_{p}\rightarrow\GL_{\bR}(T_{p}M) & h\mapsto h_{*,p}\\
\varphi:\gg_{p}\rightarrow\ggl_{\bR}(T_{p}M) &
B\mapsto(v\mapsto-[\hat{B},v]|_{p}=
(-1)^{|B||v|}v\circ \hat{B})
\end{cases}
$$
}\end{definition}
The natural extension of this representation to $T_{p}M^{r}_{s}$ is also denoted by $\phi:G_{p}\rightarrow GL_{\bR}(T_{p}M^{r}_{s})$ and the set of invariant $(r,s)$-tensor $T\in T_{p}M^{r}_{s}$ by $(T_{p}M^{r}_{s})^{\phi(G_{p})}$.
\begin{lemma}[\cite{Go}]
Let $(M_{0}, \mathcal{A}_{M}, \rho)$ be a homogeneous $G$-supermanifold. The linear isotropy representation is equivalent to the adjoint representation $\Ad:G_{p}\longrightarrow GL_{\bR}(\gg/\gg_{p})$ via the natural isomorphism $\gg/\gg_{p}\cong T_{p}M$, .
\end{lemma}
\begin{definition}
\label{first}
\rm{A superfunction $f\in \mathcal{A}(M)$ (resp. a vector field $X\in\mathcal{T}(M)$, resp. a $1$-form $\omega\in\mathcal{T}^{*}(M)$) on a $G$-supermanifold $(M_{0},\mathcal{A}_{M},\rho)$ is {\it$G$-invariant} if
$$
f\cong 1\otimes f=\rho^{*}f\in\mathcal{A}(G\times M)
\phantom{c}\phantom{c}({\rm resp.}\phantom{c}
(\Id\otimes X)\circ\rho^{*}=\rho^{*}\circ X,\phantom{c}{\rm resp.}
$$
$$
\omega(Y)\cong 1\otimes\omega(Y)=(\rho^{*}\omega)(\Id\otimes Y)\in \mathcal{A}(G\times M)\quad{\rm for}\phantom{c}{\rm every}\phantom{c}Y\in\mathcal{T}(M))\quad.
$$
}\end{definition}
The definition of $G$-invariance of a tensor field $T\in\mathcal{T}^{r}_{s}(M)$ follows naturally from Definition \ref{first}. The set of all $G$-invariant tensor fields is denoted by $\oplus_{r,s}\mathcal{T}_{s}^{r}(M)^{G}$.
\begin{definition}
\rm{A tensor field $T\in\mathcal{T}_{s}^{r}(M)$ on a $(G_{0},\gg)$-supermanifold $(\mathcal{A}_{M},\rho_{\ol{0}},\hat{\rho})$ is $(G_{0},\gg)$-{\it invariant} if it is preserved by the global action of $G_{0}$ and annihilated by the fundamental action of $\gg$, \textit{i.e.} if
\beq
\label{invariantfinale}
\begin{cases}
g_{*}(T)=T & \forall g\in G_{0}\\
\mathcal{L}_{\hat{A}}T=0 & \forall A\in\gg
\end{cases}
\eeq
The set of all $(G_{0},\gg)$-invariant tensor fields is denoted by $\oplus_{r,s}\mathcal{T}_{s}^{r}(M)^{(G_{0},\gg)}$.
}\end{definition}
\begin{theorem}
\label{primaparte}
Let $(M_{0},\mathcal{A}_{M}, \rho)$ be a homogeneous $G$-supermanifold. The evaluation at a point $p\in M_{0}$
$$
\ev_{p}:\mathcal{T}_{s}^{r}(M)^{G}\rightarrow (T_{p}M^{r}_{s})^{\phi(G_{p})}
$$
is an isomorphism between the space $\mathcal{T}_{s}^{r}(M)^{G}$ of $G$-invariant tensor fields of type $(r,s)$ and the space $(T_{p}M^{r}_{s})^{\phi(G_{p})}$ of isotropy invariant tensors of type $(r,s)$. In particular, for vector fields, the inverse map is given by
\be
\label{tenosi}
(T_{p}M)^{\phi(G_{p})}\owns v\rightarrow X_{v}:=(\Id\otimes v)\circ\rho^{*}:\mathcal{A}_{M}\longrightarrow \mathcal{A}_{G/G_{p}}\underset{\rho_{p}}{\cong} \mathcal{A}_{M}
\ee
Moreover 
$\oplus_{r,s}\mathcal{T}_{s}^{r}(M)^{G}=\oplus_{r,s}\mathcal{T}_{s}^{r}(M)^{(G_{0},\gg)}$.
\end{theorem}
\begin{pf}
See the Appendix.
\end{pf}
\section{Invariant superconnections on a homogeneous supermanifold
}
\setcounter{equation}{0}
\label{section4}
This section generalizes a theorem of Wang (\cite{KN1, KN2}), concerning invariant connections on homogeneous manifolds, to the category of supermanifolds. The main Theorem \ref{con1con2} is stated in \cite{Co} in the case of even stability subgroup.
\subsection{Superconnection}\hspace{0 cm}\newline\\
Let $\mathcal{E}$ be a locally free sheaf of $\mathcal{A}_{M}$-supermodules on $M_{0}$.
\begin{definition}[\cite{DM}]
\rm{A {\it (super)connection} on $\mathcal{E}$ is a morphism $\nabla:\mathcal{T}_{M}\otimes_{\bR}\mathcal{E}\rightarrow \mathcal{E}$ of sheaves of supervector spaces such that, for every open subset $U\subseteq M_{0}$
\begin{itemize}
\item[1)]
$
\nabla_{fX}Z=f\cdot\nabla_{X}Z\qquad
$
\item[2)]$\nabla_{X}fZ=(Xf)Z+(-1)^{|X||f|}f\nabla_{X}Z$
\end{itemize}
where $f\in \mathcal{A}_{M}(U)$, $X\in\mathcal{T}_{M}(U)$, $Z\in\mathcal{E}(U)$. 
The {\it curvature} of $\nabla$ is defined by
$$
\phantom{cccccccccc}R(X,Y)Z:=\nabla_{X}\nabla_{Y}Z-(-1)^{|X||Y|}\nabla_{Y}\nabla_{X}Z-\nabla_{[X,Y]}Z
$$
where $X,Y\in\mathcal{T}_{M}(U)$ and $Z\in\mathcal{E}(U)$.}
\end{definition}
A connection on $\mathcal{T}_{M}$ is called a {\it linear connection} on $M$. There is the usual notion of {\it Christoffel symbols}, in  
local coordinates (\ref{COORD}) on $M$,
$$
\nabla_{\frac{\partial}{\partial\eta^{i}}}\frac{\partial}{\partial\eta^{j}}=\sum_{k}\Gamma_{ij}^{k}\frac{\partial}{\partial\eta^{k}}
$$
gives elements $\Gamma_{ij}^{k}\in \mathcal{A}_{M}(U)$ of parity $|\Gamma_{ij}^{k}|=|\eta^{i}|+|\eta^{j}|+|\eta^{k}|$. 
The {\it torsion} is defined by
$$
T(X,Y):=\nabla_{X}Y-(-1)^{|X||Y|}\nabla_{Y}X-[X,Y]
$$
where $X,Y\in\mathcal{T}_{M}(U)$.
The {\it covariant derivatives} of $R$ are denoted by
$$
\nabla^{r}_{X_{r},...,X_{1}}R:=\nabla_{X_{r}}\circ\cdot\cdot\cdot\circ\nabla_{X_{1}}R
$$
where $X_{r},..,X_{1}\in\mathcal{T}_{M}(U)$ and $r\in\bN$.
\begin{definition}[\cite{Ga}]
\rm{
The {\it infinitesimal holonomy superalgebra} $\mathfrak{hol(\nabla)_{p}^{inf}}$ at a point $p\in M_{0}$ of a linear connection $\nabla:\mathcal{T}_{M}\otimes_{\bR}\mathcal{T}_{M}\rightarrow \mathcal{T}_{M}$ on M is the linear Lie superalgebra $\mathfrak{hol(\nabla)_{p}^{inf}}\subseteq\ggl_{\bR}(T_{p}M)=\ggl_{\bR}(T_{p}M)_{\ol{0}}+\ggl_{\bR}(T_{p}M)_{\ol{1}}$ spanned by 
$$
\nabla^{r}_{X_{r},..,X_{1}}R_{p}(X,Y):T_{p}M\rightarrow T_{p}M
$$
where $X_{r},..,X_{1},X,Y\in T_{p}M$, $r\in\bN$.
}\end{definition}
For the definition of the {\it holonomy supergroup} as a sHC pair see \cite{Ga}. Note that the restriction 
$$
\widetilde{\nabla}=(\nabla|_{\Gamma(TM_{0})\otimes\Gamma(TM)}):\Gamma(TM_{0})\otimes\Gamma(TM)\rightarrow\Gamma(TM)
$$
is a connection on the tangent bundle $TM$, such that $TM_{0}$ and $(TM)_{\ol{1}}$ are $\nabla$-stable. The holonomy of this connection is contained in the body of the holonomy supergroup (\cite{Ga}). This is, in general, a proper inclusion (\cite{Ga}).\\
Let $\phi=(\phi_{0},\phi^{*})\in\Mor(N,M)$ be a morphism of supermanifolds and recall that the sheaf of $\phi$-vector fields (\ref{alongv}) is denoted by $\mathcal{T}_{\phi}$. 
\begin{definition}
\rm{
A morphism $\nabla^{\phi}:(\phi_{0})_{*}\mathcal{T}_{N}\otimes_{\bR}\mathcal{T}_{\phi}\rightarrow \mathcal{T}_{\phi}$ of sheaves of supervector spaces is a {\it $\phi$-connection} if, for every open subset $U\subseteq M_{0}$
\begin{itemize}
\item[1)]
$
\nabla^{\phi}_{fY}X=f\cdot\nabla_{Y}^{\phi}X\qquad
$
\item[2)]$\nabla^{\phi}_{Y}fX=(Yf)X+(-1)^{|Y||f|}f\nabla^{\phi}_{Y}X$
\end{itemize}
where $f\in \mathcal{A}_{N}(\phi_{0}^{-1}U)$, $Y\in\mathcal{T}_{N}(\phi_{0}^{-1}U)$, $X\in\mathcal{T}_{\phi}(U)$.
}\end{definition}
\begin{example}
\label{pb}
\rm{
The {\it pull-back connection} $(\phi^{*}\nabla):(\phi_{0})_{*}\mathcal{T}_{N}\otimes_{\bR}\mathcal{T}_{\phi}\rightarrow\mathcal{T}_{\phi}$ of a linear connection $\nabla:\mathcal{T}_{M}\otimes_{\bR}\mathcal{T}_{M}\rightarrow \mathcal{T}_{M}$ is a $\phi$-connection
defined by the local expression
$$
(\phi^{*}\nabla)_{Y}(X):=\sum_{k}(Yf^{k})\cdot(\phi^{*}\circ \frac{\partial}{\partial \eta^{k}})+\sum_{k}(-1)^{|Y||f^{k}|}f^{k}\cdot (\phi^{*}\nabla)_{Y}(\phi^{*}\circ\frac{\partial}{\partial \eta^{k}})
$$
$$
\phantom{cccccccccc}:=\sum_{k}(Yf^{k})\cdot(\phi^{*}\circ \frac{\partial}{\partial \eta^{k}})+\sum_{i,j,k}(-1)^{|Y||f^{k}|}f^{k}\cdot(Y\circ\phi^{*}\eta^{j})\cdot(\phi^{*}\Gamma_{jk}^{i}) (\phi^{*}\circ\frac{\partial}{\partial\eta^{i}})\qquad 
$$
where $Y\in\mathcal{T}_{N}(\phi_{0}^{-1}U)$ and $X\in\mathcal{T}_{\phi}(U)$ is given by (\ref{along}). The usual functorial property holds, \textit{i.e.} 
for every morphism $\psi\in\Mor(L,N)$
\be
\label{antani}
\psi^{*}\{(\phi^{*}\nabla)_{Y}(X)\}=((\phi\circ\psi)^{*}\nabla)_{Z}(\psi^{*}\circ X)
\ee
where $Z\in\mathcal{T}_{L}(\psi_{0}^{-1}\phi_{0}^{-1}U)$ is a vector field $\psi$-related to $Y\in\mathcal{T}_{N}(\phi_{0}^{-1}U)$.
}\end{example}
\subsection{Nomizu map}
\begin{definition}
\label{intensificare}
\rm{A linear connection $\nabla:\mathcal{T}_{M}\otimes_{\bR}\mathcal{T}_{M}\rightarrow \mathcal{T}_{M}$ on a $G$-supermanifold $(M_{0},\mathcal{A}_{M},\rho)$ is {\it $G$-invariant} if
\be
\label{nonso}
(\Id\otimes \nabla_{X}Y)\circ\rho^{*}=(\rho^{*}\nabla)_{\Id\otimes X}((\Id\otimes Y)\circ\rho^{*})
\ee
for every $X,Y\in\mathcal{T}(M)$. 
A linear connection $\nabla:\mathcal{T}_{M}\otimes_{\bR}\mathcal{T}_{M}\rightarrow \mathcal{T}_{M}$ on a $(G_{0},\gg)$-supermanifold $(\mathcal{A}_{M},\rho_{\ol{0}},\hat{\rho})$ is {\it $(G_{0},\gg)$-invariant} if
\be
\label{con1}
g_{*}(\nabla_{X}Y)=\nabla_{g_{*}X}g_{*}Y 
\ee
and
\be
\label{con2}
\mathcal{L}_{\hat{A}}(\nabla_{X}Y)=\nabla_{\cL_{\hat{A}}X}Y+(-1)^{|X||A|}\nabla_{X}\cL_{\hat{A}}Y
\ee
for every $X,Y\in\mathcal{T}(M)$, $g\in G_{0}$ and $A\in\gg$.}
The space of all $G$ (resp. $(G_{0},\gg)$)-invariant linear connections is denoted by $\Conn(M)^{G}$ (resp. by $\Conn(M)^{(G_{0},\gg)}$).
\end{definition}
\begin{definition}
\label{Nomizudefinition}
\rm{An even linear map $L:\gg\rightarrow \ggl_{\bR}(T_{p}M)$ is a {\bf Nomizu map} in $p\in M_{0}$ if
\begin{itemize}
\item[1)] $L(B)=-\mathcal{L}_{\hat{B}}(\cdot)|_{p}\phantom{ccccccccccccccccccc}\;$ for all $B\in\gg_{p}$,
\item[2)] $L(\Ad_{h}A)=h_{*,p}\circ L(A)\circ h_{*,p}^{-1}\phantom{ccccccc}$ for all $h\in (G_{0})_{p}$, $A\in\gg$,
\item[3)] $L([B,A])=[L(B),L(A)]\phantom{cccccccccccc}$ for all $A\in\gg$, $B\in\gg_{p}$.
\end{itemize}
The space of all Nomizu maps in $p\in M_{0}$ is denoted by $\Nom_{p}(M)$.
}\end{definition}
The linear operator
$$
L_{X}:=\nabla_{X}-\mathcal{L}_{X}:\mathcal{T}(M)\longrightarrow \mathcal{T}(M)\qquad\quad X\in\mathcal{T}(M)
$$
is $\mathcal{A}(M)$-linear. Its value
$
L_{X}|_{p}:T_{p}M\rightarrow T_{p}M
$ at a point $p\in M_{0}$ is well-defined. Let $\nabla\in\Conn(M)^{(G_{0},\gg)}$ be a $(G_{0},\gg)$-invariant linear connection on $(\cA_{M},\rho_{\0},\hat{\rho})$. The even linear map
\be
\label{Nomizone}
L^{\nabla}_{p}:\gg\longrightarrow \ggl_{\bR}(T_{p}M)\phantom{cc},
\phantom{cc}A\mapsto L_{\hat{A}}|_{p}
\ee
is called the {\bf Nomizu map} in $p$ associated with $\nabla$.
\subsection{Wang's theorem for supermanifolds}
\subsubsection{Main theorem}
\begin{theorem}
\label{teoremafinale}
\label{con1con2}
Let $(M_{0},\mathcal{A}_{M}, \rho)$ be a homogeneous $G$-supermanifold and let $p\in M_{0}$ be a fixed point. The correspondence 
\be
\label{correspdan}
\Conn(M)^{G}\ni\nabla\mapsto L_{p}^{\nabla}\in\Nom_{p}(M)
\ee
is a bijection. Moreover $\Conn(M)^{G}=\Conn(M)^{(G_{0},\gg)}$.
\end{theorem}
\begin{pf}
The proof is split in five parts.
\\
{\bf i)} Every $G$-invariant linear connection on a $G$-supermanifold $(M_{0},\mathcal{A}_{M}, \rho)$ is $(G_{0},\gg)$-invariant, \textit{i.e.} $\Conn(M)^{G}\subseteq\Conn(M)^{(G_{0},\gg)}$.\\
Let $X,Y\in\mathcal{T}(M)$ be vector fields on $M$. Equation (\ref{con1}) follows from applying $(\ev_{g}\otimes \Id)$ to both sides of (\ref{nonso}) and using equation (\ref{antani}).
Equation (\ref{con2}) is a (long) calculation in local coordinates (\ref{COORD}) which uses the following equality
$$
-(-1)^{|A|(|j|+|i|)}(\nabla_{\frac{\partial}{\partial\eta^{j}}}\frac{\partial}{\partial\eta^{i}})\circ \hat{A}=-(-1)^{|A|(|j|+|i|)}(\nabla_{\frac{\partial}{\partial\eta^{j}}}\frac{\partial}{\partial\eta^{i}})\circ (A|_{e}\otimes \Id)\circ\rho^{*}
$$
$
=-(A|_{e}\otimes \Id)\circ (\Id\otimes\nabla_{\frac{\partial}{\partial\eta^{j}}}\frac{\partial}{\partial\eta^{i}})\circ\rho^{*}=-(A|_{e}\otimes \Id)(\rho^{*}\nabla)_{\Id\otimes\frac{\partial}{\partial\eta^{j}}}((\Id\otimes\frac{\partial}{\partial\eta^{i}})\circ \rho^{*})\quad.
$
{\bf ii)} The Nomizu map (\ref{Nomizone}) associated to a $(G_{0},\gg)$-invariant linear connection on a $(G_{0},\gg)$-supermanifold $(\mathcal{A}_{M},\rho_{\ol{0}},\hat{\rho})$ is an element of $\Nom_{p}(M)$.\\
Let $v\in T_{p}M$ be a vector at $p\in M_{0}$ and $X\in\mathcal{T}(M)$ be a vector field on $M$ such that $X|_{p}=v$. For every $B\in\gg_{p}$
$$
\phantom{cccc}L_{p}^{\nabla}(B)(v)=(\nabla_{\hat{B}}X-\mathcal{L}_{\hat{B}}X)|_{p}=(\nabla_{\hat{B}}X)|_{p}-\mathcal{L}_{\hat{B}}(X)|_{p}=-\mathcal{L}_{\hat{B}}(X)|_{p}\quad.
$$
Lemma \ref{obvious} together with equation (\ref{con1}) imply that
$$
L_{p}^{\nabla}(\Ad_{h}A)=(\nabla_{\hat{\Ad_{h}A}}-\mathcal{L}_{\hat{\Ad_{h}A}})|_{p}=(\nabla_{h_{*}\hat{A}}-\mathcal{L}_{h_{*}\hat{A}})|_{p}=
h_{*,p}\circ L_{p}^{\nabla}(A)\circ h^{-1}_{*,p}
$$
for every $h\in (G_{0})_{p}$ and $A\in\gg$. For every $A\in\gg$ and $B\in\gg_{p}$
\begin{align*}
&[-\mathcal{L}_{\hat{B}},L_{\hat{A}}](X)\\\
&=(\mathcal{L}_{\hat{B}}\mathcal{L}_{\hat{A}}X-(-1)^{|B||A|}\mathcal{L}_{\hat{A}}\mathcal{L}_{\hat{B}}X)+
(-\mathcal{L}_{\hat{B}}\nabla_{\hat{A}}X+(-1)^{|B||A|}\nabla_{\hat{A}}\mathcal{L}_{\hat{B}}X)\\\
&=\mathcal{L}_{[\hat{B},\hat{A}]}X-\nabla_{[\hat{B},\hat{A}]}X=(\nabla_{\hat{[B,A]}}-\mathcal{L}_{\hat{[B,A]}})X
\end{align*}
where the second to last equality follows from (\ref{Jacobi}) and (\ref{con2}). Evaluating in $p\in M_{0}$, the assertion follows.\\
{\bf iii)} Injectivity of the correspondence (\ref{correspdan}). \\
Let $\nabla^{i}:\mathcal{T}_{M}\otimes_{\bR}\mathcal{T}_{M}\rightarrow \mathcal{T}_{M}$, $i=1,2$, be two $(G_{0},\gg)$-invariant linear connections on $M$  such that the associated Nomizu maps $L^{1}_{p}, L^{2}_{p}:\gg\rightarrow \ggl_{\bR}(T_{p}M)$ coincide, \textit{i.e.}
$$
(\nabla^{1}_{\hat{A}}-\mathcal{L}_{\hat{A}})|_{p}=(\nabla^{2}_{\hat{A}}-\mathcal{L}_{\hat{A}})|_{p}
$$
for every $A\in\gg$. For every $X\in\mathcal{T}(M)$, equation (\ref{con1}) and Lemma \ref{obvious} imply that
$$
g_{*}(\nabla^{i}_{\hat{A}}X-\mathcal{L}_{\hat{A}}X)=\nabla^{i}_{g_{*}\hat{A}}g_{*}X-\mathcal{L}_{g_{*}\hat{A}}g_{*}X=(\nabla^{i}_{\hat{\Ad_{g}A}}g_{*}X-\mathcal{L}_{\hat{\Ad_{g}A}}g_{*}X)\quad.
$$
Evaluating in $g^{-1}p\in M_{0}$
$$
L^{i}_{g^{-1}p}(A)=g^{-1}_{*}|_{p}\circ L^{i}_{p}(\Ad_{g}A)\circ g_{*}|_{g^{-1}p}:T_{g^{-1}p}M\rightarrow T_{g^{-1}p}M
$$
and then
$
L^{1}_{q}=L^{2}_{q}
$
for every $q\in M_{0}$. Equations (\ref{con2}) and (\ref{Jacobi}) imply that
$$
\mathcal{L}_{\hat{B}}(L^{i}_{\hat{A}}X)=-L^{i}_{\hat{\cL_{B}A}}X+(-1)^{|A||B|}L^{i}_{\hat{A}}\cL_{\hat{B}}X
$$
for every $A,B\in\gg$ and $X\in\mathcal{T}(M)$. In particular
$$
\mathcal{L}_{\gg}(L^{i}_{\gg}\mathcal{T}(M))\subseteq L^{i}_{\gg}\mathcal{T}(M)
$$
and iterating (for every $0\leq k\lneq \infty $)
$$
\mathcal{L}_{\hat{A}_{k}}\cdot\cdot\cdot\mathcal{L}_{\hat{A}_{1}}(L^{i}_{\hat{A}}X)\subseteq L^{i}_{\gg}\mathcal{T}(M)
$$
where $A_{i}\in\gg$ for $1\leq i\leq k$. For every $0\leq k\lneq \infty $
$$
(\mathcal{L}_{\hat{A}_{k}}\cdot\cdot\cdot\mathcal{L}_{\hat{A}_{1}}(L^{1}_{\hat{A}}X))|_{q}=
(\mathcal{L}_{\hat{A}_{k}}\cdot\cdot\cdot\mathcal{L}_{\hat{A}_{1}}(L^{2}_{\hat{A}}X))|_{q}\qquad.
$$
Lemma \ref{appendixA2} applied to the vector field 
$
L^{1}_{\hat{A}}X-L^{2}_{\hat{A}}X=\nabla^{1}_{\hat{A}}X-\nabla^{2}_{\hat{A}}X
$
implies that $\nabla^{1}_{\hat{A}}X=\nabla^{2}_{\hat{A}}X$
and, from Nakayama's Lemma, $\nabla^{1}=\nabla^{2}$. \\
{\bf iv)} In order to prove surjectivity of the correspondence (\ref{correspdan}), some preliminaries are needed. \\
Denote $G_{p}$ by $H$ and let $G/H$ be the associated homogeneous supermanifold (recall Theorem \ref{Fio}). Say that 
a $\pi$-connection $\nabla^{\pi}:(\pi_{0})_{*}\mathcal{T}_{G}\otimes_{\bR}\mathcal{T}_{\pi}\rightarrow \mathcal{T}_{\pi}$ is {\it projectable} if
\begin{eqnarray}
\label{firstproperty}
R_{H_{0}}-\textbf{Invariance}\qquad 
R_{h}^{*}\circ(\nabla^{\pi}_{Y}X) & = & \nabla^{\pi}_{(R^{-1}_{h})_{*}Y}(R_{h}^{*}\circ X) \\
\label{secondproperty}
\gh-\textbf{Invariance}\,\,\,\,\,\qquad 
\mathcal{L}_{B}(\nabla^{\pi}_{Y}X) & = & \nabla^{\pi}_{\mathcal{L}_{B}Y}X+(-1)^{|Y||B|}\nabla^{\pi}_{Y}\mathcal{L}_{B}X\,\,\,\,\,\,\,\,\,\
\end{eqnarray}
for every $h\in H_{0}$, $Y\in\mathcal{T}_{G}(\pi_{0}^{-1}U)$, $X\in\mathcal{T}_{\pi}(U)$, $B\in\gh$ and
\begin{eqnarray}
\label{thirdproperty}
\textbf{Horizontality}\qquad
\nabla^{\pi}_{Y}(f\cdot \pi^{*}\circ X) & = & (Yf)\cdot(\pi^{*}\circ X)\phantom{ccccccccccccc}
\end{eqnarray}
for every $Y\in \mathcal{A}_{G}(\pi_{0}^{-1}U)\otimes\gh$, $f\in \mathcal{A}_{G}(\pi_{0}^{-1}U)$ and $X\in\mathcal{T}_{G/H}(U)$.
The pull-back
\be
\label{tensorialcon}
\nabla\overset{\pi^{*}}{\rightarrow} (\pi^{*}\nabla)
\ee
is a bijection between the set of linear connections on $G/H$ and the set of projectable $\pi$-connections. We sketch a proof of this fact.
Let $\nabla:\mathcal{T}_{G/H}\otimes_{\bR}\mathcal{T}_{G/H}\rightarrow\mathcal{T}_{G/H}$ be a linear connection on $G/H$. One can check that the pull-back connection $\nabla^{\pi}:=(\pi^{*}\nabla)$ is projectable. Here we only remark that (\ref{firstproperty})
is a consequence of (\ref{antani}).
Consider a projectable $\pi$-connection $\nabla^{\pi}:(\pi_{0})_{*}\mathcal{T}_{G}\otimes_{\bR}\mathcal{T}_{\pi}\rightarrow \mathcal{T}_{\pi}$.
Using equation (\ref{thirdproperty}) and adapted coordinates near $h\in H_{0}\subseteq G_{0}$ and $o\in G_{0}/H_{0}$, it is possible to construct (locally near $o\in U\subseteq G_{0}/H_{0}$) a map
$
\nabla:\mathcal{T}_{G/H}\otimes_{\bR}\mathcal{T}_{G/H}\rightarrow\mathcal{T}_{\pi}
$.
By definition 
$$
\nabla_{Z}X:=\nabla^{\pi}_{Y}(\pi^{*}\circ X)
$$
where $X,Z\in\mathcal{T}_{G/H}(U)$ and $Y\in\mathcal{T}_{G}(\pi^{-1}_{0}U)$ is $\pi$-related to $Z$.
From (\ref{firstproperty}) and (\ref{secondproperty}), check that locally
$
\nabla:\mathcal{T}_{G/H}\otimes_{\bR}\mathcal{T}_{G/H}\rightarrow\mathcal{T}_{G/H}
$.
Repeating this construction for every $gH\in G_{0}/H_{0}$, we get a globally well-defined linear connection on $G/H$. \\
Due to Theorem \ref{Fio} and isomorphism (\ref{tensorialcon}), we need to prove that to any even map $L:\gg\rightarrow \ggl_{\bR}(\gg/\gh)$
satisfying
\begin{itemize}
\item[1)] $L(B)=\ad_{B}$ for every $B\in\gh$,
\item[2)] $L(\Ad_{h}A)=\Ad_{h}\circ L(A)\circ \Ad_{h^{-1}}$ for every $h\in H_{0}$, $A\in\gg$,
\item[3)] $L([B,A])=[L(B),L(A)]$ for every $A\in\gg$, $B\in\gh$,
\end{itemize}
there exists a corresponding projectable $\pi$-connection.
Every $X\in\mathcal{T}_{\pi}$ can be locally written as
$$
X=\sum_{j}(g^{j}\otimes B_{j})\circ\pi^{*}
$$
where $g^{j}\otimes B_{j}\in \mathcal{A}_{G}\otimes\gg$ and two such expressions differ by an element of $\mathcal{A}_{G}\otimes\gh$. The following $\pi$-connection is well-defined
$$
\nabla^{\pi}:\mathcal{T}_{G}(\pi_{0}^{-1}U)\otimes_{\bR}\mathcal{T}_{\pi}(U)\rightarrow\mathcal{T}_{\pi}(U)
$$ 
\beq
\label{formula}
(\sum_{i} (f^{i}\otimes A_{i})\otimes_{\bR} X)\mapsto\sum_{i,j} f^{i}\{(A_{i}g^{j})\cdot(B_{j}\circ\pi^{*})+(-1)^{|A_{i}||g^{j}|}g^{j}(L(A_{i})B_{j})\circ\pi^{*}\}
\eeq
where $\sum_{i}^{i}f^{i}\otimes A_{i}\in\mathcal{A}_{G}(\pi_{0}^{-1}U)\otimes\gg$. It is projectable, for equation (\ref{firstproperty})
is a consequence of 2), equation (\ref{secondproperty}) of 1) and 3) while equation (\ref{thirdproperty}) of 1). The condition of $G$-invariance follows directly from (\ref{formula}).\\
{\bf v)} The inclusion $\Conn(M)^{(G_{0},\gg)}\subseteq\Conn(M)^{G}$ is now easy to prove.
\qedhere
\end{pf}
In the reductive case, the set of $G$-invariant linear connections is in bijective correspondence with the set of even linear mappings
$
L_{\gm}:\gm\rightarrow\ggl_{\bR}(\gm)
$
satisfying
\begin{itemize}
\item[$\cdot$] $L_{\gm}(\Ad_{h}A)=\Ad_{h}\circ L_{\gm}(A)\circ \Ad_{h^{-1}}$ for every $h\in H_{0}$, $A\in\gm$,
\item[$\cdot$] $L_{\gm}([B,A])=[\ad_{B},L_{\gm}(A)]$ for every $A\in\gm$, $B\in\gh$,
\end{itemize}
where $L|_{\gm}=L_{\gm}$, $L|_{\gh}=\ad_{\gh}$ is the Nomizu map. By abuse of language, we refer to the map $L_{\gm}:\gm\rightarrow\ggl_{\bR}(\gm)$ as the {\it Nomizu map}. The following corollaries are direct consequences of the main Theorem \ref{teoremafinale}.
\begin{corollary}
Let $M=G/H$, be a homogeneous supermanifold with reductive decomposition (\ref{reductived}). The $G$-invariant linear connection associated with the Nomizu map $L_{\gm}:\gm\rightarrow\ggl_{\bR}(\gm)$ is given by
$$
\nabla:\mathcal{T}(M)\otimes_{\bR}\mathcal{T}(M)\longrightarrow\mathcal{T}(M)\phantom{cccccccccccccccccc}
$$ 
$$
\phantom{ccccccccccccc}(f^{i}\otimes A_{i})\otimes_{\bR} (g^{j}\otimes B_{j})\mapsto f^{i}\{(A_{i}g^{j})B^{j}+(-1)^{|A_{i}||g^{j}|}g^{j}(L_{\gm}(A_{i})B_{j})\}
$$
where $f^{i}\otimes A_{i}$, $g^{j}\otimes B_{j}\in[\cA(G)\otimes\gm]^{R_{H}}\cong\mathcal{T}(M)$ and we have used the Einstein summation convention.
\end{corollary}
\begin{pf}
Lemma \ref{utile} gives an equivalent description of the space of vector fields
of a homogeneous reductive supermanifold $M=G/H$. 
\end{pf}
\begin{definition}
\rm{Let $G/H$ be a homogeneous supermanifold with reductive decomposition (\ref{reductived}). The {\it canonical} connection $\nabla^{can}$ (resp. the {\it natural torsion-free} connection $\nabla^{free}$) is the $G$-invariant linear connection associated with $$L_{\gm}=0\qquad\qquad ({\rm resp.}\,\, L_{\gm}(A)=\frac{1}{2}[A,\cdot]_{m}\,)\,\qquad\qquad.$$
}\end{definition}
If the stability subgroup is even, \textit{i.e.} $H=H_{0}$, the following corollary holds
\begin{corollary}
\label{gravity}
Let $H=H_{0}\subseteq G$ be a connected closed Lie subgroup of a Lie supergroup $G$ such that $M_{0}=G_{0}/H_{0}$ is a homogeneous manifold with reductive decomposition $\gg_{0}=\gh+\gm_{\0}$. Every pair of even $\gh$-invariant linear map
$$
L_{\gm_{\0}}|_{\gg_{\ou}}\in\Hom_{\bR}(\gm_{\0},\ggl_{\bR}(\gg_{\ou}))^{\gh}_{\0}\quad,\quad L_{\gg_{\ou}}\in\Hom_{\bR}(\gg_{\ou},\ggl_{\bR}(\gm_{\0}+\gg_{\ou}))^{\gh}_{\0}
$$
defines an extension of a fixed Nomizu map $L_{0}\in\Hom_{\bR}(\gm_{\0},\ggl_{\bR}(\gm_{\0}))^{\gh}$ on $M_{0}$
to a Nomizu map $L_{\gm}\in\Hom_{\bR}(\gm,\ggl_{\bR}(\gm))^{\gh}$ on the homogeneous supermanifold $M=G/H_{0}$ with reductive decomposition $\gg=(\gh+\gm_{\0})+\gg_{\ou}=(\gh+\gm_{\0})+\gm_{\ou}=\gh+\gm$.
\end{corollary}
Via Corollary \ref{gravity}, the Nomizu map $L_{\gm}\in\Hom_{\bR}(\gm,\ggl_{\bR}(\gm))^{\gh}$ associated with
\be
\label{supersymmetryconnection}
L_{0}=0\qquad,\qquad L_{\gm_{\0}}|_{\gg_{\ou}}=\ad_{\gm_{\0}}|_{\gg_{\ou}}\quad,\quad L_{\gg_{\ou}}=0
\ee
defines a $G$-invariant linear connection $\nabla^{\cS}$ on $G/H_{0}$ called the {\bf supersymmetry} connection.
\subsubsection{Invariant superconnections on a Poincare' superspacetime}\hfill\newline\\
Invariant linear connections on a Poincare' superspacetime of signature $(r,s)$ are in bijective correspondence with $\Spin_{r,s}^{0}$-equivariant even maps
$$
L:\bR^{r,s}+S\rightarrow\ggl_{\bR}(\bR^{r,s}+S)\qquad,
$$
\textit{i.e} with $\spin_{r,s}$-equivariant maps 
$$
\bR^{r,s}\rightarrow (({\bR^{r,s}})^{*}\otimes \bR^{r,s})\oplus(S^{*}\otimes S)\quad\rm{and}\quad
S\rightarrow ((\bR^{r,s})^{*}\otimes S)\oplus (S^{*}\otimes \bR^{r,s})\quad.
$$
The existence of non-degenerate invariant bilinear forms on $\bR^{r,s}$ and $S$ (see \cite{AC}) implies that $\bR^{r,s}$ and $S$ are self-dual $\spin_{r,s}$-modules, \textit{i.e.} $(\bR^{r,s})^{*}\cong\bR^{r,s}$ and $S^{*}\cong S$.
\begin{theorem}
The dimension $D$ of the vector space of invariant connections on a Poincare' superspacetime of signature $(r,s)$ depends on $r-s \mod 8$ as follows
$$\left[
\begin{array}{c|c|c|c|c|c|c|c|c}
\ r-s \mod 8 & 1 & 2 & 3 & 4 & 5 & 6 & 7 & 8 \\
\hline
\ D & 12 & 24 & 12 & 24 & 12 & 6 & 3 & 6 \\
\end{array}
\right]$$
\end{theorem}
\begin{pf}
We want to determine the dimension of the following three spaces
$$
\Hom(\bR^{r,s},\bR^{r,s}\otimes \bR^{r,s})^{\spin_{r,s}}\phantom{c},\phantom{c}\Hom(\bR^{r,s},S\otimes S)^{\spin_{r,s}}\phantom{c},\phantom{c} \Hom(S,\bR^{r,s}\otimes S)^{\spin_{r,s}}\quad.
$$
$\dim_{\bR}\Hom(\bR^{r,s},\bR^{r,s}\otimes \bR^{r,s})^{\spin_{r,s}}=0$ because $\bR^{r,s}$ is not a submodule of $\bR^{r,s}\otimes \bR^{r,s}$; Theorem 1.1 of \cite{AC} implies that
$$
\dim_{\bR}\Hom(\bR^{r,s},S\otimes S)^{\spin_{r,s}}=\dim_{\bR}\Hom(S,\bR^{r,s}\otimes S)^{\spin_{r,s}}=\dim_{\bR}\mathcal{C}_{r,s}
$$ 
where $\mathcal{C}_{r,s}:=\Hom(S,S)^{\spin_{r,s}}$ is the Schur algebra of $S$. The assertion of the theorem is then a direct consequence of Corollary 1.3 of \cite{AC}.
\end{pf}
\subsubsection{Reconstruction of an invariant connection in local coordinates of $M$}\hfill\newline\\
Let $(M_{0},\mathcal{A}_{M}, \rho)$ be a homogenous $G$-supermanifold of dimension $\dim M=m|n$. When $n$ is small the $G$-invariant linear connection $\nabla:\mathcal{T}_{M}\otimes_{\bR}\mathcal{T}_{M}\rightarrow \mathcal{T}_{M}$ associated with a Nomizu map $L_{p}:\gg\rightarrow\ggl_{\bR}(T_{p}M)$ can be recostructed in local coordinates $\left\{x^{i},\xi_{j}\right\}$ of $M$. For every $q\in M_{0}$, consider
$$
L_{q}(\cdot):=g^{-1}_{*}|_{p}\circ L_{p}(\Ad_{g}\cdot)\circ g_{*}|_{g^{-1}p}:\gg\rightarrow\ggl_{\bR}(T_{q}M)
$$
where $q=g^{-1}p\in M_{0}$ for some $g\in G_{0}$. By definition of transitive action, every $Y\in\mathcal{T}(M)$ can be locally written (non-uniquely) on an open set $U\subseteq M_{0}$ as
$$
Y=\sum_{i}f^{i}\hat{A_{i}}
$$
where $\{A_{i}\}$ is a basis of $\gg$ and $f^{i}\in \mathcal{A}_{M}(U)$.
Write
$$
\nabla_{Y}X|_{q}=f^{i}(q)(L_{q}(A_{i})(X|_{q})+(\mathcal{L}_{\hat{A_{i}}}X)|_{q})
$$
for every $X\in\mathcal{T}_{M}(U)$ and $q\in U$. This is not sufficient to reconstruct the connection unless $n=0$. However, a Taylor expansion with respect to the odd coordinates
$$
\nabla_{Y}X=(\nabla_{Y}X)|_{q}+\sum_{k=1}^{n}\xi_{k}(\mathcal{L}_{\frac{\partial}{\partial\xi_{k}}}\nabla_{Y}X)|_{q}+\sum_{k\lneq j}\xi_{k}\xi_{j}(\mathcal{L}_{\frac{\partial}{\partial\xi_{j}}}\mathcal{L}_{\frac{\partial}{\partial\xi_{k}}}\nabla_{Y}X)|_{q}+...
$$
allows to construct $\nabla$. Indeed the condition of $(G_{0},\gg)$-invariance implies that
\begin{align*}
&(\mathcal{L}_{\frac{\partial}{\partial\xi_{k}}}\nabla_{Y}X)|_{q}\\\
&=\sum_{j}\{g_{k}^{j}(\nabla_{[\hat{A_{j}},Y]}X+(-1)^{|A_{i}||Y|}\nabla_{Y}[\hat{A_{j}},X])
-(-1)^{(|Y|+|X|)}((\nabla_{Y}X)g_{k}^{j})\hat{A_{j}}\}|_{q}
\end{align*}
where 
$
\frac{\partial}{\partial\xi_{k}}=\sum_{j}g_{k}^{j}\hat{A_{j}}
$
is a local (non-unique) description of the odd derivation $\frac{\partial}{\partial\xi_{k}}$.
\subsubsection{Curvature, Torsion and Holonomy}\hfill\newline\\
Let $(M_{0},\mathcal{A}_{M},\rho)$ be a homogenous $G$-supermanifold, $p\in M_{0}$. The isotropy subgroup $G_{p}$ is denoted by $H$. The identifications 
$$(M_{0},\mathcal{A}_{M},\rho)\cong (G_{0}/H_{0}, G/H, \mu)\qquad,\qquad\pi_{*,e}:\gg/\gh\rightarrow T_{o}G/H$$
are tacitly used. The curvature and the torsion of a $G$-invariant linear connection $\nabla:\mathcal{T}_{M}\otimes_{\bR}\mathcal{T}_{M}\rightarrow \mathcal{T}_{M}$ are invariant tensors satisfying for every $A,B\in\gg$
$$
R_{o}(A,B)=[L(A),L(B)]-L([A,B])
$$
$$
\phantom{ccccccccccc}T_{o}(A,B)=L(A)B-(-1)^{|A||B|}L(B)A-[A,B]
$$
where $L:\gg\rightarrow\ggl_{\bR}(\gg/\gh)$ is the Nomizu map associated with $\nabla$. The following lemma then holds. 
\begin{lemma}
A $G$-invariant linear connection of a homogeneous $G$-supermanifold 
is flat if and only if the associated Nomizu map 
is a Lie superalgebra morphism.
\end{lemma}
In the reductive case, for every $A,B\in\gm$
$$
R_{o}(A,B)=[L_{\gm}(A),L_{\gm}(B)]-L_{\gm}([A,B]_{\gm})-[[A,B]_{\gh},\cdot]\quad,
$$
$$
T_{o}(A,B)=L_{\gm}(A)B-(-1)^{|A||B|}L_{\gm}(B)A-[A,B]_{\gm}\quad\quad\phantom{c}.
$$
A $G$-invariant linear connection preserves a $G$-invariant tensor field $T\in\mathcal{T}_{s}^{r}(M)$, \textit{i.e.} $\nabla T=0$, if and only if $L_{p}(\gg)\subseteq \left\{ L\in\ggl_{\bR}(T_{p}M)|L\cdot T|_{p}=0\right\}$.
\begin{definition}
\rm{Let $M=G/H$ be a homogenous manifold, $\gk\subseteq\ggl_{\bR}(\gg/\gh)$ a linear Lie superalgebra. A $G$-invariant linear connection is {\it $\gk$-compatible} if $L(\gg)\subseteq\gk$.
}\end{definition}
\begin{example}
\rm{
Let $G/H$, $\gg=\gh+\gm$ be a homogeneous supermanifold with reductive decomposition (\ref{reductived}). The Levi-Civita connection $\nabla^{LC}$ (see \cite{MS} for its existence and unicity) of an invariant metric associated with an $H$-invariant scalar product $g:\gm\times\gm\rightarrow\bR$ is given by
$$
L_{\gm}(A)B=\frac{1}{2}[A,B]_{\gm}+U(A,B)
$$
where $U:\gm\times\gm\rightarrow\gm$ is the supersymmetric bilinear map given by ($A,B,C\in\gm$)
$$
2g(U(A,B),C):=(-1)^{|B||C|}g(A,[C,B]_{\gm})+(-1)^{|C|(|A|+|B|)}g([C,A]_{\gm},B)\qquad.
$$}
\end{example}
In our context, Theorems 5.1 and 8.1 of \cite{Ga} read as follows. 
\begin{theorem}[\cite{Ga}]
\label{holonomycorresp}
Let $\nabla:\mathcal{T}_{M}\otimes_{\bR}\mathcal{T}_{M}\rightarrow \mathcal{T}_{M}$ be a $G$-invariant linear connection on a simply connected homogeneous supermanifold $M=G/H$. Then there exists a natural bijective correspondence between 
\begin{itemize}
\item[$1)$] Parallel tensors $T\in\mathcal{T}_{s}^{r}(M)$,
\item[$2)$] Tensors $T_{o}\in T_{o}M^{r}_{s}$ annihilated by the representation of the infinitesimal holonomy algebra $\mathfrak{hol(\nabla)_{o}^{inf}}$ on $T_{o}M^{r}_{s}$.
\end{itemize}
\end{theorem}
We calculate the infinitesimal holonomy algebra of a $G$-invariant connection.
\begin{theorem}
The infinitesimal holonomy algebra of a $G$-invariant linear connection $\nabla$ on a homogeneous supermanifold $M=G/H$ is given by
$$
\mathfrak{hol(\nabla)_{o}^{inf}}=\gr+[L(\gg),\gr]+[L(\gg),[L(\gg),\gr]]+\cdot\cdot\cdot
$$
where $L:\gg\rightarrow\ggl_{\bR}(\gg/\gh)$ is the Nomizu map associated with $\nabla$ and 
$$
\gr:=\rm{\Span}_{\bR}\left\{[L(A),L(B)]-L([A,B])\,|\,A,B\in\gg\right\}\subseteq\ggl_{\bR}(\gg/\gh)\quad.
$$
\end{theorem}
\begin{pf}
Consider
$
\nabla^{r}_{A_{r},..,A_{1}}R_{o}(A,B)\in\ggl_{\bR}(\gg/\gh)
$
for $A_{r},..,A_{1},A,B\in\gg$. The first two cases are illustrated, the other ones follows similarly.
\\
For $r=0$, $R_{o}(A,B)=[L(A),L(B)]-L([A,B])\in\gr$.
\\
For $r=1$,
$$
\nabla^{1}_{A_{1}}R_{o}(A,B)=L(A_{1})\circ(R_{o}(A,B))-R_{o}(L(A_{1})A,B)-(-1)^{|A||A_{1}|}R_{o}(A,L(A_{1})B)
$$
$$
\phantom{cccccccccccccccc}-(-1)^{|A_{1}|(|A|+||B||)}R_{o}(A,B)\circ L(A_{1})=[L(A_{1}),R_{o}(A,B)]\quad\rm{mod}\,\gr\quad.
$$
\end{pf}
The proof of the following corollary is as in \cite{KN2} and thus omitted.
\begin{corollary}
\label{holonomyreductive}
The infinitesimal holonomy algebra of a $G$-invariant linear connection $\nabla$ on a homogeneous supermanifold $M=G/H$ with reductive decomposition (\ref{reductived}) is given by
$$
\mathfrak{hol(\nabla)_{o}^{inf}}=\gr+[L_{\gm}(\gm),\gr]+[L_{\gm}(\gm),[L_{\gm}(\gm),\gr]]+\cdot\cdot\cdot
$$
where $L_{\gm}:\gm\rightarrow\ggl_{\bR}(\gm)$ is the Nomizu map associated with $\nabla$ and
$$
\gr:=\rm{\Span}_{\bR}\left\{[L_{\gm}(A),L_{\gm}(B)]-L_{\gm}([A,B]_{m})-[[A,B]_{\gh},\cdot]\,|\,A,B\in\gm\right\}\subseteq\ggl_{\bR}(\gm)\quad.
$$
\end{corollary}
In particular, if $\nabla^{can}$ is the canonical connection
and $M_{0}$ is simply connected, 
every invariant tensor field is parallel 
(and viceversa if $[\gm,\gm]=\gh$).

\section{Superizations of a homogeneous spin manifold 
}
\label{sezione5}
\setcounter{equation}{0}
\subsection{Homogeneous spin manifold}\hfill\newline\\
Let $M_{0}=G_{0}/H$ be a (pseudo)-Riemannian reductive homogeneous space. The Lie algebra $\gg_{0}$ of the Lie group $G_{0}$ splits into $\gg_{0}=\gh+\gm_{\0}$, $\gh$ and $\gm_{\0}$ being the Lie algebra of the Lie group $H$ and a $\Ad_{H}$-invariant complement respectively. The (pseudo)-Riemannian metric on $M_{0}$ is identified with a $H$-invariant inner product $g\in\Bil_{\bR}(\gm_{\0})^{H_{0}}$ on $\gm_{\0}$.
Recall that 
$
\varphi_{0}:\gg_{0}\rightarrow\mathcal{T}(G_{0})
$
is the realization of the (abstract) Lie algebra $\gg_{0}$ as left-invariant vector fields of $G_{0}$.
It is known that the $H$-principal bundle $\pi_{0}:G_{0}\rightarrow M_{0}$
is a (non-canonical) reduction of the bundle of orthonormal frames of $M_{0}$.
The tangent bundle of $M_{0}$ is then given by
$$
G_{0}\times_{\Ad(H)}\varphi_{0}(\gm_{\0})\cong T(M_{0})
$$
where an element $[g,\varphi_{0}(w)]\in G_{0}\times_{\Ad(H)}\varphi_{0}(\gm_{\0})$ corresponds to a vector via
$$
G_{0}\times_{\Ad(H)}\varphi_{0}(\gm_{\0})\ni [g,\varphi_{0}(w)]\mapsto (\pi_{0})_{*,g}(\varphi_{0}(w)|_{g})\in T_{gH}(M_{0})\quad.
$$
Hence a vector field on $M_{0}$ is an $H$-equivariant map
$
X:G_{0}\rightarrow\varphi_{0}(\gm_{\0})
$,
\textit{i.e.} a map $X\in\Map(G_{0},\varphi_{0}(\gm_{\0}))$ satisfying the equivariance property
$$
X(g)=\Ad_{h}X(gh)
$$
for all $g\in G_{0}$, $h\in H$. In other words, the following identifications hold
\be
\label{vargas}
[\mathcal{C}^{\infty}(G_{0})\otimes\varphi_{0}(\gm_{\0})]^{R_{H_{0}}}\cong\Map(G_{0},\varphi_{0}(\gm_{\ol{0}}))^{H_{0}}\cong\mathcal{T}(M_{0})
\ee
These identifications are a special case of Lemma \ref{utile}.
\begin{definition}[\cite{Bar, Fr}]
\rm{
A \textbf{homogeneous spin structure} of a (pseudo)-Riemannian reductive homogeneous space $(M_{0}=G_{0}/H,g)$ is a homomorphism $$\tilde{\Ad}:H\rightarrow\Spin(\gm_{\0})$$
lifting the isotropy representation to $\Spin(\gm_{\0})$, \textit{i.e.} such that the diagram  
\begin{equation*}
\begin{CD}
H @>\tilde{\Ad}>> \Spin(\gm_{\0}) \\
@V \Id_{H} VV               @VV \xi V    \\
H     @>\Ad>>    \SO(\gm_{\0}) 
\end{CD}
\end{equation*}
commutes. The triple $(M_{0},g,\tilde{\Ad})$ is called a \textbf{homogeneous spin manifold}. The associated {\bf spin bundle} is defined by
$$S(M_{0}):=G_{0}\times_{\tilde{\Ad}(H)}S
$$
where $S$ is an irreducible real representation $$\Delta:\mathrm{Cl}(\gm_{\0})\rightarrow\End_{\bR}(S)$$ of the Clifford algebra $\mathrm{Cl}(\gm_{\0})\supseteq\Spin(\gm_{\0})$.
}\end{definition}
If $G_{0}$ is simply connected, the homogeneous spin structures of $M_{0}$ are in one-to-one correspondence with the spin structures of $M_{0}$ (\cite{Bar}). This condition is tacitly assumed together with the connectedness of $H$.\\
Henceforth denote the spin module $S$ by $\gm_{\ol{1}}$.
A spinor field on $M_{0}$ is an $H$-equivariant map
$
\psi:G_{0}\rightarrow\gm_{\ou}
$,
\textit{i.e.} a map $\psi\in\Map(G_{0},\gm_{\ou})$ satisfying
\be
\label{equi2}
\psi(g)=\Delta\circ\tilde{\Ad}_{h}\psi(gh)
\ee
for all $g\in G_{0}$, $h\in H$. Write $\psi\in\mathcal{S}(M_{0})=\Map(G_{0},\gm_{\ou})^{H}$. Let $X,Y\in\Map(G_{0},\varphi_{0}(\gm_{\0}))^{H}$ (resp. $\psi,\zeta\in\mathcal{S}(M_{0})$) be vector (resp. spinor) fields of a homogeneous spin manifold $(M_{0},g,\tilde{\Ad})$. 
The canonical connection derivatives $$\nabla^{can}_{X}Y\in\Map(G_{0},\varphi_{0}(\gm_{\0}))^{H}\quad,\quad\nabla^{can}_{X}\psi\in\mathcal{S}(M_{0})$$ 
are the usual derivatives along $X$, \textit{i.e.}
$$
(\nabla^{can}_{X}Y)|_{g}=(\partial_{X|_{g}}Y)\quad,\quad(\nabla^{can}_{X}\psi)|_{g}=(\partial_{X|_{g}}\psi)
$$
for every $g\in G_{0}$. Recall the following important definition.
\begin{definition}[\cite{K}]
\rm{The {\it Kosmann Lie derivative} $\mathcal{L}_{X}\psi\in \cS(M_{0})$ is defined by
$$
\label{Kosmann1}
\mathcal{L}_{X}\psi:=\nabla_{X}^{LC}\psi-\cA(\nabla^{LC}X)\cdot\psi
$$
where $\cA(\nabla^{LC}X)\in\so(T_{p}M_{0},g)\cong\spin(T_{p}M_{0},g)$ is the alternation of the Levi-Civita covariant derivative $\nabla^{LC}X\in\ggl_{\bR}(T_{p}M_{0})$ of $X$.
}\end{definition}
The following definition gives extra-structures on $(M_{0},g,\tilde{\Ad})$.
\begin{definition}\rm{
Every $\gh$-invariant bilinear form $C:\gm_{\0}\otimes\gm_{\ou}\rightarrow\gm_{\ou}$ defines a {\bf generalized Clifford multiplication}
$$
\mathcal{T}(M_{0})\otimes\cS(M_{0})\rightarrow\cS(M_{0})
$$ 
$$
\phantom{cccccccccccccccccccccccccc}X\otimes\psi\mapsto C(X,\psi):=(g\rightarrow C(X(g),\psi(g)))
$$
where $g\in G_{0}$. A spinor field $\psi\in\cS(M_{0})$ is a {\bf generalized Killing spinor} if
\be
\label{killingequation}
\nabla^{\cS}_{X}\psi:=\nabla^{can}_{X}\psi+C(X,\psi)=0
\ee
where $X\in\mathcal{T}(M_{0})$.
Every $\gh$-invariant symmetric bilinear form $\Gamma:\gm_{\ou}\otimes\gm_{\ou}\rightarrow\gm_{\0}$ defines a {\bf Dirac current bracket}
$$
\cS(M_{0})\otimes\cS(M_{0})\rightarrow\mathcal{T}(M_{0})
$$ 
$$
\phantom{cccccccccccccccccccccc}\psi\otimes\zeta\mapsto\Gamma(\psi,\zeta):=(g\rightarrow\Gamma(\psi(g),\zeta(g))
$$
where $g\in G_{0}$.
}\end{definition}
The generalized Clifford multiplication $C$ and the Dirac bracket $\Gamma$ are some of the objects involved in Definition \ref{supersymmetric}.
\subsection{Algebra of supersymmetry}\hfill\newline\\
\label{quelladopo}
Recall that to a pseudo-Riemannian spin manifold $(M_{0},g,S)$ is naturally associated a split supermanifold $M=(M_{0},\cA_{M})$ whose sheaf of superfunctions $\cA_{M}$ is the exterior algebra of the (dual) of the spin bundle $S$, \textit{i.e.}
\be
\label{sheafsheaf}
\cA(M)\cong\Lambda(\cS^{*}(M_{0}))
\ee
Supermanifolds of this type have been studied in \cite{ACDS, K1, K2}. In our setting, it is natural to ask whether (\ref{sheafsheaf}) can be endowed with a structure of homogeneous supermanifold.
\begin{definition}
\label{superizationhomogeneous}
\rm{
A \textbf{superization} of a homogeneous spin manifold 
\be
\label{quasitre}
(M_{0}=G_{0}/H,g,\tilde{\Ad})
\ee
with reductive decomposition $\gg_{0}=\gh+\gm_{\0}$, is a homogeneous structure on the supermanifold $M=(M_{0},\cA_{M})$ whose sheaf of superfunctions is given by (\ref{sheafsheaf}). 
}\end{definition}
Constructing a superization can be reduced to the following algebraic problem. Recall that the Lie morphism induced by a homogeneous spin structure is given by $$\tilde{\ad}:=\xi_{*}^{-1}\circ\ad:\gh\rightarrow\spin(\gm_{\0})\quad.$$
\begin{definition}
\label{supersymmetric}
\rm{
A Lie superalgebra $(\gg,[\cdot,\cdot])$, where
\be
\label{algsuper}
\gg=\gg_{\ol{0}}+\gg_{\ol{1}}:=(\gh+\gm_{\ol{0}})+\gm_{\ol{1}}
\ee and the adjoint action of $\gh$ on $\gm_{\ol{1}}$ is given by 
$$
\label{azione}
\Delta\circ\tilde{\ad}:\gh\rightarrow\ggl_{\bR}(\gm_{\ou})
$$
is called an {\bf algebra of supersymmetry} adapted to the spin manifold $(M_{0},g,\tilde{\Ad})$.
}\end{definition}
\begin{theorem}
Every algebra of supersymmetry (\ref{algsuper}) adapted to a homogeneous spin manifold (\ref{quasitre}) defines a superization $M=G/H$ of (\ref{quasitre}).
\end{theorem}
\begin{proof}
Consider the unique sHC pair associated with the algebra of supersymmetry (\ref{algsuper}) ($\pi_{1}(G_{0})=\left\{0\right\}$ is used to prove the existence of the adjoint action of the sHC pair) and denote by $G$ the associated Lie supergroup via Theorem \ref{kostantone}. The homogeneous supermanifold $G/H$ is a superization of $G_{0}/H$. Indeed, using the Koszul realization of the sheaf $\cA_{G}\cong\cC^{\infty}_{G_{0}}\otimes\Lambda(\gg_{\ou}^{*})$, for every superfunction 
$$
f\in\Hom_{\bR}(\Lambda(\gg_{\ou}),\cC^{\infty}(G_{0}))
$$
the right action of the Lie group $H\ni h$ is given by
$$
(R_{h}^{*}f)(a)(g)=f(\Ad_{h}^{-1}a)(gh)
$$
where $g\in G_{0}$ and $a\in\Lambda(\gg_{\ou})$. Then the structure sheaf $\cA_{G/H}$ satisfies (\ref{sheafsheaf}).
\end{proof}
Henceforth every superization is tacitly assumed to come from an adapted algebra of supersymmetry (\ref{algsuper}) and the notation
$$C:=[\cdot,\cdot]|_{\gm_{\0}\otimes\gm_{\ou}}:\gm_{\0}\otimes\gm_{\ou}\rightarrow\gm_{\ou}\phantom{ccc},\phantom{ccc}\Gamma:=\pi_{\gm_{\0}}\circ [\cdot,\cdot]|_{\gm_{\ou}\otimes\gm_{\ou}}:\gm_{\ou}\otimes\gm_{\ou}\rightarrow\gm_{\0}$$
is used. Note that the adjoint representation $\Ad:G_{0}\rightarrow\Aut(\gg)$ of the sHC pair $(G_{0},\gg)$ associated with an adapted algebra of supersymmetry satisfies
$$
\Ad_{h}(a)=\Delta\circ\tilde{\Ad}_{h}(a)
$$
where $h\in H$ and $a\in\gm_{\ol{1}}$. 
Thus equation (\ref{equi2}) can be re-written as
$$
\psi(g)=\Ad_{h}\psi(gh)
$$
where $g\in G_{0}$ and $h\in H$. In other words, the identifications
$$
[\mathcal{C}^{\infty}(G_{0})\otimes\gm_{\ou}]^{R_{H}}\cong\Map(G_{0},\gm_{\ol{1}})^{H}\cong \mathcal{S}(M_{0})
$$
are the analogues of ($\ref{vargas}$) and show that vector and spinor fields can be treated equally inside the superization $M=G/H$.
\begin{lemma}
\label{importante2}
The tangent bundle $TM$ of a superization $M=G/H$ is isomorphic to the direct sum of the tangent bundle $T(M_{0})$ with the spin bundle $S(M_{0})$ of the body $M_{0}=G_{0}/H$. Indeed the projections (\ref{even}) and (\ref{odd}) are given by 
$$
\label{evenpro}
\,\ev_{\ol{0}}:[\mathcal{A}(G)\otimes\varphi(\gm)]^{R_{H}}\rightarrow[\mathcal{C}^{\infty}(G_{0})\otimes\varphi_{0}(\gm_{\ol{0}})]^{R_{H}}\quad,
$$
$$
\label{oddpro}
\ev_{\ol{1}}:[\mathcal{A}(G)\otimes\varphi(\gm)]^{R_{H}}\rightarrow[\mathcal{C}^{\infty}(G_{0})\otimes\gm_{\ol{1}}]^{R_{H}}\phantom{\rho_{cc0}}\,\quad.
$$
\end{lemma}
\begin{proof}
The assertion follows directly from Lemma \ref{utile} and the following remark. Equation (\ref{tangentedispari}) implies that the odd-value at a point $gH\in G_{0}/H$ of a vector field
$$
\sum_{i}f^{i}\otimes\varphi(a_{i})\in[\mathcal{A}(G)\otimes\varphi(\gm_{\ol{1}})]^{R_{H}}
$$
is given by
$\sum_{i}f^{i}(g)\otimes\varphi(a_{i})|_{g}\cong f(g)\otimes a_{i}$, \textit{i.e}
$\ev_{\ou}(\sum_{i}f^{i}\otimes\varphi(a_{i}))=\sum_{i}\tilde{f^{i}}\otimes a_{i}$.
\end{proof}
\subsection{Even and odd Killing fields}\hfill\newline\\
\label{domani}
In this subsection we assume, for simplicity, that $(M_{0},g,\tilde{\Ad})$ is a (pseudo)-Riemannnian symmetric spin manifold. In this case, the Levi-Civita connection 
coincides 
with the canonical connection. Fix an adapted algebra of supersymmetry (\ref{algsuper}) and denote by $M=G/H$ the associated superization of $M_{0}=G_{0}/H$. The space of spinor fields which satisfy the generalized Killing equation (\ref{killingequation}) is denoted by
\be
\label{killingequation2}
\mathcal{KS}:=\left\{\psi\in\cS(M_{0})\,|\,\nabla^{\cS}_{X}\psi=\nabla^{LC}_{X}\psi+C(X,\psi)=0\right\}
\ee
Similarly the image of the representation of the Lie algebra $\gg_{\0}$ 
$$
\hat{\varphi}_{0}:\gg_{\0}\rightarrow\mathcal{T}(M_{0})
$$
$$
\phantom{cccccccccccccccccccccccccccccc}x\mapsto\hat{\varphi}_{0}(x)=(g\mapsto \varphi_{0}((\Ad_{g^{-1}}x)_{\gm_{\0}}))
$$
by Killing vector fields of $M_{0}$ is denoted by
$
\mathcal{KV}:=\left\{\hat{\varphi}_{0}(a)|a\in\gg_{\0}\right\}
$.
\begin{definition}\rm{
The anti-homomorphism of Lie superalgebras $$\hat{\varphi}:\gg\rightarrow\Der_{\bR}(\Lambda(\cS(M_{0}^{*})))$$ given by (the passage to the quotient of) the map (\ref{representation2}) is the {\bf Killing representation} of the algebra of supersymmetry $\gg=\gg_{\0}+\gg_{\ou}$ on $\Lambda(\cS^{*}(M_{0}))$. For every $a\in\gg$, the fundamental vector field $\hat{\varphi}(a)$ is called a {\bf Killing field}.
}\end{definition} 
Recall that the Killing representation of the algebra of supersymmetry is explicitly described by Proposition \ref{parikill} and Theorem \ref{Andrea}. The following Theorem is the main result of this subsection. In particular it shows that our construction is consistent with the {\it symmetry superalgebra construction} in the theoretical physics literature (see all the references cited in Example \ref{esempioni}). The interpretation of generalized Killing spinors (\ref{killingequation2}) as odd-values of the odd Killing fields (\ref{Andreadestra}) is new. Moreover, vector fields and spinor fields are embedded into a bigger space, naturally endowed with a Lie superalgebra structure, namely $\mathcal{T}(M)\cong\Der_{\bR}(\Lambda(\cS^{*}(M_{0})))$.
\begin{theorem}
\label{differentialspinor}
Let $(M_{0},g,\tilde{\Ad})$ be a (pseudo)-Riemannian symmetric spin manifold $M_{0}=G_{0}/H$ together with an adapted algebra of supersymmetry $\gg=\gg_{\0}+\gm_{\ou}=(\gh+\gm_{\0})+S$ and let $M=G/H$ be the associated superization such that $\cA(M)\cong\Lambda(\cS^{*}(M_{0}))$. For every $x\in\gg_{\0}$ and $s\in S$, the value of the associated Killing field on $M$ is given by
\be
\label{quasifine}
\ev_{\0}\circ\hat{\varphi}(x)=\hat{\varphi}_{0}(x)\in\mathcal{KV}\subseteq\mathcal{T}(M_{0})
\ee
\be
\label{quasifine2}
\,\ev_{\ou}\circ\hat{\varphi}(s)=\psi^{s}\in\mathcal{KS}\subseteq\mathcal{S}(M_{0})\phantom{cccc}
\ee
where 
\be
\label{figu}
\psi^{s}:=(g\mapsto\Ad_{g^{-1}}s)
\ee
is the (unique) generalized Killing spinor with value $s\in S$ at the point $o\in M_{0}$. The direct sum of (\ref{quasifine}) together with (\ref{quasifine2}) is an isomorphism of Lie superalgebras
\be
\label{geometry}
\ev:\hat{\varphi}(\gg)\rightarrow\mathcal{KV}+\mathcal{KS}
\ee
where the structure of Lie superalgebra of $\mathcal{KV}+\mathcal{KS}$ is defined through the Kosmann Lie derivative and the Dirac current bracket via
$$
[\hat{\varphi}_{0}(x),\psi^{s}]:=\mathcal{L}_{\hat{\varphi}_{0}(x)}\psi^{s}\qquad,\qquad[\psi^{s},\psi^{t}]:=-\Gamma(\psi^{s},\psi^{t})
$$
where $x\in\gg_{\0}$ and $s,t\in S$. Moreover, there exists a canonical embedding
\be
\label{embeddingsuper}
\mathcal{F}:\mathcal{T}(M_{0})\oplus\mathcal{S}(M_{0})\hookrightarrow\Der_{\bR}(\Lambda(\cS^{*}(M_{0})))
\ee
$$
X+\psi\mapsto\mathbf{X}+\mathbf{\Psi}
$$
such that 
\begin{itemize}
\item[i)] The linear map $\mathcal{T}(M_{0})\ni X\mapsto\mathbf{X}\in\Der_{\bR}(\Lambda(\cS^{*}(M_{0})))$ is a morphism of Lie superalgebras, 
\item[ii)] The linear map (\ref{embeddingsuper}) satisfies $(\ev_{\0}+\ev_{\ou})\circ\mathcal{F}=\Id$,
\item[iii)]
$\ev_{\0}\circ[\mathbf{X},\mathbf{Y}]\,\,=[X,Y]\in\mathcal{T}(M_{0})$,
\item[iv)]
$\ev_{\ou}\circ[\mathbf{X},\mathbf{\Psi}]\,\,=\nabla^{\cS}_{X}\psi\in\mathcal{S}(M_{0})
$,
\item[v)]
$\ev_{\0}\circ[\mathbf{\Psi},\mathbf{\Psi^{'}}]=\Gamma(\psi,\psi^{'})\in\mathcal{T}(M_{0})$,
\end{itemize}
for every $X,Y\in\mathcal{T}(M_{0})$ and $\psi,\psi^{'}\in\mathcal{S}(M_{0})$.
\end{theorem}
\begin{proof}
Equations (\ref{quasifine}) and (\ref{quasifine2}) are direct consequences of Proposition \ref{parikill}, Theorem \ref{Andrea} and related comments
(\textit{e.g.} equation (\ref{andiamoeven})).
Note that $\hat{\varphi}(s)$ and $\psi^{s}$ bijectively correspond to each other. Equation (\ref{killingequation}) is obviously satisfied,
we differentiate equation (\ref{figu}) along the vector field 
$$
X=\sum_{i}f^{i}\otimes\varphi_{0}(a_{i})\in[\mathcal{C}^{\infty}(G_{0})\otimes\varphi_{0}(\gm_{\ol{0}})]^{R_{H}}
$$
so that
$$
\partial_{X}\psi^{s}=
\sum_{i}f^{i}\otimes\frac{d}{dt}|_{t=0}\psi^{s}(g\cdot \exp(ta_{i}))=\sum_{i}f^{i}\otimes\frac{d}{dt}|_{t=0}\Ad_{\exp(-ta_{i})}\circ\Ad_{g^{-1}}(s)
$$
$$
=-\sum_{i}f^{i}\otimes[a_{i},\Ad_{g^{-1}}(s)]=-\sum_{i}f^{i}\otimes[a_{i},\psi^{s}(g)]=-C(X,\psi^{s})\phantom{ccc}
$$
We prove that the bijective map (\ref{geometry}) is a Lie superalgebra morphism in the Riemannian case, the pseudo-Riemannian case being analogous.
Fix an orthonormal linear frame $\left\{ a_{i}\right\}_{i=1}^{\dim M_{0}}\subseteq \gm_{\0}$, then
$$
\mathcal{L}_{X}\psi=\sum_{i}f^{i}\partial_{\varphi_{0}(a_{i})}\psi-\frac{1}{8}\sum_{i,j}(\partial_{\varphi_{0}(a_{i})}f^{j}-\partial_{\varphi_{0}(a_{j})}f^{i})a_{i}\cdot a_{j}\cdot\psi\quad.
$$
If $X$ is Killing, the previous formula reduces to
$$
\mathcal{L}_{X}\psi=\sum_{i}f^{i}\partial_{\varphi_{0}(a_{i})}\psi+\frac{1}{4}\sum_{i,j}(\partial_{\varphi_{0}(a_{i})}f^{j})a_{j}\cdot a_{i}\cdot\psi
$$
In particular, when 
$$X=\hat{\varphi}_{0}(x)=\left\{g\mapsto \varphi_{0}((\Ad_{g^{-1}}x)_{\gm_{\0}})\right\}\qquad,\qquad\psi=\psi^{s}$$
we get that
$$
\mathcal{L}_{\hat{\varphi}_{0}(x)}\psi^{s}=-\sum_{i}C((\Ad_{g^{-1}}x)_{\gm_{\0}},\Ad_{g^{-1}}s)+\frac{1}{4}\sum_{i}(\partial_{\varphi_{0}(a_{i})}\hat{\varphi}_{0}(x))\cdot a_{i}\cdot\psi^{s}
$$
$$
\phantom{ccccccccc}=-\sum_{i}C((\Ad_{g^{-1}}x)_{\gm_{\0}},\Ad_{g^{-1}}s)-\frac{1}{4}\sum_{i}[a_{i},\Ad_{g^{-1}}x]_{\gm_{\0}}\cdot a_{i}\cdot\psi^{s}
$$
$$
\phantom{ccccccccc}=-\sum_{i}C((\Ad_{g^{-1}}x)_{\gm_{\0}},\Ad_{g^{-1}}s)-\frac{1}{4}\sum_{i}[a_{i},(\Ad_{g^{-1}}x)_{\gh}]\cdot a_{i}\cdot\psi^{s}
$$
$$
=-\sum_{i}[(\Ad_{g^{-1}}x)_{\gm_{\0}},\Ad_{g^{-1}}s]-[(\Ad_{g^{-1}}x)_{\gh},\Ad_{g^{-1}}s]\phantom{c}
$$
$$
\phantom{cc}=-[\Ad_{g^{-1}}x,\Ad_{g^{-1}}s]=-\Ad_{g^{-1}}[x,s]=-\psi^{[x,s]}\qquad\phantom{ccccc}
$$
It then follows that
$$
\ev_{\ou}\circ[\hat{\varphi}(x),\hat{\varphi}(s)]=-\ev_{\ou}\circ\hat{\varphi}[x,s]=-\psi^{[x,s]}=\mathcal{L}_{\hat{\varphi}(x)}\psi^{s}
$$
for every $x\in\gg_{\0}$ and $s\in\gm_{\ou}=S$. The equations
$$
\ev_{\0}\circ[\hat{\varphi}(x),\hat{\varphi}(y)]=[\hat{\varphi}_{0}(x),\hat{\varphi}_{0}(y)]\qquad,\qquad
\ev_{\0}\circ[\hat{\varphi}(s),\hat{\varphi}(t)]=-\Gamma(\psi^{s},\psi^{t})
$$
for every $x,y\in\gg_{\0}$ and $s,t\in\gm_{\ou}=S$ imply that (\ref{geometry}) is an isomorphism of Lie superalgebras. The embedding (\ref{embeddingsuper}) is defined by
$$
[\mathcal{C}^{\infty}(G_{0})\otimes\varphi_{0}(\gm_{\0})]^{R_{H}}\ni X=\sum_{i}f^{i}\otimes\varphi_{0}(a_{i})\mapsto\sum_{i}f^{i}\otimes\varphi(a_{i})=:\mathbf{X}\in[\mathcal{A}(G)\otimes\varphi(\gm)]^{R_{H}}
$$
$$
[\mathcal{C}^{\infty}(G_{0})\otimes\gm_{\ou}]^{R_{H}}\ni\psi=\sum_{i}f^{i}\otimes a_{i}\mapsto \sum f^{i}\otimes\varphi(a_{i})=:\Psi\in[\mathcal{A}(G)\otimes\varphi(\gm)]^{R_{H}}\phantom{cccccc}
$$
The proof of the remaining properties is straightforward using the remarks after Lemma \ref{utile}.
\end{proof}
The following example clarifies the situation.
\begin{example}
\label{riprendo}
\rm{The notation is the one of Example \ref{examplepoincarè} and Example \ref{spst}. The even-value of the odd Killing field
$$
\hat{\varphi}(s_{\beta})=y^{\alpha}_{\beta}(g^{-1})[-\frac{\partial}{\partial s^{\alpha}}+
\frac{1}{2}s^{\eta}\varphi_{0}(\Gamma_{\alpha\eta}^{k}e_{k})]\in\mathcal{T}(M)
$$
is obviously zero. The equality 
$$
\hat{\varphi}(s_{\beta})=y^{\alpha}_{\beta}(g^{-1})[\varphi(s_{\alpha})+s^{\eta}\varphi_{0}(\Gamma_{\alpha\eta}^{k}e_{k})]\qquad \mod\phantom{c}\mathcal{A}(G)\otimes\gh
$$
implies that the odd-value is the parallel spinor 
$$
y^{\alpha}_{\beta}(g^{-1})s_{\alpha}\in\cS(M_{0})
$$
whose image under the embedding (\ref{embeddingsuper}) is given by
$$
y^{\alpha}_{\beta}(g^{-1})\varphi(s_{\alpha})=y^{\alpha}_{\beta}(g^{-1})[-\frac{\partial}{\partial s^{\alpha}}-\frac{1}{2}s^{\eta}\varphi_{0}(\Gamma_{\alpha\eta}^{k}e_{k})]\in\mathcal{T}(M)
\quad.
$$
}\end{example}
The following example describes the symmetry superalgebras of the maximally supersymmetric solutions of (bosonic) 11-dimensional supergravity.
\begin{example}[\cite{CJS, F, F1, F2, FP, HKS}]
\label{esempioni}
\label{max1}
\rm{
A Lorentzian spin manifold $(M_{0},g,\mathcal{S})$ together with a closed {\it flux} $F\in\Lambda^{4}(M_{0})$ is a bosonic solution of 11-dimensional {\bf supergravity} if
$$
\label{Ein}
\Ric(X,Y)-\frac{1}{2}sg(X,Y)=-\frac{1}{2}g(i_{X}F,\i_{Y}F)+\frac{1}{6}g(X,Y)|F|^{2}\qquad\quad\textbf{(Einstein)}
$$
$$
\label{Max}
d*F=-\frac{1}{2}F\wedge F\phantom{ccccccccccccccccccccccccccccccccccccccccc}\qquad\quad\textbf{(Maxwell)}
$$
for all $X,Y\in\mathcal{T}(M_{0})$. A spinor field $\psi\in\mathcal{S}(M_{0})$ is {\bf supergravity-Killing} if
$$
\label{superkill}
\nabla_{X}^{\cS}\psi:=\nabla^{LC}_{X}\psi+(-\frac{1}{12}X\wedge F^{\sharp}+\frac{1}{6}(\i_{X}F)^{\sharp})\cdot\psi=0
$$
for every $X\in\mathcal{T}(M_{0})$. The spin connection $\nabla^{\cS}$ is the main object of the theory. Indeed the Einstein and Maxwell equations are equivalent to
$$
\label{curv}
\sum_{i}e_{i}^{*}\cdot R^{\nabla^{\cS}}(X,e_{i})=0
$$
where $X\in\mathcal{T}(M_{0})$, $\left\{e_{i}\right\}$ is a local pseudo-orthonormal frame of $M_{0}$ and $\cdot$ denotes Clifford multiplication.
The symmetry superalgebras $\gg=\gg_{\0}+\gg_{\ou}=(\gh+\gm_{\0})+S$ of maximally supersymmetric solutions, \textit{i.e.} plane-wave, Freund-Rubin and flat backgrounds, are recalled. Note that the spin connection $\nabla^{\cS}$ is encoded in the bracket $[\gm_{\0},S]$. Theorem \ref{differentialspinor} implies that the supergravity Killing spinors are given by the odd-value (\ref{figu}) of the odd-Killing fields of the supermanifold $\Lambda(\mathcal{S}^{*}(M_{0}))$. 
\\ 
\\
\textit{Plane-wave.}
Let $\bR^{1,10}$ be the vector space $\bR^{11}$, together with the standard inner product $\left\langle\cdot ,\cdot\right\rangle $ of signature $(1,10)=(+,-)$. Fix a Witt decomposition
$$\bR^{1,10}=(\bR\mathrm{p} \oplus\bR\mathrm{q})\bigoplus E=(\bR\mathrm{p}\oplus\bR\mathrm{q})\bigoplus(\oplus_{i=1}^{9}\bR e_{i})$$
with $\left\langle e_{i} , e_{j}\right\rangle=-\delta_{ij}$, $\left\langle \mathrm{p},E\right\rangle=\left\langle \mathrm{q},E\right\rangle=\left\langle \mathrm{p},\mathrm{p}\right\rangle=\left\langle \mathrm{q},\mathrm{q}\right\rangle=0$ and $\left\langle \mathrm{p},\mathrm{q}\right\rangle=1$.
The non-trivial Lie brackets of the Cahen-Wallach symmetric space
$$
\label{CWfig}
\gg_{0}=\gh+\gm_{\0}=(E^{*}+\so(3)+\so(8))+\bR^{1,10}
$$
associated with
$$
B=\begin{pmatrix} \frac{1}{9} & 0 & 0 & 0 & 0 & \cdot\cdot\cdot & 0 \\ 0 & \frac{1}{9} & 0 & 0 & 0 & \cdot\cdot\cdot & 0 \\
0 & 0 & \frac{1}{9} & 0 & 0 &\cdot\cdot\cdot & 0 \\ 0 & 0 & 0 & \frac{1}{36} & 0 & \cdot\cdot\cdot & 0 \\
0 & 0 & 0 & 0 & \frac{1}{36} & \cdot\cdot\cdot & 0 \\ 0 & 0 & 0 & 0 & 0 & \cdot\cdot\cdot & 0 \\
0 & 0 & 0 & 0 & 0 & \cdot\cdot\cdot & \frac{1}{36}
\end{pmatrix}\in\End_{\bR}(E)
$$
are given by
\begin{itemize}
\item[$i)$]
$[M_{ij},e_{k}^{*}]=-\delta_{ik}e_{j}^{*}+\delta_{jk}e_{i}^{*}$,
\item[$ii)$]
$[M_{ij},e_{k}]=-\delta_{ik}e_{j}+\delta_{jk}e_{i}$,
\item[$iii)$]
$[e_{i}^{*},e_{i}]=-B(e_{i})\mathrm{p}$,
\item[$iv)$]
$[e_{i}^{*},\mathrm{q}]=-B(e_{i})e_{i}$,
\item[$v)$]
$[\mathrm{q},e_{i}]=-e_{i}^{*}$,
\end{itemize}
where, for $1\leq i,j \leq 3$ or $4\leq i,j \leq 11$, the infinitesimal generators of $\so(3)\oplus\so(8)$  
$$M_{ij}=e_{i}\wedge e_{j}=\left\langle e_{i},\cdot\right\rangle e_{j}-\left\langle e_{j},\cdot\right\rangle e_{i}\in\so(3)\oplus\so(8)$$
satisfy the usual mutual relations.
The flux is the invariant four form associated with
$F=-\mathrm{q}^{*}\wedge e_{1}^{*}\wedge e_{2}^{*}\wedge e_{3}^{*}\in\Lambda^{4}\gm_{0}^{*}$.
The non-trivial even-odd brackets are 
\begin{itemize}
\item[$i)$] 
$[M_{ij},Q_{\pm}]=\frac{1}{2}e_{i}e_{j}\cdot Q_{\pm}$,
\item[$ii)$]
$[e_{i}^{*},Q_{+}]=\frac{1}{18}e_{i}\mathrm{p}\cdot Q_{+}$ if $1\leq i\leq 3$,
\item[$iii)$]
$[e_{i}^{*},Q_{+}]=\frac{1}{72}e_{i}\mathrm{p}\cdot Q_{+}$ if $4\leq i\leq 11$,
\item[$iv)$]
$[\mathrm{q},Q_{+}]=\frac{1}{4}I\cdot Q_{+}\qquad,\qquad[\mathrm{q},Q_{-}]=\frac{1}{12}I\cdot Q_{-}\qquad\qquad\qquad I:=e_{1}e_{2}e_{3}$
\item[$v)$]
$[e_{i},Q_{+}]=\frac{1}{6}Ie_{i}\mathrm{p}\cdot Q_{+}$ if $1\leq i\leq 3$,
\item[$vi)$]
$[e_{i},Q_{+}]=-\frac{1}{12}Ie_{i}\mathrm{p}\cdot Q_{+}$ if $4\leq i\leq 11$.
\end{itemize}
The $32$-dimensional vector space $S$ is decomposed into two 16-dimensional vector spaces $S_{\pm}$ whose elements are denoted by $Q_{\pm}\in S_{\pm}$. The odd-odd bracket is given by 
$$
[Q_{+},Q_{+}]=(Q_{+},\mathrm{p}\cdot Q_{+})\mathrm{q}+\sum_{i,j\leq 3}\frac{1}{6}(Q_{+},Ie_{i}e_{j}\mathrm{p}\cdot Q_{+})M_{ij}+\sum_{4\leq i,j}\frac{1}{12}(Q_{+},Ie_{i}e_{j}\mathrm{p}\cdot Q_{+})M_{ij},
$$
$$
[Q_{+},Q_{-}]=-\sum_{i=1}^{9}(Q_{+},e_{i}\cdot Q_{-})e_{i}-3\sum_{i\leq 3}(Q_{+},Ie_{i}\cdot Q_{-})e_{i}^{*}-6\sum_{4\leq i}(Q_{+},Ie_{i}\cdot Q_{-})e_{i}^{*},\phantom{ccccc}
$$
$$
[Q_{-},Q_{-}]=(Q_{-},\mathrm{q}\cdot Q_{-})\mathrm{p}\quad.\phantom{ccccccccccccccccccccccccccccccccccccccccccccccccccccccccc}
$$
\cite{FP} gives examples of non-maximally supersymmetric plane-wave solutions. The associated Lie superalgebras are algebra of supersymmetries with odd part the $\gh$-submodule $S_{-}$.\\
\\
\textit{Freund-Rubin.} The Freund-Rubin backgrounds $AdS_{4}\times S_{7}$ and $AdS_{7}\times S_{4}$ are (maximally supersymmetric) solutions of 11-dimensional supergravity once the radii of curvature of the two factors are in ratio $2:1$ and the flux is, up to real constant, the volume form of the 4-dimensional factor. We describe the symmetry superalgebra $\gg=\gg_{\0}+\gg_{\ou}$ in the first case (the second is similar). The action of
$$
\gg_{\0}=\so(2,3)\oplus\so(0,8)=(\so(1,3)+\bR^{1,3})\oplus(\so(0,7)+\bR^{0,7})
$$
on $\gg_{\ou}=S$ is given by
\begin{itemize}
\item[$i)$] $[M_{ij},Q]=\frac{1}{2}e_{i}e_{j}\cdot Q$ 
\item[$ii)$] $[v,Q]=-\frac{1}{2}Iv\cdot Q\quad,\quad[w,Q]=\frac{1}{2}Iw\cdot Q\qquad\qquad\qquad I:=\rm{dvol}(AdS_{4})^{\sharp}$
\end{itemize}
where $M_{ij}\in\gh=\so(1,3)\oplus\so(0,7)$, $v\in\bR^{1,3}$, $w\in\bR^{0,7}$. The odd-odd bracket is given by 
$$
[Q,Q]=\sum_{i=1}^{11}(Q,e_{i}\cdot Q)e_{i}+\sum_{i,j}(Q,Ie_{i}e_{j}\cdot Q)M_{ij}\quad.
$$
\\
\textit{Minkowski.}
The symmetry superalgebra of the flat solution $(\bR^{1,10},\left\langle \cdot,\cdot\right\rangle)$, $F=0$, is the Poincare' Lie superalgebra in signature $(1,10)$.
}\end{example}
The following example deals with the superconformal algebra. This could be an indication that our setting can be developed in more general situations.
\begin{example}[\cite{WZ}]
\rm{
The even part of the {\it superconformal algebra $\gg=\gg_{\0}+\gg_{\ou}$ of Wess and Zumino}
is a central extension of $\so(2,4)$:
$$
\gg_{\0}=\bR\cdot 1 + \so(2,4)=\bR\cdot 1 + (V^{'}+\bR\cdot d +\so(1,3) + V)
$$
where $V$ and $V^{'}$ are two copies of Minkowsky spacetime $(\bR^{1,3},\left\langle \cdot,\cdot\right\rangle)$ and the second equality is the usual decomposition of $\so(2,4)$ into infinitesimal generators of special conformal transformations, dilations, rotations and translations of $\bR^{1,3}$. The spin module in signature $S_{2,4}$ is the direct sum of two copies of the $\spin(1,3)$-module
$$
\gg_{\ou}=S_{2,4}=S_{1,3}+S_{1,3}^{'}
$$
on which $\so(1,3)$ acts diagonally via the $\spin(1,3)$-representation. The Schur algebra of the $\so(1,3)$-representation $S_{1,3}$ is generated by a complex structure $J$. The adjoint actions of the central charge $1\in\gg_{\0}$ and of the dilation generator $d\in\gg_{\0}$ on the odd part $\gg_{\ou}$ are diagonal 
$$
[1,(s,s^{'})]=(Js,-Js^{'})\quad,\quad
[d,(s,s^{'})]=(\frac{1}{2}s,-\frac{1}{2}s^{'})
$$
where $s\in S_{1,3}$ and $s^{'}\in S_{1,3}^{'}$, while those of $v\in V$ and $v^{'}\in V^{'}$ are given by
$$
[v,(s,s^{'})]=\frac{1}{\sqrt{2}}(v\cdot s^{'},0)\quad,\quad
[v^{'},(s,s^{'})]=-\frac{1}{\sqrt{2}}(0,v^{'}\cdot s)
$$
where $\cdot$ denotes Clifford multiplication of a vector with a spinor. The odd-odd brackets are defined as follows. \cite{ACDV} defines, for any $0\leq k\leq 4$, an isomorphism 
$$
\Gamma^{k}:\mathrm{Bil}(S_{1,3})^{\so(1,3)}\rightarrow\mathrm{Bil}^{k}(S_{1,3})^{\so(1,3)}
$$
$$
\phantom{cccc}\beta\mapsto \Gamma_{\beta}^{k}
$$
of the vector space of $\so(1,3)$-invariant bilinear forms on $S_{1,3}$ with the vector space of $\so(1,3)$-invariant $\Lambda^{k}\bR^{1,3}$-valued bilinear forms on $S_{1,3}$, where
$$
\left\langle \Gamma^{k}_{\beta}(s\otimes t),v_{1}\wedge\cdot\cdot\cdot\wedge v_{k}\right\rangle:=\sum_{\pi\in S_{k}}\sgn(\pi)\beta(v_{\pi(1)}\cdot\cdot \cdot v_{\pi(k)}\cdot s,t)
$$
for every $s,t\in S_{1,3}$ and $v_{1},...,v_{k}\in \bR^{1,3}$.
The vector space $\mathrm{Bil}(S_{1,3})^{\so(1,3)}$ is two dimensional and an admissible basis (see \cite{AC}) is given by two skew-symmetric bilinear forms $\beta$ and $\beta_{J}:=\beta(J\cdot,\cdot)$. Clifford multiplication of a vector with a spinor is a skew-symmetric (resp. symmetric) operation with respect to $\beta$ (resp. $\beta_{J}$). The odd-odd brackets are given by
$$
\Gamma_{S}=\Gamma^{1}_{r\beta}\cdot:S\vee S\rightarrow V\quad,\quad\Gamma_{S^{'}}=-\Gamma^{1}_{r\beta}:S^{'}\vee S^{'}\rightarrow V^{'}
$$
$$
\Gamma_{S\otimes S^{'}}=\Gamma_{1}+\Gamma_{d}+\Gamma_{\Lambda^{2}}:S\otimes S^{'}\rightarrow\bR\cdot 1+\bR\cdot d+\Lambda^{2}\bR^{1,3}
$$
where
$$
\Gamma_{1}\cong-\frac{3}{2\sqrt{2}}\Gamma^{0}_{r\beta_{J}}\quad,\quad\Gamma_{d}\cong-\frac{1}{\sqrt{2}}\Gamma^{0}_{r\beta}\quad,\quad\Gamma_{\Lambda^{2}}\cong\frac{1}{4\sqrt{2}}\Gamma^{2}_{r\beta}
$$
and $r\in\bR$ is an arbitrary non-zero real constant. The body of the Wess-Zumino Lie supergroup $G$ is given by
$$
G_{0}=\Spin^{0}(2,4)\times\bR\cong SU(2,2)\times \bR
$$
and the Poincare' group $\Spin^{0}(1,3)\ltimes\bR^{1,3}$ is a Lie subgroup. The projection of the odd Killing fields of $G$ on $\bR^{1,3}\subseteq G_{0}$ gives maps
$$
\psi^{s^{'}}:\bR^{1,3}\rightarrow S_{1,3}+S_{1,3}\quad,\quad v\mapsto -\frac{1}{\sqrt{2}}(v\cdot s^{'},0)
$$
where $s^{'}\in S_{1,3}$. These are all the (non-parallel) {\it twistor spinors}, also called {\it conformal Killing spinors} in the literature, of Minkowsky spacetime.
}\end{example}
\section{Appendix
}
\label{lemmone1}
\setcounter{equation}{0}
\subsubsection{Proof of Lemma \ref{appendixA2}}\hfill\newline\\
We prove the Lemma only in the case of a vector field $X\in\mathcal{T}(N)$. The general case is similar.
It is enough to prove that $X=0$ on a coordinate patch $\left\{x^{r},\xi_{s}\right\}$. Fix a point $p\in U\subseteq N_{0}$. By $\bR$-linearity of relation $ii)$, assume that $$Y_{i}|_{p}=\frac{\partial}{\partial\xi_{i}}|_{p}\quad,$$
\textit{i.e.} there exist $f^{r}_{i}, g^{s}_{i}\in\mathcal{A}_{M}(U)$ such that
$$
Y_{i}=\sum_{r=1}^{m}f^{r}_{i}\frac{\partial}{\partial x^{r}}+\sum_{s=1}^{n}g^{s}_{i}\frac{\partial}{\partial\xi_{s}}
$$
where $f^{r}_{i}(p)=g^{s}_{i}(p)=0$ for every $1\leq r\leq m$, $1\leq s\neq i \leq n$ and $g^{i}_{i}(p)=1$.
Every superfunction $f\in\mathcal{A}_{M}(U)$ can be expressed as
$$
f=\sum_{s=0}^{n}\sum_{\alpha_{1}<\cdot\cdot<\alpha_{s}}f_{\alpha_{1}\cdot\cdot\alpha_{s}}\xi_{\alpha_{1}}\cdot\cdot\xi_{\alpha_{s}}
$$
where $f_{\alpha_{1}\cdot\cdot\alpha_{s}}\in\mathcal{C}^{\infty}(U)$. Say that $f$ is {\it zero up to order} $q\in\bN$ if 
$
f_{\alpha_{1}\cdot\cdot\alpha_{s}}=0 
$
whenever $s\leq q$. The local expression of a vector field $X$ is given by
\be
\label{glottide}
X|_{U}=\sum_{r=1}^{m}f^{r}\frac{\partial}{\partial x^{r}}+\sum_{s=1}^{n}g^{s}\frac{\partial}{\partial\xi_{s}}
\ee
where $f^{r}, g^{s}\in \mathcal{A}_{M}(U)$. Say that $X$ is {\it zero up to order} $q\in\bN$ if so are
all $f^{r}$, $g^{s}$. Given two local vector fields $X,Z$, the non-associative operation
$$
Z\cdot X:=\sum_{r=1}^{m}(Zf^{r})\frac{\partial}{\partial x^{r}}+\sum_{s=1}^{n}(Zg^{s})\frac{\partial}{\partial\xi_{s}}
$$
is well-defined ($X$ is as in (\ref{glottide})). It satisfies the relation
$
\mathcal{L}_{Z}X=Z\cdot X-(-1)^{|Z||X|}X\cdot Z
$.
If $Z$ is zero up to order $q$, then $Z\cdot X$ is zero up to order $q$ as well.
Define 
$$
X_{i_{k}\cdot\cdot i_{1}}:=\mathcal{L}_{Y_{i_{k}}}\cdot\cdot\cdot\mathcal{L}_{Y_{i_{1}}}X\quad.
$$
By induction on $q\in\bN$, we prove that if $X$ satisfies i) and ii) for $1\leq k\leq q$ then $X$ is zero up to order $q$.\\
{\bf $(q=0)$}
From i), $\tilde{f^{r}}=\tilde{g^{s}}=0$, \textit{i.e.} $X$ is zero up to order $0$.\\
{\bf $(q=1)$} $X_{i_{1}}=\mathcal{L}_{Y_{i_{1}}}X=Y_{i_{1}}\cdot X\pm X\cdot Y_{i_{1}}$ and evaluating in $p$ 
$$
0=X_{i_{1}}|_{p}=(Y_{i_{1}}\cdot X)|_{p}
$$
because $X$ is zero up to order 0. Then
$$
0=(Y_{i_{1}}\cdot X)|_{p}=\sum_{r=1}^{m}(Y_{i_{1}}f^{r})(p)\frac{\partial}{\partial x^{r}}+\sum_{s=1}^{n}(Y_{i_{1}}g^{s})(p)\frac{\partial}{\partial\xi_{s}}
$$
$$
\phantom{cccccccccccccccccccc}\,=\sum_{r=1}^{m}(\frac{\partial}{\partial\xi_{i_{1}}}f^{r})(p)\frac{\partial}{\partial x^{r}}+\sum_{s=1}^{n}(\frac{\partial}{\partial\xi_{i_{1}}}g^{s})(p)\frac{\partial}{\partial\xi_{s}}\quad,
$$
\textit{i.e.} $X$ is zero up to order 1.\\
{\bf $(q=2)$} $X_{i_{2}i_{1}}=\mathcal{L}_{Y_{i_{2}}}X_{i_{1}}=Y_{i_{2}}\cdot X_{i_{1}}\pm X_{i_{1}}\cdot Y_{i_{2}}$ and evaluating in $p$ 
$$
0=X_{i_{2}i_{1}}|_{p}=(Y_{i_{2}}\cdot X_{i_{1}})|_{p}
$$
because $X_{i_{1}}$ is zero up to order 1 (and so 0). Moreover
$$
(Y_{i_{2}}\cdot X_{i_{1}})|_{p}=(\frac{\partial}{\partial\xi_{i_{2}}}\cdot X_{i_{1}})|_{p}
=(\frac{\partial}{\partial\xi_{i_{2}}}\cdot (Y_{i_{1}}\cdot X))|_{p}\pm(\frac{\partial}{\partial\xi_{i_{2}}}\cdot (X\cdot Y_{i_{1}}))|_{p}
=(\frac{\partial}{\partial\xi_{i_{2}}}\cdot (Y_{i_{1}}\cdot X))|_{p}
$$
because $X$ is zero up to order 1. Then
$$
0=(\frac{\partial}{\partial\xi_{i_{2}}}\cdot (Y_{i_{1}}\cdot X))|_{p}=(\frac{\partial}{\partial\xi_{i_{2}}}\cdot (\frac{\partial}{\partial\xi_{i_{1}}}\cdot X))|_{p}
$$
$$
\phantom{ccccccccccccc}=\sum_{r=1}^{m}(\frac{\partial^{2}}{\partial\xi_{i_{2}}\partial\xi_{i_{1}}}f^{r})(p)\frac{\partial}{\partial x^{r}}+\sum_{s=1}^{n}(\frac{\partial^{2}}{\partial\xi_{i_{2}}\partial\xi_{i_{1}}}g^{s})(p)\frac{\partial}{\partial\xi_{s}}=0
$$
where the second equality follows from the fact that $X$ is zero up to order 1.\\
{\bf $(q\geq 3)$} $X_{i_{q}\cdot\cdot i_{1}}=\mathcal{L}_{Y_{i_{q}}}X_{i_{q-1}\cdot\cdot i_{1}}=Y_{i_{q}}\cdot X_{i_{q-1}\cdot\cdot i_{1}}\pm X_{i_{q-1}\cdot\cdot i_{1}}\cdot Y_{i_{q}}$ and evaluating in $p$
$$
0=X_{i_{q}\cdot\cdot i_{1}}|_{p}=(Y_{i_{q}}\cdot X_{i_{q-1}\cdot\cdot i_{1}})|_{p}
$$
because $X_{i_{q-1}\cdot\cdot i_{1}}$ is zero up to order 1 (and so 0). Moreover
$$
(Y_{i_{q}}\cdot X_{i_{q-1}\cdot\cdot i_{1}})|_{p}=(\frac{\partial}{\partial\xi_{i_{q}}}\cdot X_{i_{q-1}\cdot\cdot i_{1}})|_{p}
=(\frac{\partial}{\partial\xi_{i_{q}}}\cdot (Y_{i_{q-1}}\cdot X_{i_{q-2}\cdot\cdot i_{1}}))|_{p}
$$
$$
\phantom{ccccccccccccccccccccccc}\pm(\frac{\partial}{\partial\xi_{i_{q}}}\cdot (X_{i_{q-2}\cdot\cdot i_{1}}\cdot Y_{i_{q}}))|_{p}
=(\frac{\partial}{\partial\xi_{i_{q}}}\cdot (Y_{i_{q-1}}\cdot X_{i_{q-2}\cdot\cdot i_{1}}))|_{p}
$$
because $X_{i_{q-2}\cdot\cdot i_{1}}$ is zero up to order 2 (and so 1) by the induction hypothesis. Using the fact that $X_{i_{k}\cdot\cdot i_{1}}$ is zero up to order $q-k$ (and so $q-k-1$) by the induction hypothesis we get that
$$
0=(Y_{i_{q}}\cdot X_{i_{q-1}\cdot\cdot i_{1}})|_{p}=(\frac{\partial}{\partial\xi_{i_{q}}}\cdot (Y_{i_{q-1}}\cdot (Y_{i_{q-2}}\cdot\cdot\cdot(Y_{i_{1}}\cdot X)\cdot\cdot\cdot)))|_{p}
$$
$$
\phantom{cccccccccccccccccccccccc}=(\frac{\partial}{\partial\xi_{i_{q}}}\cdot (\frac{\partial}{\partial\xi_{i_{q-1}}}\cdot (\frac{\partial}{\partial\xi_{i_{q-2}}}\cdot\cdot\cdot(\frac{\partial}{\partial\xi_{i_{1}}}\cdot X)\cdot\cdot\cdot)))|_{p}
$$
$$
\phantom{ccccccccccccccccccccccccccccc}=\sum_{r=1}^{m}(\frac{\partial^{q}}{\partial\xi_{i_{q}}\cdot\cdot\partial\xi_{i_{1}}}f^{r})(p)\frac{\partial}{\partial x^{r}}+\sum_{s=1}^{n}(\frac{\partial^{q}}{\partial\xi_{i_{q}}\cdot\cdot\partial\xi_{i_{1}}}g^{s})(p)\frac{\partial}{\partial\xi_{s}}
$$
The second to last equality follows from the fact that $X$ is zero up to order $q-1$ by the induction hypothesis. $\square$

\subsubsection{Proof of Theorem \ref{primaparte}}\hfill\newline\\
Every $G$-invariant superfunction $f$ is constant;
$$
\rho_{p}^{*}f=(\Id\otimes \ev_{p})\circ\rho^{*}f=(\Id\otimes \ev_{p})(1\otimes f)=\tilde{f}(p)1_{G}=\rho_{p}^{*}(\tilde{f}(p)1_{M})
$$
implies the assertion. 
\\
Let $X\in\mathcal{T}(M)^{G}$ be a $G$-invariant vector field. For every $g\in G_{0}$, applying $(\ev_{g}\otimes \Id)$ to both sides of $(\Id\otimes X)\circ\rho^{*}=\rho^{*}\circ X$,
we get that
$
X\circ g^{*}=g^{*}\circ X
$ while
$$
X\circ (A|_{e}\otimes \Id)\circ\rho^{*}=(-1)^{|X||A|}(A|_{e}\otimes \Id)\circ(\Id\otimes X)\circ\rho^{*}
$$
and
$$
(A|_{e}\otimes \Id)\circ\rho^{*}\circ X=(A|_{e}\otimes \Id)\circ(\Id\otimes X)\circ\rho^{*}
$$
imply that $\mathcal{L}_{\hat{A}}X=0$ for every $A\in\gg$. We have proved that $\mathcal{T}(M)^{G}\subseteq\mathcal{T}(M)^{(G_{0},\gg)}$. The inclusion $$\ev_{p}(\mathcal{T}(M)^{G})\subseteq \ev_{p}(\mathcal{T}(M)^{(G_{0},\gg)})\subseteq (T_{p}M)^{\phi(G_{p})}$$ follows from evaluating (\ref{invariantfinale}) at $p$
when $g\in (G_{0})_{p}$ and $A\in\gg_{p}$. 
Let $X_{1},X_{2}\in\mathcal{T}(M)^{G}$ be two $G$-invariant vector fields on $M$ with the same value $X_{1}|_{p}=X_{2}|_{p}=v$.
Applying $(\Id\otimes \ev_{p})$ to both sides of $(\Id\otimes X)\circ\rho^{*}=\rho^{*}\circ X$, we get that
$$
(\Id\otimes v)\circ\rho^{*}=\rho_{p}^{*}\circ X_{1}=\rho_{p}^{*}\circ X_{2}\quad,
$$ 
\textit{i.e.} $X_{1}=X_{2}$. This proves injectivity of the correspondence at the level of vector fields. We prove surjectivity.
For every $v\in (T_{p}M)^{\phi(G_{p})}$ equation (\ref{tenosi})
defines a $G$-invariant vector field through the identification $\rho_{p}:G/G_{p}\rightarrow M$.
First, equation (\ref{homogeneous}) together with
$$
R_{h}^{*}\circ(\Id\otimes v)\circ\rho^{*}=(\Id\otimes v)\circ (R_{h}^{*}\otimes \Id)\circ\rho^{*}=(\Id\otimes v)\circ (\Id\otimes h^{*})\circ\rho^{*}=
(\Id\otimes v)\circ\rho^{*}\, ,
$$
$$
(-1)^{|B||v|}B\circ(\Id\otimes v)\circ\rho^{*}=(\Id\otimes v)\circ(B\otimes \Id)\circ\rho^{*}=(\Id\otimes v)\circ(\Id\otimes \hat{B})\circ\rho^{*}=0\, ,
$$
where $h\in (G_{0})_{p}$ and $B\in\gg_{p}$, implies that
$\Im((\Id\otimes v)\circ\rho^{*})\subseteq\mathcal{A}_{G/G_{p}}$.
Equation
$$
(\Id\otimes v)\circ\rho^{*}(fg)=(\Id\otimes v)\rho^{*}(f)\cdot\rho_{p}^{*}(g)+(-1)^{|f||v|}\rho_{p}^{*}(f)\cdot(\Id\otimes v)\rho^{*}(g)
$$
implies that (\ref{tenosi}) is a derivation. It is $G$-invariant, for $$\rho^{*}\circ(\rho_{p}^{*})^{-1}\circ((\Id\otimes v)\circ\rho^{*})$$ equals
$$
(\Id\otimes(\rho_{p}^{*})^{-1})\circ(\Id\otimes (\Id\otimes v)\circ\rho^{*})\circ\rho^{*}=(\Id\otimes(\rho_{p}^{*})^{-1})\circ(\Id\otimes \Id\otimes v)\circ(\Id\otimes \rho^{*})\circ\rho^{*}
$$
$$
=(\Id\otimes(\rho_{p}^{*})^{-1})\circ(\Id\otimes \Id\otimes v)\circ(m^{*}\otimes \Id)\circ\rho^{*}=(\Id\otimes(\rho_{p}^{*})^{-1})\circ \mu^{*}\circ((\Id\otimes v)\circ\rho^{*})\quad.
$$
The inclusion $\mathcal{T}(M)^{(G_{0},\gg)}\subseteq \mathcal{T}(M)^{G}$ is proved as in Lemma \ref{allora3}. 
\\ 
Let $\omega\in\mathcal{T}^{*}(M)^{G}$ be a $G$-invariant 1-form. The inclusion $$\mathcal{T}^{*}(M)^{G}\subseteq\mathcal{T}^{*}(M)^{(G_{0},\gg)}$$ follows from applying $(\ev_{g}\otimes \Id)$ and $(A|_{e}\otimes \Id)$ to both sides of $\omega(Y)=(\rho^{*}\omega)(\Id\otimes Y)$. Thus $$\ev_{p}(\mathcal{T}^{*}(M)^{G})\subseteq \ev_{p}(\mathcal{T}^{*}(M)^{(G_{0},\gg)})\subseteq (T_{p}^{*}M)^{\phi(G_{p})}\quad.$$
Let $\omega_{1},\omega_{2}\in\mathcal{T}^{*}(M)^{G}$ be two $G$-invariant 1-forms with the same value $\omega_{1}|_{p}=\omega_{2}|_{p}=\omega_{p}$.
Applying $(\Id\otimes \ev_{p})$ to both sides of
$$
1\otimes\omega_{i}(\hat{A})=(\rho^{*}\omega_{i})(\Id\otimes \hat{A})
$$
we get that
$$
\omega_{p}(\hat{A}|_{p})\cdot 1_{G}=(\Id\otimes \ev_{p})(\rho^{*}\omega_{i})(\Id\otimes \hat{A})=(\rho_{p}^{*}\omega_{i})(A)\quad,
$$
\textit{i.e.} $\omega_{1}=\omega_{2}$. Denote $G_{p}$ by $H$.
A $1$-form $\omega\in\mathcal{T}^{*}(G)$ is {\it projectable} if
$$
\label{stramba1}
\textbf{$R_{H_{0}}$-Invariance\,\,\,\,}\phantom{ccccccccccccccccccccccccccccccccccccc}\omega=R_{h}^{*}\omega\phantom{\qquad Y\in\mathcal{T}(G)}
$$
$$
\label{stramba2}
\textbf{$\gh$-Invariance}\phantom{ccccccccccccccccccc} B(\omega(Y))=(-1)^{|\omega||B|}\omega(\mathcal{L}_{B}Y)\qquad (Y\in\mathcal{T}(G))
$$
$$
\label{stramba3}
\textbf{Horizontality\,\,\,\,\,\,}\phantom{ccccccccccccccccccccccccccccccccccc}  \omega|_{\mathcal{A}_{G}\otimes\gh}=0\phantom{\qquad Y\in\mathcal{T}(G)}
$$
for every $h\in H_{0}$,  $B\in\gh$. The set of all projectable 1-forms is denoted by $\mathcal{T}^{*}(G)_{hor}^{R_{H}}$ and its intersection with the set of $G$-invariant 1-forms by $$\mathcal{T}^{*}(G)_{hor}^{L_{G}\times R_{H}}:=\mathcal{T}^{*}(G)_{hor}^{R_{H}}\cap\mathcal{T}^{*}(G)^{G}\quad.$$ 
The pull-back
$$\mathcal{T}^{*}(G/H)\owns\nu\overset{\pi^{*}}{\longrightarrow} \pi^{*}\nu\in\mathcal{T}^{*}(G)_{hor}^{R_{H}}$$
is a bijection such that
$\pi^{*}(\mathcal{T}^{*}(G/H)^{G})=\mathcal{T}^{*}(G)_{hor}^{L_{G}\times R_{H}}$.
A 1-form $\omega\in\mathcal{T}^{*}(G)_{hor}^{R_{H}}$ uniquely defines a 1-form $\nu\in\mathcal{T}^{*}(G/H)$ by
$$
\nu(X):=\omega(Y)
$$
where $X\in\mathcal{T}_{G/H}(U)$, $Y\in\mathcal{T}_{G}(\pi_{0}^{-1}U)$ are $\pi$-related vector fields. These assertions depend on the existence of adapted coordinates. 
For every $H$-invariant covector $\omega_{o}\in (\gg/\gh)^{*}$, the formula
$$
\omega(f\otimes A):=(-1)^{|f||\omega_{o}|}f\omega_{o}(A)\in \mathcal{A}(G)\qquad\quad f\otimes A\in \mathcal{A}(G)\otimes \gg
$$
defines a 1-form $\omega\in\mathcal{T}^{*}(G)_{hor}^{L_{G}\times R_{H}}$ and so a $G$-invariant 1-form on $G/H\cong M$. \\
These proofs can be generalized naturally to tensor fields of type $(r,s)$. $\square$
\medskip
\\
\font\smallit = cmti8
{\smallit Acknowledgement}\/. {\small The author is grateful to the University of Florence and I.N.d.A.M. for financial support during his stay in Edinburgh and to the University of Edinburgh for hospitality. This paper is essentially one chapter of the Ph.D. thesis the author wrote under the supervision of D.\ V.\ Alekseevsky. The author would like to thank him and A.\ Spiro for their constant encouraging support. He also would like to thank his friends A.\ D'Avanzo and D.\ Sepe for their careful reading of the paper.}

\bigskip
\bigskip
\font\smallsmc = cmcsc8
\font\smalltt = cmtt8
\font\smallit = cmti8
\hbox{\parindent=0pt\parskip=0pt
\vbox{\baselineskip 9.5 pt \hsize=3.1truein
\obeylines
{\smallsmc
Andrea Santi
Dept. of Mathematics
University of Florence
Viale Morgagni 67/a  
Florence 50134 
Italy
Phone +390554237111
Fax +390554222695
}\medskip
{\smallit E-mail}\/: {\smalltt santi@math.unifi.it}}}

\end{document}